%% file: thesis.tex
\documentclass[12pt,offcenter,edeposit]{uiucthesis07}
\usepackage[top=1.0in, bottom=1.0in, left=1.5in, right=1.0in]{geometry} 
\usepackage{amsmath,amssymb}
\usepackage{graphicx}
\usepackage{setspace}
\usepackage{enumerate}
\newtheorem{theorem}{Theorem}[section]
\newtheorem{remark}[theorem]{Remark}

\begin{document}

\title{A STABILIZED FINITE ELEMENT FORMULATION OF NON-SMOOTH CONTACT}
\author{Ghadir Haikal}
\department{Civil Engineering}
\schools{B.C.E, Tishreen University, 1998\\
         M.S., University of Illinois at Urbana-Champaign, 2004}
\phdthesis
\advisor{Keith D. Hjelmstad}
\committee{Professor Keith D. Hjelmstad, Chair\\ Professor Robert H. Dodds, Jr\\ Professor Daniel A. Tortorelli\\ Associate Professor Arif Masud\\ Assistant Professor Ilinca Stanciulescu Panea}
\degreeyear{2009}
\maketitle

\frontmatter

\begin{abstract} 

The computational modeling of many engineering problems using the Finite Element
method involves the modeling of two or more bodies that meet through an
interface. The interface can be physical, as in multi-physics and contact problems,
or purely numerical, as in the coupling of non-conforming meshes. The
most critical part of the modeling process is to ensure geometric compatibility
and a complete transfer of surface tractions between the different components
at the connecting interfaces. Contact problems are a special family of interaction
problems where the bodies on either side of the interface may separate freely or
connect with each other, depending on the direction of motion. This type of behavior
can be observed in complex civil, mechanical, bio-mechanical or aerospace
structural components, and, on a smaller scale, in the interaction of different
constituents in heterogeneous and composite materials and in the opening
and closing of cracks in fracture mechanics.

Popular contact modeling techniques rely on geometric projections to detect and
resolve overlapping or mass interpenetration between two or more contacting
bodies. Such approaches have been shown to have two major drawbacks: they
are not suitable for contact at highly nonlinear surfaces and sharp corners
where smooth normal projections are not feasible, and they fail to guarantee a
complete and accurate transfer of pressure across the interface. This dissertation
presents a novel formulation for the modeling of contact problems that possesses
the ability to resolve complicated contact scenarios effectively, while being simpler
to implement and more widely applicable than currently available methods. We show that the formulation boils down to a node-to-surface gap function that works effectively for non-smooth contact. The numerical implementation using the midpoint rule shows the need to guarantee the conservation
of the total energy during impact, for which a Lagrange multiplier method is used. We propose a local enrichment of the interface and a simple stabilization
procedure based on the discontinuous Galerkin method to guarantee an accurate transfer
of the pressure field. The result is a robust interface
formulation for contact problems and the coupling of non-conforming meshes.

\end{abstract}

\begin{dedication}
To the spirit of my father, for being my biggest fan and my inspiration, \\
and to my mother, for being my rock and my best friend.
\end{dedication}

\chapter*{Acknowledgments}
First, I want to thank my advisor, Professor Keith D. Hjelmstad, for his direction and guidance throughout my graduate education at Illinois. I have learned a lot from him, both as an academic and as a person, and I could not have asked for a better supervisor or mentor. His passion for mechanics and education is very inspiring and I greatly appreciate his patience and availability despite his busy schedule. His feedback was instrumental in the completion of this project and the success of my search for an academic position.

I would also like to thank the members of my dissertation committee, Professor Robert H. Dodds Jr., Professor Daniel A. Tortorelli, Professor Arif Masud, and Professor Ilinca Stanciulescu, for their involvement in this project and their continuing support. I especially appreciate their patience in the final stages of the preparation of this thesis. Professor Dodds' class on computer methods has shaped my approach to computer programming and numerical methods. I have thoroughly enjoyed the discussions with Professor Tortorelli during our research group meetings. Thanks to Professor Masud for his great classes on finite element methods and to Professor Ilinca Stanciulescu for her valuable career advice and encouragement.

I am grateful to the members of our research group, Drs. Alireza Namazifard, Arun Prakash, Kalyana Babu (Kabab) Nakshatrala, Kristine Cochran, Daniel Turner, and Daqing Xu for enriching my academic experience at Illinois and for their close friendship. I will always cherish the great times I have spent with Dan and Kristine and their families, and will never forget the heated discussions I have had with Kabab and Arun about mechanics, life, and academia at Moonstruck. 

Special thanks to my family away from home: Liz and Christopher Smiley, Chris, Bill and Caroline Denniston, and to my close friends: George Alhaj, Maha El Choubassi, Ava Zeineddin, James Sobotka and Jason Patrick for sharing the good, the bad, and the not-so-occasional cup of coffee or lunch. Thanks to the students in CEE471 for making my life harder, yet much more enjoyable, as a teaching assistant at Illinois.

I am forever indebted to my parents, Dr. Ghassan and Souraya, and to my siblings, Zein, Razan and Ghassan Jr., for their unwavering love. This journey has been met by many difficulties, and although the academic obstacles were stimulating and inspiring, the personal challenges would have been paralyzing without their support. 


\tableofcontents
\listoffigures





\mainmatter

\include{introduction}

\include{contact}

\include{motivation}
\include{DG}

\include{enrichment}

\chapter{Conclusions and future directions}

Interface problems in general, and contact problems in particular, are a great challenge in computational mechanics. The difficulty in modeling interface problems lies in representing the conditions of equilibrium and geometric compatibility along the interface faithfully. Although these conditions are straightforward in the continuum setting, the discretization with finite elements often restricts the ability of the numerical model to enforce both conditions simultaneously. Enforcing geometric compatibility becomes especially challenging for contact on non-smooth surfaces where surface normals are not uniquely defined. The simple approach of imposing the geometric compatibility constraints to the interface degrees of freedom, located at the nodes, does not yield a complete projection of the interface forces and therefore exhibits numerical instabilities. Alternative dual coupling approaches that employ a field of Lagrange multipliers to enforce weak geometric compatibility are sensitive to the choice of the Lagrange multiplier discretization, which is based on the user-defined \textit{master} side of the interface. As a result, such methods are sensitive to the choice of the master surface. In the particular case of unilateral contact, weak continuity can lead to a physically inadmissible configuration with large mass overlap, especially if the contacting surfaces are not smooth. More importantly, the geometry of the non-smooth surfaces has a great influence on the response and therefore needs to be captured accurately. 

In this dissertation, we have proposed a novel interface formulation that is suitable for contact problems and the coupling of non-conforming meshes. The primary feature of the proposed formulation is that it enables the enforcement of a set of discrete compatibility constraints at the interface nodes. Such feature is essential to the contact problem since it helps prevent mass interpenetration or overlap, especially at non-smooth locations. We presented a new formulation of the contact constraints that is suitable for contact at non-smooth locations such as corners. The contact constraint is based on the calculation of an oriented volume and does not require a unique definition of the surface normal. We have shown that the contact constraint formulation applies to elements of any order, both in 2D and in 3D.

To enable an unbiased two-pass approach where the geometric compatibility conditions can be enforced at both sides of the interface, we proposed a local enrichment of the interface elements that would transform the node-to-surface constraint function to a node-to-node one. We have shown that the enrichment prevents the  \textit{locking} or artificial stiffening of the interface as a result of imposing the compatibility constraints on the discretized surfaces. A nice feature of the enrichment is that it is local to the interface elements and does not affect the remainder of the mesh. For the coupling of non-conforming meshes, the additional degrees of freedom are handled directly by the assembly procedure without any increase in the size of the problem.

Next, we examined the issue of traction transfer across the interface. We have shown that for the transfer of forces between two surfaces to be complete, the variational displacement field has to be continuous, at least in a weak sense. In the standard Bubnov-Galerkin finite element discretization, continuity of the variational displacement field is closely coupled with the continuity of the real field. For a conforming mesh, this condition is satisfied by design. For the more general non-conforming interface, we have shown that it is possible to achieve point-wise continuity using a set of Multi-Point-Constraints (or gap functions for contact) when the meshes are discretized with bilinear elements and a two-pass approach is used. This scenario, however, exhibits severe surface locking.

The proposed surface enrichment enforces displacement continuity at the nodes only and does not generally lead to weak continuity along the interface. To guarantee a complete transfer of forces along the non-conforming part of the interface, we implement a stabilization procedure based on the discontinuous Galerkin method. The stabilization terms are in fact derived from the continuum statement of the coupled problem and aim to enforce weak equilibrium of interface tractions. These terms vanish in a conforming mesh and the proposed DG formulation includes the standard Galerkin method as a subset.

In sum, the main contribution of this work is a robust two-pass unbiased interface model that accommodates non-smooth contact conditions, guarantees a complete transfer of interface forces and strongly enforces the non-penetration constraint at all interface nodes without inducing surface locking. This formulation has a clear advantage over the biased single-pass methods that enforce the interpenetration condition only in a weak sense, and the two-pass methods that require a user-defined mesh-dependent parameter for stability. The proposed interface formulation yields much improved estimates of the interface stresses and shows great promise and potential for application to more complex coupling problems. A few possible future developments for this method include the following
\begin{itemize}

\item {\textbf{3D interfaces:} As discussed in Section \ref{sec_DG:implementation} the extension to 3D only involves implementing a robust interface detection and integration algorithm and does not involve any changes in formulation.}
\item {\textbf{Dynamic coupling:} We have shown in Chapter 2 that the implementation of the contact constraints in a dynamic setting using the mid-point rule requires enforcing the conservation of energy as an additional constraint. This is primarily due to the fact that the mid-point scheme assumes a smooth velocity profile that is not generally true of contact/impact scenarios. We have developed the interface formulation for the contact problem in a quasi-static setting to simplify the presentation of the method. As mentioned previously, extending this formulation to a dynamic framework where the discretization is non-conforming in space but continuous in time is fairly straightforward. The fact that the stabilization terms are derived from the DG method, however, is a very important aspect that we expect to be key to the dynamic contact problem. DG methods that are discontinuous in time are especially suited to impact scenarios. Therefore we plan to explore an interface stabilization procedure that accommodates discontinuity in both time and space along the contact interface.}
\item{\textbf{Frictional contact:} The extension to the frictional case involves adding the tangential frictional forces to the interface equilibrium term.} The proposed interface stabilization shows great promise for this particular application, given its better estimates of the interface stress fields, especially the shear stress.
\item{\textbf{The coupling of multi-physical domains:} In multi-physics problems such as fluid-structure and soil-structure interaction each physical domain is typically discretized using an appropriate numerical method such as finite volume methods for fluids, boundary elements methods or meshless methods for soil, and finite element methods for the structural domain. Thus, the modeling of such problems would require developing an efficient interface model that could couple different spatial discretization schemes for multi-physical domains, including effects such as friction and heat transfer.}
\end{itemize}
\appendix

\include{appendix}

\backmatter

\nocite{dxu08,kulkarniDDMSE04,kulkarniIJNME07,haberCMAME84,huertaIJNME00,duarteIJNME07,belytschkoCM07,liuCMAME04,chenIJNME00}
\bibliographystyle{plain}
\bibliography{thesisreferences}

\chapter{Author's Biography}

Ghadir Haikal was born on October 24, 1975, in Damascus, Syria. She holds a Bachelor's degree (B.C.E.) in Civil Engineering and a Masters degree (M.S.) in Structural Engineering from Tishreen University in Lattakia, Syria, in 1998 and 2002, respectively. She graduated from the University of Illinois with a Masters degree (M.S.) in Civil Engineering in December, 2004. 

\end{document}

%% file: introduction.tex
\chapter {Introduction}

Unilateral contact constraints are typically employed in finite element analysis of structures to prevent overlapping or mass interpenetration when solid bodies collide. The early works on contact relied on geometric distance to characterize the potential for contact between two solid bodies. The contact constraint function in such a case is the oriented distance or gap between a candidate node and its normal projection on the contacting body. The impenetrability constraint is enforced between the so-called \textit{slave} node and \textit{master} surface using a discrete node-to-surface gap function evaluated at the slave node.  A set of equivalent nodal forces is computed at the slave nodes that represent the pressure acting along the contact surface. 

Discrete gap functions have been extensively used as contact constraints due to their relative simplicity and applicability to all types of finite element meshes. However, the discrete node-to-surface gap formulation does not pass the contact patch test. The patch test is designed to verify that a given contact formulation is capable of representing a state of constant pressure, thus ensuring the completeness of the pressure field. The discrete gap function satisfies this condition only when the two contact surfaces are discretized with linear elements at the contact events happen at the nodes. For quadratic and higher order elements, the contacting nodes have to be of the same type (edge node with edge node, ..., etc.).
In the general case of non-matching elements, the transfer of the pressure field is not complete.

Papadopoulos and Taylor \cite{papadopoulosCMAME92} introduced an averaged node-to-surface gap function, in which the impenetrability constraint is enforced in an average sense along the whole contact surface and the contact pressure is integrated along the contact slide line via Simpson's rule. Jones and Papadopoulos \cite{jonesIJNME01} later proposed a 3D contact formulation that employs an isoparametric interpolation of the contact pressure on both sides of the contact surface. Equilibrium is then enforced in a weak sense between the two contacting surfaces, and the interpenetration condition is enforced between a set of slave points (typically integration points) and the master surface. This family of formulations passes the contact patch test by design.

The aforementioned methods suffer from major drawbacks. Firstly, the robustness of the solution procedure can be affected by the choice of the master/slave pairing. Better results have been observed when the coarser of two contacting surfaces is designated as the master surface. This bias can be eliminated via a two-pass procedure where each of the two surfaces serves as a \textit{master} to the nodes of the other. However, two-pass methods have been shown to fail the Ladyzhenskaya-Babu$\breve{s}$ka-Brezzi (LBB) condition and may therefore exhibit \textit{surface locking} \cite{papadopoulosCMAME92}. This phenomenon is an artifact of the finite element discretization of the contact surfaces and occurs when enforcing the impenetrability condition at multiple locations leads to an artificial stiffening of the contact surface. Consequently, convergence of the solution cannot be obtained in the limit of mesh refinement. Locking can be a major handicap when the contact surfaces are not smooth, either due to the actual geometry of the problem or as a result of the finite element discretization. Secondly, the formulation of a gap function requires a unique definition of the surface normal at the location of contact. Therefore, traditional gap function models are only applicable to contact on smooth surfaces and are usually referred to as \textit{smooth} contact constraints.  

More recently, non-smooth contact models have been developed \cite{kaneCMAME99}. The contact between two bodies is said to be non-smooth when it can occur along the smooth boundaries of the bodies as well as at non-smooth locations such as corners. This possibility disables the treatment of the problem using smooth analysis tools such as gap functions and surface normals and requires an appropriate definition of the contact potential regardless of its location. Modeling of non-smooth contact is essential in many engineering applications such as granular flow and fragmentation, and non-smooth contact formulations that are suitable to such applications have been the focus of recent research \cite{kaneCMAME99, pusoCMAME04}.

The non-smooth contact formulations go back at least to Kane et al. \cite{kaneCMAME99} who used the intersection area (or volume in the 3D case) between the contacting bodies as a contact constraint function. The main drawback of this approach is that it is restricted to geometrically linear triangular elements such as 3-node triangles in 2D and 4-node tetrahedra in 3D. Moreover, the resulting constraints are highly nonlinear and the implementation requires special care in defining the orientation of the elements in space, since the sign of the function is determined by the relative locations of the nodes of the contacting elements. Belytschko et al. \cite{belytschkoIJNME02} proposed computing the gap function with respect to a smoothed surface that represents a least-squares fit to the original, non-smooth one.

Adopting the concept of pressure interpolation introduced by Papadopoulos and Taylor, Puso and Laursen \cite{pusoCMAME04} developed a segment-to-segment formulation, called the mortar method. In this single-pass method, the gap function is averaged along the contacting
segments and the pressure at the slave contact points is interpolated in terms of the nodal pressures of the master surface. This approach uses an \textit{averaged nodal normal} to address the non-uniqueness of the normal to the contacting segments at non-smooth locations. The mortar approach can be applied to non-smooth contact scenarios and proved to overcome the over-constraint associated with discrete node-to-surface gap functions. However, since it is based on a weighted average of the contact gap, this method can leave some unresolved mass interpenetration at non-smooth contact locations. Yang et al. \cite{yangIJNME05} extended the mortar method to large deformations. This approach has also been applied to curved surfaces by Flemish et al. \cite{flemischIJNME05} and to quadratic elements by Puso et al. \cite{pusoCMAME08}. 

A mortar-based two-pass formulation has recently been developed by Solberg et al. \cite{solbergCMAME07}. Unlike the single-pass mortar method, this approach strongly enforces the contact constraints at a number of select points while continuity of pressure is satisfied in a weak sense along the contact surface. As a result, a penalty-based stabilization term needs to be imposed to minimize the pressure jump across the contact interface. It is suggested that on each surface, nodes [be] a priori identified as active or inactive, via a \textit{binary patterning scheme} where the interface nodes are alternatively designated as active contact locations. This ad-hoc approach is not guaranteed to work for arbitrary meshes.

In this dissertation, we present a new formulation of nonsmooth contact constraints and their implementation in a dynamic nonlinear finite element framework. Based on the calculation of an oriented volume, the suggested contact constraint formulation allows for a simple and unified treatment of all potential contact scenarios in the presence of large deformations. The elements of particular interest to this study are quadrilateral and hexahedral elements, although the results can be easily extended to triangular elements. The proposed approach is equally applicable to bilinear and higher-order elements, for which no nonsmooth contact formulation has been developed to date, both in 2D and 3D. We show that this formulation is a modified version of the discrete node-to-surface gap function that retains its advantages and avoids its main shortcoming.

Furthermore, we propose a stabilized interface formulation, based on the non-smooth node-to-surface gap function, that would cure locking due to over-constraint. This stabilization is in the form of a local enrichment of the contact surface that would transform the node-to-surface gap function to a node-to-node one that passes the patch test for all types of configurations, thereby eliminating the need for a master-slave definition while still satisfying the LBB condition. The importance of this work stems from the fact that the procedures suggested in the literature to address the issue of surface locking in two-pass-methods are mostly ad-hoc. Unlike the single-pass mortar method, the robustness of the solution in the proposed approach is not affected by the choice of the master surface. The contact constraints are enforced strongly at the nodes and therefore no regularization is needed for either the contact gap or pressure across the interface. The computational cost is greatly reduced since the contact effects can be treated locally and no additional fields are introduced. The scope of this work is limited to elastic frictionless contact. 

Contact problems are a special case of the larger class of interface and coupling problems, with the distinction that the bodies on either side of the interface are free to move apart or come in contact with each other. While this unilateral behavior adds to the complexity of the mathematical model, some of the numerical issues encountered in contact problems, such as patch-test performance and surface locking, can also be found in similar interface problems such as the coupling of non-conforming meshes. Therefore, the solutions to these issues are equally applicable to both problems as well. In fact, some of the current contact modeling techniques, namely the mortar \cite{bernardiNPDEA92} and Nitsche methods \cite{nitsche71}, have originated in the domain decomposition literature. Therefore, to simplify the presentation of our proposed interface formulation, we apply it first to the domain decomposition problem and the tying of non-conforming meshes, before introducing the contact problem where other more involved effects such as sliding and large deformations come in play.

The outline of the dissertation is as follows. In Chapter 2, we outline the mathematical formulation of the equations of motion (including the finite element discretization and numerical solution procedures) leading to the statement of the nonlinear constrained optimization problem. We then describe the suggested approach for the formulation of the contact constraints. Chapter 4 discusses the background and motivation for the stabilized interface formulation. In Chapter 5, we present the proposed interface formulation in the context of coupling non-conforming meshes. Chapter 6 extends the formulation to contact interfaces. Finally, in Chapter 7, we present our conclusions and discuss future applications.

%% file: contact.tex
\chapter{A finite element formulation of non-smooth contact} 
%
%
%
%
%
%

\section{Introduction} 

This chapter concerns the finite element modeling of contact between solid bodies, with a special emphasis on the treatment of nonsmooth conditions.  We propose a new formulation of the contact constraints that does not require the explicit definition of a surface normal and therefore works effectively for non-smooth contact scenarios. We restrict our attention to quadrilateral and hexahedral elements in frictionless contact, although the results can readily extended to triangular elements and to the general framework of frictional contact. 
\footnote{The early stages of this work were initiated by the author in \textit{A new approach for the finite element formulation of contact based on intersecting volumes for quadrilateral and hexahedral elements}, MS thesis, University of Illinois, 2004. The contents of this chapter are part of a published article \cite{haikalCMAME07}.}

Although contact can be treated statically or quasi-statically, dynamic analysis enables a more general solution where the motion of the interacting bodies during and after collision can be simulated. Moreover, high impact forces often arise due to contact and these forces are most adequately accounted for within a dynamic framework. To avoid the numerical instabilities that are typically encountered in the numerical simulation of contact problems, we employ an implicit energy-preserving time stepping scheme. This scheme is a variant of the mid-point rule in which energy conservation is enforced via a Lagrange multiplier method. 

The outline of the chapter is as follows. Section \ref{sec_con:MathForm} introduces the mathematical formulation of the equations of motion (the finite element discretization and numerical solution procedures), leading to the statement of the nonlinear constrained optimization problem. In Section \ref{sec_con:OrVol}, we present the suggested approach for the formulation of the contact constraints and describe the implementation procedure. Section \ref{sec_con:ImplSol} outlines the solution procedure and the computational algorithm. The numerical examples are presented and discussed in Section \ref{sec_con:NumEx}.
%
%
%
%
\section{Mathematical formulation} \label{sec_con:MathForm} 
%
%
%
\subsection{Dynamic equations of motion} \label{sec_con:DynEq}
	\begin{figure}
	\centering
	\includegraphics[clip]{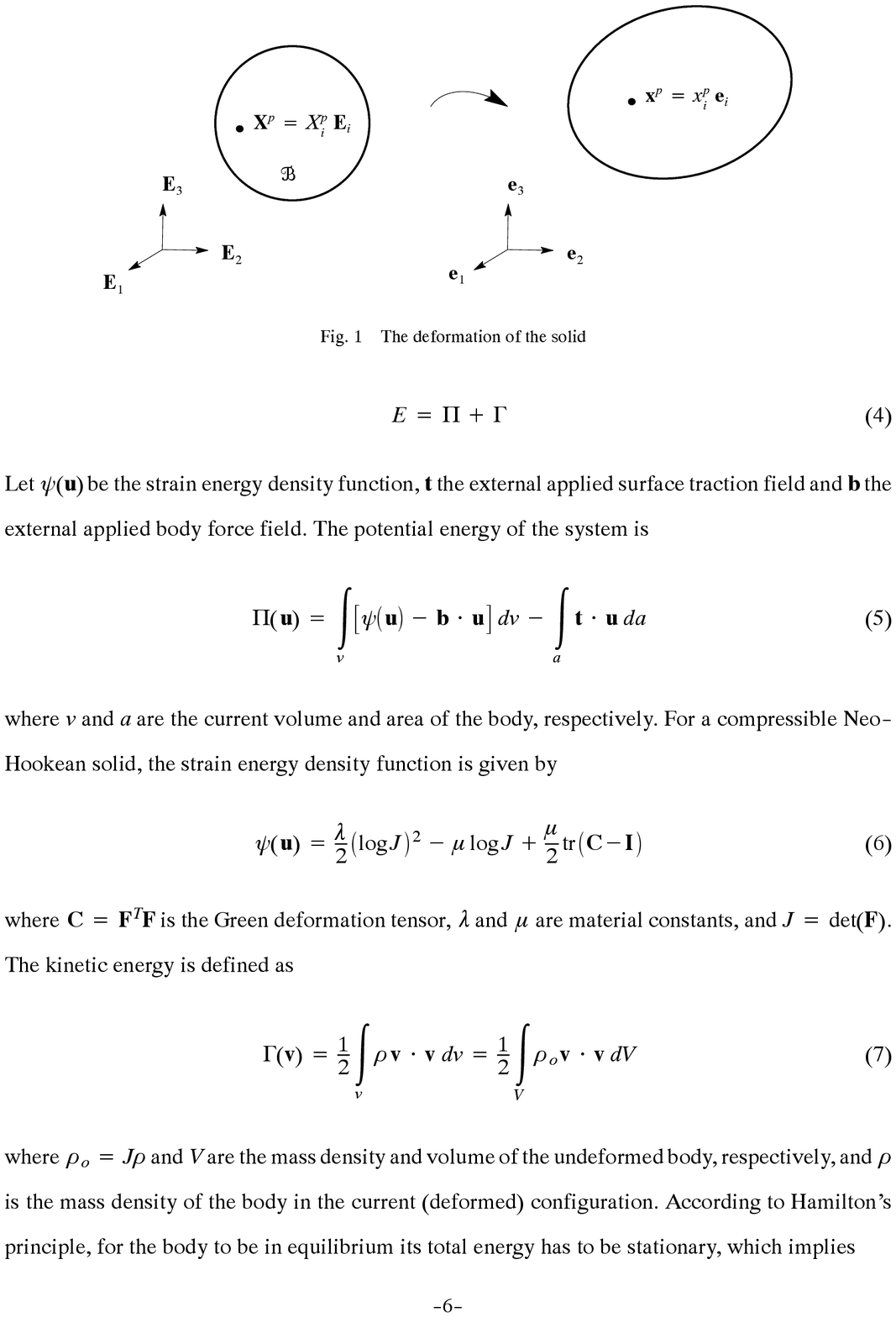}
	\caption{The deformation of the solid  \label{fig:con_potato}}
	\end{figure}
Consider the solid body $\mathcal{B}$ shown in Figure \ref{fig:con_potato}. The initial configuration of the body is defined by the material coordinates of its points $\mathbf{X}^{p} = X_{i}^{p}\mathbf{E}_{i}$ where $p\in\mathcal{B}$ and $\mathbf{E}_{i}$ are the base vectors. The configuration changes as the body moves and deforms. The current configuration is designated by the set of spatial coordinates $\mathbf{x}^{p} = x_{i}^{p}\mathbf{e}_{i}$. In the subsequent derivations, lowercase characters will be used for quantities in the deformed (spatial) configuration of the body whereas uppercase characters will denote quantities in the undeformed (material) configuration. The Einstein summation convention applies. The displacement field of the body is
\begin{equation}
\mathbf{u} = \mathbf{x} - \mathbf{X}.
\end{equation}
The rate of change of this field with time is the velocity $\mathbf{v}={d\mathbf{u}}/dt\equiv\dot{\mathbf{u}}$, and the rate of change of the velocity is the acceleration $\mathbf{a}={d\mathbf{v}}/dt={d^{2}\mathbf u}/dt^{2}\equiv\ddot{\mathbf{u}}$. The deformation of the body is characterized by the deformation gradient $\mathbf{F}$, which represents the rate of change of the position of the body with respect to its material coordinates $\mathbf{X}$ 
\begin{equation}
\mathbf{F} = \nabla_{\mathbf{X}}\mathbf{x}=\mathbf{I}+\nabla_{\mathbf{X}}\mathbf{u},
\end{equation}
where $\nabla_{\mathbf{X}}\mathbf{u}={\partial u_i}/{\partial X_j}\,\mathbf{e}_i\otimes\mathbf{E}_j$. Let $\Gamma$ be the kinetic energy of the system and $\Pi$ be its potential energy,
\begin{align}  
\Gamma = \int_v \frac{\rho}{2}\mathbf{v}\cdot\mathbf{v}\,dv,  \ \ \ \ \
\Pi = \int_v \left[\psi\left(\mathbf{u}\right)-\mathbf{b}\cdot\mathbf{u}\right]\,dv  - \int_a\mathbf{h}\cdot\mathbf{u}\,da, 
\end{align}
where $\psi(\mathbf{u})$ is the strain energy density function, $\mathbf{b}$ is the body force vector and $\mathbf{h}$ is the traction vector acting on the surface of the body. Let $E$ be the total enery of the system,
\begin{align}  \label{eq:E}
E = \Gamma + \Pi &=  \int_v \left[\frac{\rho}{2}\,\mathbf{v}\cdot\mathbf{v}+\psi\left(\mathbf{u}\right)-\mathbf{b}\cdot\mathbf{u}\right]dv  - \int_a\mathbf{h}\cdot\mathbf{u}\,da.
\end{align}
According to Hamilton's principle, for the body to be in equilibrium its total energy has to be stationary, which implies
\begin{align}  \label{eq_con:Edot}
\dot{E} &= \frac{dE}{dt} = \int_V \left[\rho_0\ddot{\mathbf{u}}\cdot\dot{\mathbf{u}}+\mathbf{P}\cdot\dot{\mathbf{F}}-\mathbf{b}_0\cdot\dot{\mathbf{u}}\right]\,dV - \int_A\mathbf{h}_0\cdot\dot{\mathbf{u}}\,dA = 0,
\end{align}
where $\mathbf{b}_0 = J\mathbf{b}$  is the body force vector acting on the body in its undeformed configuration, $\mathbf{h}_0=\mathbf{h}\,da/dA$ is the traction vector acting on the surface of the body in its undeformed configuration and $\mathbf{P}$ is the first Piola-Kirchhoff stress tensor.

Equation \eqref{eq_con:Edot} gives the necessary conditions for equilibrium. These equations are discretized spatially via an isoparametric finite element formulation. Let the body be subdivided into a set of finite elements $V = \cup_{e}V_{e}$. The displacement and position vectors of any point within an element $e$ are expressed in terms of the local nodal quantities as
\begin{equation}
\mathbf{x}^h\equiv N^\alpha\mathbf{x}^\alpha , \ \ \ \  \mathbf{X}^h\equiv N^\alpha\mathbf{X}^\alpha , \ \ \ \  \mathbf{u}^h\equiv N^\alpha\mathbf{d}^\alpha,
\end{equation}
where the vectors $\mathbf{d}^\alpha$, $\mathbf{X}^\alpha$, and $\mathbf{x}^\alpha$ represent the displacements, material coordinates and spatial coordinates of node $\alpha$ in element $e$. A repeated Greek symbol implies summation from $1$ to the number of nodes in an element $n$. In these equations, the shape functions $N^\alpha(\mathbf{\zeta})$ define the mapping of the actual element to a reference element in the coordinates $\mathbf{\zeta}$ (with $i=1,\cdots,$ number of spatial dimensions). Substituting the interpolation into equation \eqref{eq_con:Edot} at the element level yields
\begin{equation} \label{eq_con:dadotFa}
\dot{\mathbf{d}}^\alpha\cdot\left[\mathbf{V}^\alpha+\mathbf{T}^\alpha-\mathbf{F}^\alpha\right]=0,
\end{equation}
where
\begin{equation}
\mathbf{T}^\alpha = \int_{V_e}\mathbf{P}\nabla_{\mathbf{X}}N^\alpha\,dV
\end{equation}
is the element internal force vector and
\begin{equation}
\mathbf{V}^\alpha = \int_{V_e}N^\alpha\rho_0\ddot{\mathbf{u}}\,dV = \left[\int_{V_e}N^\alpha\rho_0N^\beta\, dV\right]\ddot{\mathbf{d}}^\beta = \mathbf{m}^{\alpha\beta}\ddot{\mathbf{d}}^\beta
\end{equation}
is the element nodal inertia force vector, in which $\mathbf{m}$ is the (symmetric) consistent mass matrix of the element. A lumped mass matrix can be obtained by adding the off-diagonal terms on each row to the diagonal term. Lastly,
\begin{equation}
\mathbf{F}^\alpha = \int_{V_e}N^\alpha\mathbf{b}_0\,dV+\int_{A_e}N^\alpha\mathbf{h}_0\,dA
\end{equation}
is the equivalent external element nodal force vector.

Let $\mathbf{d}=\left[\mathbf{d}^1,\mathbf{d}^2,\cdots,\mathbf{d}^n_{total}\right]\in\mathbb{R}^{1+3n_{total}}$ be the global displacement vector in the body, where $n_{total}$ is the total number of nodes in the finite element mesh. We define the $(3\times3n_{total})$ boolean matrix $\mathbf{D}_\alpha$ that extracts the displacements of each node $\alpha$ in element $e$ from the global vector $\mathbf{d}$, thus $\mathbf{d}^\alpha = \mathbf{D}_\alpha \mathbf{d}$. Using this transformation, equation \eqref{eq_con:dadotFa} can be written in terms of the global kinematic variables as
\begin{equation} \label{eq_con:ddotF}
\dot{\mathbf{d}}\cdot\mathbf{D}_\alpha^T\left[\mathbf{V}^\alpha+\mathbf{T}^\alpha-\mathbf{F}^\alpha\right]=0.
\end{equation}
For equation \eqref{eq_con:ddotF} to be valid for $\dot{\mathbf{d}}\neq 0 $, we must have
\begin{equation}
\mathbf{D}_\alpha^T\left[\mathbf{V}^\alpha+\mathbf{T}^\alpha-\mathbf{F}^\alpha\right]=\mathbf{0},
\end{equation}
which constitutes the system of nonlinear equations governing the equilibrium of the element $e$. Assembling these equations over all the elements constituting the body, and, in the case of multiple bodies, over all the bodies in the domain, yields the global equilibrium equations
\begin{equation}
\mathbf{M}\ddot{\mathbf{d}}+\mathbf{T}\left(\mathbf{d}\right)-\mathbf{F}=\mathbf{0},
\end{equation}
where $\mathbf{M}$ is the global mass matrix.
%
%
\subsection{Implementation of the contact constraints} \label{sec_con:ConCons}

The equations obtained above can be extended to constrained systems by incorporating the contact constraints into the energy formulation via Lagrange multipliers. The solution of the constrained problem corresponds to the extremum of a modified energy functional
\begin{equation} \label{eq_con:Ehat}
\hat{E}=E - \sum_{i\in N_c}\lambda_i^c g_i^c\left(\mathbf{d}\right)
\end{equation}
with the condition
\begin{equation}
g_i^c\left(\mathbf{d}\right) = 0 \ \ \ \ \forall \, i\in N_c,
\end{equation}
where $g_i^c(\mathbf{d})$ are the contact constraint functions and $N_c$ is the set of active (binding) constraints.

The stationary point of the modified energy functional of equation \eqref{eq_con:Ehat} corresponds to the solution of the system of nonlinear equations:
\begin{align} \label{eq_con:equil}
\mathbf{M}\ddot{\mathbf{d}}+\mathbf{T}\left(\mathbf{d}\right)-\mathbf{F}+\mathbf{F}^c\left(\mathbf{d}\right) &= \mathbf{0} \\ \label{eq_con:gc}
g_i^c\left(\mathbf{d}\right) &= 0 \ \ \ \ \forall \, i\in N_c,
\end{align}
where
\begin{equation}
\mathbf{F}^c\left(\mathbf{d}\right)\equiv - \sum_{i\in N_c}\lambda_i^c \nabla_{\mathbf{d}} g_i^c\left(\mathbf{d}\right)
\end{equation}
is the vector of active contact forces.

Following the active set strategy, the contact constraints are included one by one into $N_c$, starting with the most violated one. The Lagrange multipliers represent the negative of the contact pressure and must satisfy the optimality condition $\lambda_i^c\geq0 \ \forall \ i \in N_c$. When this condition is violated, the corresponding constraint has to be removed from the active set. The solution to equations \eqref{eq_con:equil} and \eqref{eq_con:gc} needs to be updated each time a constraint is added/removed from the active set. Note that the gradient of the contact constraints
yields discrete contact forces at the contact nodes. The contact pressure distribution can be computed using the standard ways of obtaining surface tractions.
%
%
%
%
\section{The oriented volume approach for the formulation of non-smooth contact} \label{sec_con:OrVol} 
%
%
\subsection{The oriented volume contact constraint} \label{sec_con:OrVolCons}
	\begin{figure}
	\centering
	\includegraphics[clip]{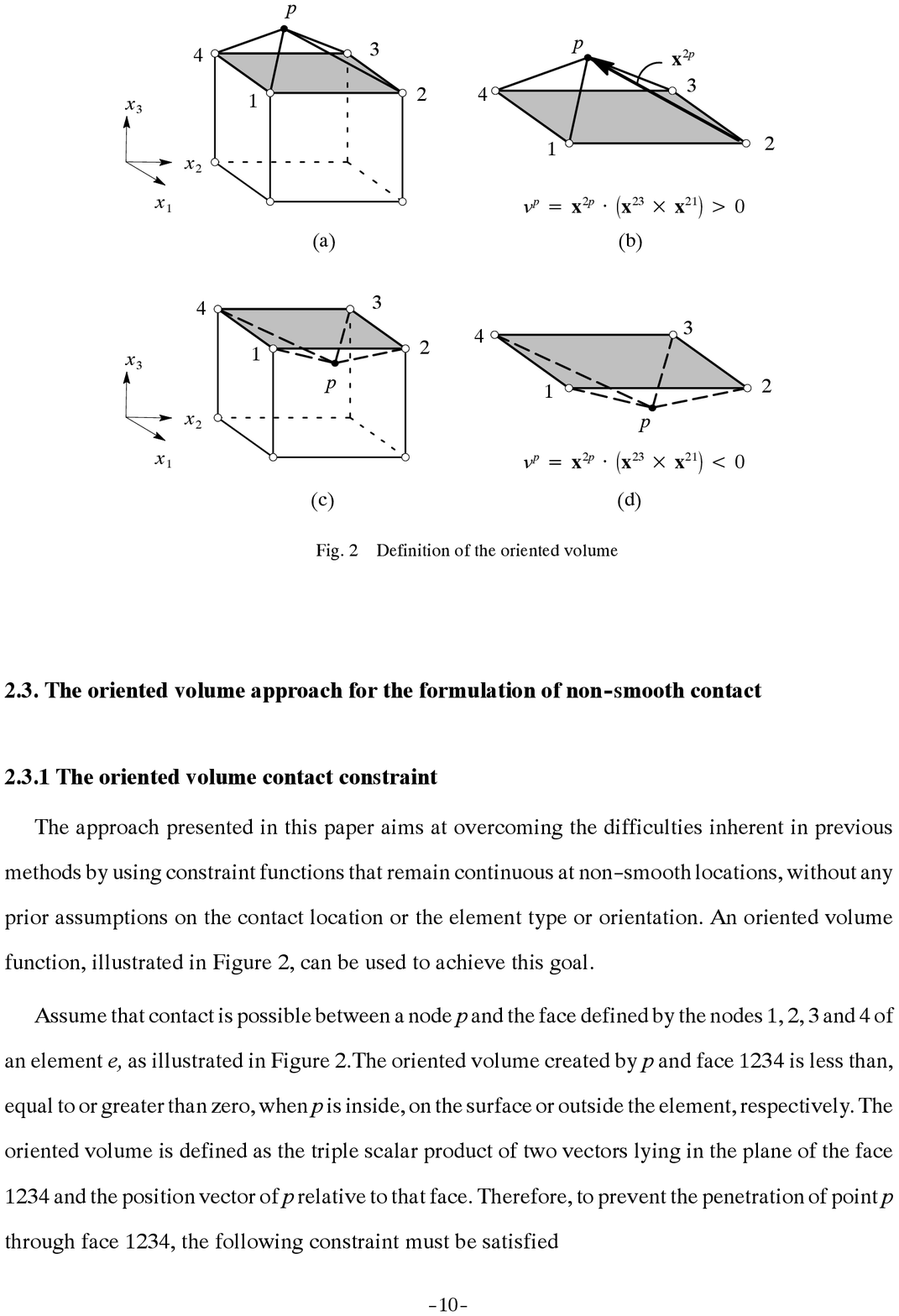}
	\caption{Definition of the oriented volume for a point $p$ with respect to an element surface $1234$ \label{fig:con_orvol}}
	\end{figure}
The approach presented here aims at overcoming the difficulties inherent in previous methods by using constraint functions that remain continuous at non-smooth locations, without any prior assumptions on the contact location or the element type or orientation. An oriented
volume function, illustrated in Figure \ref{fig:con_orvol}, is used to achieve this goal.

Assume that contact is possible between a node $p$ and the facet defined by the nodes $1,2,3,4$ of an element $e$, as illustrated in Figure \ref{fig:con_orvol}. The oriented volume created by $p$ and facet $1234$ is less than, equal to or greater than zero, when $p$ is inside, on the surface or outside the element, respectively. The oriented volume is defined as the triple scalar product of two vectors lying in the plane of the facet $1234$ and the position vector of $p$ relative to that facet $\mathbf{x}^{2p}$. To prevent the penetration of point $p$ through facet $1234$, the following constraint must be satisfied
\begin{equation} \label{eq_con:vp}
v^p = \mathbf{x}^{2p}\cdot\left(\mathbf{x}^{12}\times\mathbf{x}^{23}\right) \geq 0,
\end{equation}
where $\mathbf{x}^{ij} = \mathbf{x}^j - \mathbf{x}^{i}$ is the vector pointing from point $i$ to $j$. Note that the cross product defines the normal to the surface. Equation \eqref{eq_con:vp} computes the oriented volume of the paralelipiped created by the point $p$ and the surface $1234$, not that of the actual pyramid. The volume of the pyramid is proportional to that of the paralelipiped, and the constant of proportionality is not relevant to the result, and will therefore be ignored. In the case where contact occurs through an edge or a corner of the element, the oriented volume created by $p$ and each of the surfaces connected at the edge/corner becomes negative. The resulting oriented volumes
can be used as independent contact constraint functions.

The oriented volume, as defined above, can only be computed for 3-node triangular and 4-node quadrilateral faces, and in the latter case only if none of the four nodes defining the surface displaces outside the plane defined by the remaining three. This assumption may not hold in the presence of large deformations. Also, for T6 and Q8 elements, the higher-order interpolation, enforced by the presence of a central node on each edge, rules out the possibility of using this approach to compute the oriented volume, as the lines connecting each two nodes will not necessarily remain straight.

\begin{figure}
	\centering
	\includegraphics[clip]{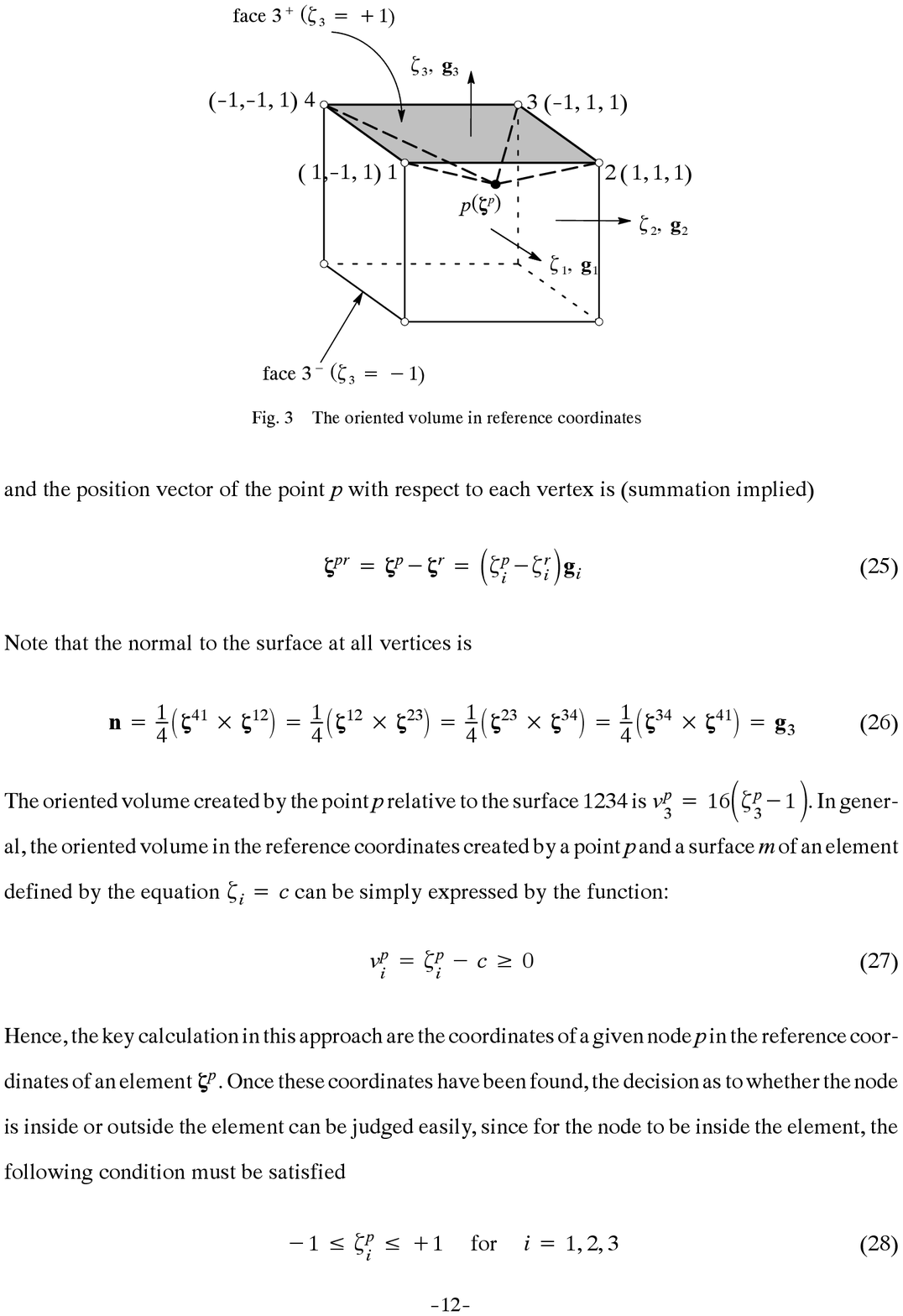}
	\caption{The oriented volume in reference coordinates  \label{fig:con_refvol}}
	\end{figure}

Alternatively, we can compute the oriented volume in the reference (parent) coordinates of the element, where the straight-edge and flat-facet assumptions hold regardless of element order and deformation, as shown in Figure \ref{fig:con_refvol}. An added benefit in this case would be that the coordinates of the nodes in the reference geometry are known a priori, which reduces the effort required to compute the volume.

\begin{figure}
	\centering
	\includegraphics[clip]{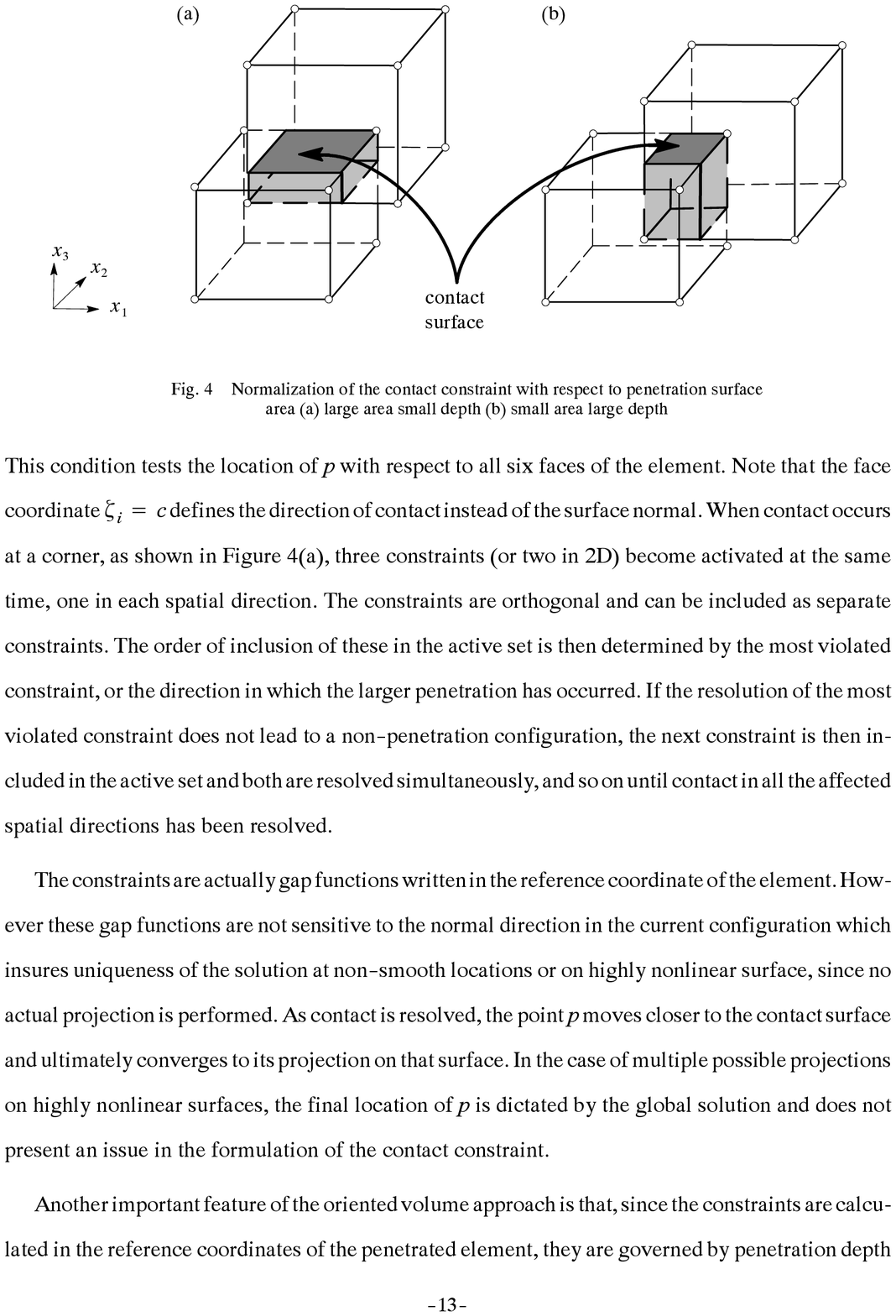}
	\caption{Normalization of the contact constraint with respect to penetration surface area (a) large area small depth (b) small area large depth  \label{fig:con_surf}}
	\end{figure}

Consider the solid block in the reference coordinates $\zeta_i$, with unit vectors $\mathbf{g}_i$, shown in Figure \ref{fig:con_refvol}. Let us compute the
oriented volume of the node $p$ relative to the surface $1234$. The coordinates $\zeta_i$ of the surface vertices are $(1,-1,1)$, $(1,1,1)$, $(-1,1,1)$ and $(-1,-1,1)$, respectively. The surface edges are denoted by the vectors
\begin{equation}
\begin{array}{ll}
\mathbf{\zeta}^{41}=\mathbf{\zeta}^1-\mathbf{\zeta}^4=2\,\mathbf{g}_1, & \mathbf{\zeta}^{12}=\mathbf{\zeta}^2-\mathbf{\zeta}^1=2\,\mathbf{g}_2, \\
\mathbf{\zeta}^{23}=\mathbf{\zeta}^3-\mathbf{\zeta}^2=-2\,\mathbf{g}_1, & \mathbf{\zeta}^{34}=\mathbf{\zeta}^4-\mathbf{\zeta}^3=-2\,\mathbf{g}_2,
\end{array}
\end{equation}
and the position vector of the point $p$ with respect to each vertex is (summation implied)
\begin{equation}
\mathbf{\zeta}^{rp} = \mathbf{\zeta}^p-\mathbf{\zeta}^r = \left(\zeta_i^p-\zeta_i^r\right)\,\mathbf{g}_i.
\end{equation}
Note that the normal to the surface at all vertices is
\begin{equation}
\mathbf{n}=\frac{1}{4}\left(\mathbf{\zeta}^{41}\times\mathbf{\zeta}^{12}\right)=\frac{1}{4}\left(\mathbf{\zeta}^{12}\times\mathbf{\zeta}^{23}\right)=\frac{1}{4}\left(\mathbf{\zeta}^{23}\times\mathbf{\zeta}^{34}\right)=\frac{1}{4}\left(\mathbf{\zeta}^{34}\times\mathbf{\zeta}^{41}\right)=\mathbf{g}_3.
\end{equation}
The oriented volume created by the point $p$ relative to the surface $1234$ is $v_3^p=16(\zeta_3^p - 1)$. In general, the oriented volume
in the reference coordinates created by a point $p$ and a surface $m$ of an element defined by the equation $\zeta_i = c$ can be simply expressed by the function:
\begin{equation}
v_i^p=\zeta_i^p-c \geq 0.
\end{equation}
Hence, the key step in this approach is the calculation of the coordinates of a given node $p$ in the reference coordinates $(\mathbf{\zeta}^p)$. Once these coordinates are found, the decision as to whether the node is inside or outside the element can be judged easily, since for the node to be inside the element, the following condition must be satisfied
\begin{equation}
-1 \leq \zeta_i^p-c \leq 1 \ \ \ \ \mathrm{for} \ \ i=1,2,3.
\end{equation}
This condition tests the location of $p$ with respect to all six facets of the element. Note that, instead of the surface normal, the face coordinate $\zeta_i=c$ defines the direction of contact. When contact occurs at a corner, as shown in Figure \ref{fig:con_surf} (a), three constraints (or two in 2D) become activated at the same time, one in each spatial direction. The constraints are independent and can be treated as separate constraints. The order of inclusion of these in the active set is then determined by the most violated constraint, or the direction
in which the larger penetration has occurred. If the resolution of the most violated constraint does not lead to a all-positive-volume configuration, the next constraint is then included in the active set and both are resolved simultaneously, and so on until contact in all the affected spatial directions has been resolved. This approach recalls the treatment of multiple yield constraints in non-smooth plasticity \cite{simoIJNME88}.

The constraints are, in fact, \textit{gap} functions written in the reference coordinate of the element. However these gap functions are not sensitive to the normal direction in the current configuration which ensures uniqueness of the solution at non-smooth locations or on highly nonlinear surface, since no actual projection is performed. As contact is resolved, the point $p$ moves closer to the contact surface and ultimately converges to its projection on that surface. In the case of multiple possible projections on highly nonlinear surfaces, the final location of $p$ is dictated by the global solution and does not present an issue in the formulation of the contact constraint.

Another important feature of the oriented volume approach is that, since the constraints are calculated in the reference coordinates of the penetrated element, they are governed by penetration depth only, and the area of contact is not a contributing factor. Consider, for example
the two contact scenarios depicted in Figure \ref{fig:con_surf}. The penetration volume is the same in both cases. However, the depth is larger in case (b), and therefore the mass penetration in case (b) is more critical than it is in case (a).

\begin{figure}
	\centering
	\includegraphics[clip]{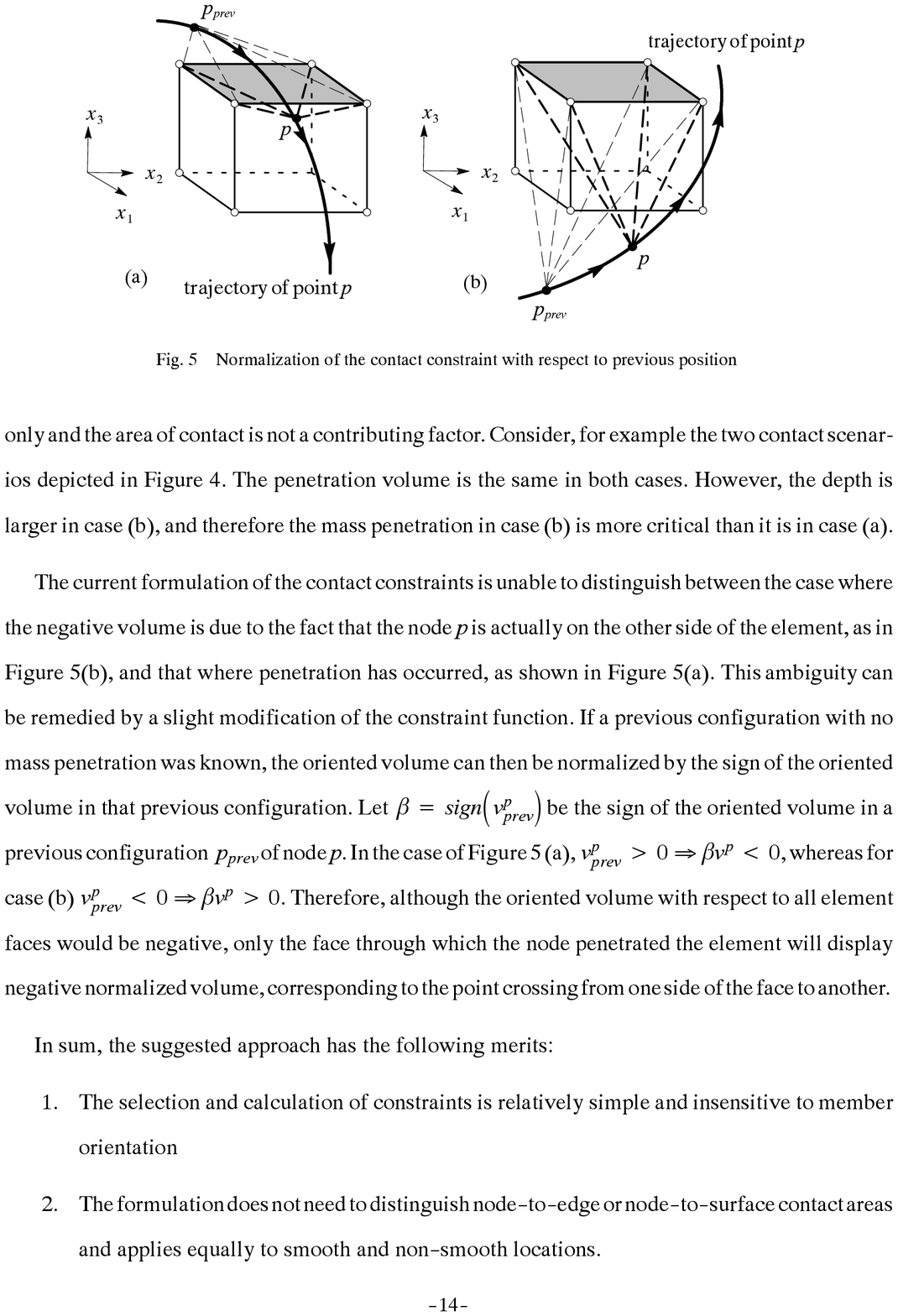}
	\caption{Normalization of the contact constraint with respect to previous position  \label{fig:con_norm}}
	\end{figure}

The current formulation of the contact constraints is unable to distinguish between the case where the negative volume is due to the fact that the node $p$ is actually on the other side of the element, as in Figure \ref{fig:con_norm} (b), and that where penetration has occurred, as shown in Figure \ref{fig:con_norm} (a). This ambiguity can be remedied by a slight modification of the constraint function. If a previous configuration with no mass penetration
was known, the oriented volume can then be normalized by the sign of the oriented volume in that previous configuration. Let $\beta = \mathrm{sign} \left(v_{prev}^p\right)$ be the sign of the oriented volume in a previous configuration $p_{prev}$ of node $p$. In the case of Figure \ref{fig:con_norm} (a), $v_{prev}^p>0 \Rightarrow \beta v^p<0$, whereas for case (b) $v_{prev}^p<0 \Rightarrow \beta v^p>0$. Therefore, although the oriented volume with respect to all element facets would be negative, only the facet through which the node penetrated the element will display negative normalized volume, corresponding to the point crossing from one side of the face to another.

In summary, the suggested approach has the following merits:
\begin{enumerate}
\item The selection and calculation of constraints is relatively simple and insensitive to member orientation.
\item The formulation does not need to distinguish node-to-edge or node-to-surface contact areas and applies equally to smooth and non-smooth locations. 
\item The constraint works well for critical contact scenarios such as at edges and corners or on highly nonlinear surfaces.
\item The constraint has a wider range of applicability than the volume/surface penetration approach of Kane et al. \cite{kaneCMAME99}.
\end{enumerate}
%
%
%
	\begin{figure}
	\centering
	\includegraphics[clip]{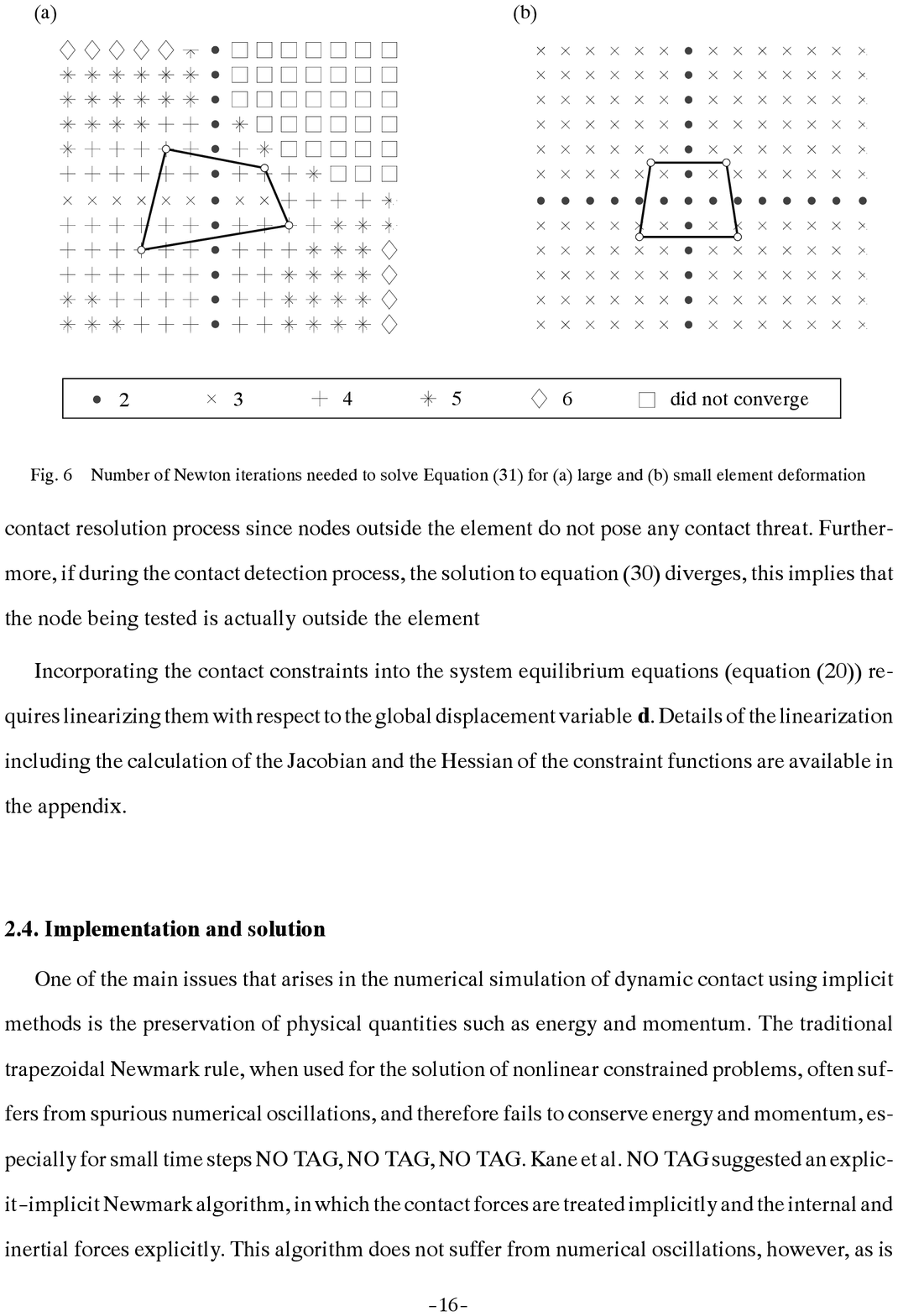}
	\caption{Number of Newton iterations needed to solve Equation (31) for (a) large and (b) small element deformation  \label{fig:con_sol}}
	\end{figure}
\subsection{Finding the $\mathbf{\zeta}^p$ coordinates} \label{sec_con:Zcoord}
The position of a given node $p$ in the reference coordinates of a given element can be obtained, for the three-dimensional case, from the solution of the system of nonlinear equations
\begin{equation} \label{eq_con:xp}
\mathbf{x}^p = N^\alpha\left(\mathbf{\zeta}^p\right)\mathbf{x}^\alpha
\end{equation}
for $\mathbf{\zeta}^p$. If the shape functions $N^\alpha\left(\mathbf{\zeta}^p\right)$ are linear, the solution can be reached in one Newton iteration. The computational cost increases for higher-order elements and in the presence of large deformations, but it remains within
reasonable bounds if the shape of the element does not deteriorate. 

Consider, for example, the 4-node quadrilateral bilinear element, for which the shape functions are given by
\begin{equation} \label{eq_con:Na}
N^\alpha\left(\mathbf{\zeta}^p\right)=\frac{1}{4}\left(1+\zeta_1^p\,\zeta_1^\alpha\right)\left(1+\zeta_2^p\,\zeta_2^\alpha\right)
\end{equation}
(summation over the repeated index $\alpha$ from $1$ to the total number of nodes in the element is implied in equation \eqref{eq_con:xp} but not in equation \eqref{eq_con:Na}). Figure \ref{fig:con_sol} shows the number of Newton iterations needed to reach a solution for $\mathbf{\zeta}^p$ for various locations of $p$, in two different element configurations. The element in Figure \ref{fig:con_sol} (a) is chosen in a highly deformed state leading to an obtuse angle at one of its vertices, whereas the deformation of the element in Figure \ref{fig:con_sol} (b) is relatively small. Notice that, in either case, the solution can be reached easily when $p$ is inside the element. As we move further outside the element, and especially in the case of Figure \ref{fig:con_sol} (a), the computational cost increases gradually. In particular, in the neighborhood of an obtuse vertex the algorithm may diverge altogether because the gradient of equation \eqref{eq_con:xp} with respect to $\zeta_i^p$ becomes singular. Based on the observation that this may only happen outside of badly distorted elements, it will not be considered an issue in the contact resolution process since nodes outside the element do not pose any contact threat. Furthermore, if during the contact detection process, the solution to equation \eqref{eq_con:Na} diverges, this implies that the node being tested is actually outside the element.

Incorporating the contact constraints into the system equilibrium equation \eqref{eq_con:equil} requires linearizing them with respect to the global displacement variable $\mathbf{d}$. Details of the linearization including the calculation of the Jacobian and the Hessian of the constraint functions are provided in Appendix \ref{AppA:ConsLin}.
%
%
%
%
\section{Implementation and solution} \label{sec_con:ImplSol} 
One of the main issues arising in the numerical simulation of dynamic contact using implicit methods is the preservation of physical quantities such as energy and
momentum. The traditional trapezoidal Newmark rule, when used for the solution of nonlinear constrained problems, often suffers from spurious numerical oscillations, and therefore fails to conserve energy and momentum, especially for small time steps \cite{chungJAM93}-\cite{laursenIJNME97}. Kane et al. \cite{kaneCMAME99} suggested an explicit-implicit Newmark algorithm, in which the contact forces are treated implicitly and the internal and inertial forces explicitly. This algorithm does not suffer from numerical oscillations, however, as is the case with all explicit schemes, a very small time step must be used to ensure stability. Alternatively, implicit schemes that possess numerical dissipation such as the HHT method \cite{hilberEESD77}, or more recently, the Generalized $\alpha$ method \cite{chungJAM93} have been used.

In the HHT method, numerical dissipation is introduced by writing the equations of motion in terms of the displacement vector $\mathbf{d}_{n+\alpha}$, calculated as a linear combination of $\mathbf{d}_n$ and $\mathbf{d}_{n+1}$
\begin{equation}
\mathbf{d}_\alpha = \alpha \,\mathbf{d}_{n+1} + (1-\alpha)\mathbf{d}_{n},
\end{equation}
where $n$ denotes the time step index and $\mathbf{d}_m$ should be read ''$\mathbf{d}$ evaluated at time $t = m \Delta t$''. The parameter $\alpha$ controls the level of high-frequency dissipation. For nonlinear systems, the interpolation can also be applied to the internal forces and to the forcing term while the acceleration values
at state $n + 1$ are used. The Generalized $\alpha$ scheme \cite{chungJAM93} extends HHT by applying the same concept to both the displacement/forcing and the acceleration vectors using different values of the interpolation parameter, namely $\alpha_H$ for the displacement/forcing and $\alpha_B$ for the acceleration. Both methods are unconditionally stable and second-order accurate for linear systems. Czekanzi et al. \cite{czekanziCNME01} proposed optimal values for the Generalized $\alpha$ parameters based on a user-defined level of high-frequency dissipation for contact problems.

Dissipative schemes have proven very efficient in filtering high-frequency numerical oscillations for unconstrained dynamic problems. However, since the equations
of equilibrium include asynchronous displacements and accelerations, energy and angular momentum conservation for constrained systems could not been proven rigorously. Moreover, the features of second-order accuracy and unconditional stability of these schemes cannot be guaranteed for nonlinear systems. These issues prompted researchers to develop new families of methods, such as the high-frequency dissipative schemes of Armero et al. \cite{armeroCMAME01pI}\cite{armeroCMAME01pII} or the variational time integrators \cite{kaneIJNME00}\cite{lewIJNME04}, that conserve energy and momenta for nonlinear systems. These methods can be substantially more computationally intensive as compared to ''standard'' schemes such as Newmark or HHT.

The main handicap in the standard implicit time integration schemes is their inability to accommodate the sudden change in the momentum of the contacting bodies at the
moment of impact. Based on this observation, the Decomposition Contact Response method \cite{cirakIJNME05} suggests decomposing the solution in two separate phases, before and after contact. The exchange of momentum and energy due to impact is then accounted for in the solution. Although very robust, this method can be very costly in the case of multiple collisions in a single time step. Therefore, it has only been applied explicitly, using a predictor-corrector approach to estimate the final configuration assuming all contact events occur at the end of the time step.

An implicit scheme that was found to experience less numerical instability than the Newmark family of methods while providing a more robust approach towards the verification of the conservation of energy is the midpoint rule. Unlike the trapezoidal rule, the midpoint rule conserves angular momentum exactly \cite{simoCMAME92}, but energy conservation is only guaranteed for linear systems. The energy conservation property of the midpoint rule for constrained systems was investigated by Bachau \cite{bachauMSD00} and Lenz \cite{lensMSD04}, and was shown to be satisfied provided that the contact constraints produce no work over the time step. For static and persistent contact, this condition is satisfied automatically, but for general dynamic simulations that may include impact and separation of the contacting bodies, additional temporal discretization of the contact constraints is needed.

In the following section we examine the energy conservation properties of the suggested formulation using the midpoint time integration scheme. The choice of midpoint
is motivated by the relative simplicity and efficiency of this scheme in simulating nonlinear dynamic systems in general.
%
%
\subsection{Solution procedure using the midpoint rule} \label{sec_con:SolProc}
The midpoint rule is a one-step time integration scheme in which equilibrium is enforced at the middle of each time interval $\left[ t_n , t_{n+1} \right]$, assuming the following relationships hold:
\begin{align}
\dot{\mathbf{d}}_{n+1/2} &= \frac{1}{2}\left(\dot{\mathbf{d}}_{n+1}+\dot{\mathbf{d}}_{n}\right)=\frac{1}{\Delta t}
\left(\mathbf{d}_{n+1}-\mathbf{d}_{n}\right)  \\
\ddot{\mathbf{d}}_{n+1/2} &= \frac{1}{2}\left(\ddot{\mathbf{d}}_{n+1}+\ddot{\mathbf{d}}_{n}\right)=\frac{1}{\Delta t}
\left(\dot{\mathbf{d}}_{n+1}-\dot{\mathbf{d}}_{n}\right)  \\
\mathbf{d}_{n+1/2} &= \frac{1}{2}\left(\mathbf{d}_{n+1}+\mathbf{d}_{n}\right).
\end{align}
For the unconstrained system, given a known configuration at step $n$, and enforcing equilibrium at time $t_{n+1/2}$, the following system of nonlinear equations can be solved for the unknown configuration at step $n + 1$,
\begin{equation}
\mathbf{M}\ddot{\mathbf{d}}_{n+1/2}+\mathbf{T}_{n+1/2}-\mathbf{F}_{n+1/2}=\mathbf{0}.
\end{equation}
This system of equations corresponds to the extremum of the discrete Lagrangian
\begin{equation}
\mathcal{L}\left(\mathbf{d}_{n+1/2}\right) = \frac{\Delta t^2}{2}\ddot{\mathbf{d}}_{n+1/2}\cdot\mathbf{M}\ddot{\mathbf{d}}_{n+1/2}+\mathit{\Phi}\left(\mathbf{d}_{n+1/2}\right)-\mathbf{d}_{n+1/2}\cdot\mathbf{F}_{n+1/2},
\end{equation}
where $\mathbf{T}_{n+1/2}=\partial\mathit{\Phi}_{n+1/2} / \partial\mathbf{d}_{n+1/2}$ and
\begin{equation}
\ddot{\mathbf{d}}_{n+1/2}=\frac{4}{\Delta t^2}\left(\mathbf{d}_{n+1/2}-\mathbf{d}_n\right)-\frac{2}{\Delta t}\dot{\mathbf{d}}_n.
\end{equation}
Therefore, the solution to the constrained system corresponds to the extremum of the modified discrete functional
\begin{equation}
\hat{\mathcal{L}}\left(\mathbf{d}_{n+1/2},\mathbf{\lambda}^c\right) = \mathcal{L}\left(\mathbf{d}_{n+1/2}\right) -\frac{1}{2}\sum_{i\in N_c}\lambda_i^c g_i^c\left(\mathbf{d}_{n+1}\right).
\end{equation}
Taking the variations of $\hat{\mathcal{L}} \left( \mathbf{d}_{n+1/2},\mathbf{\lambda}^c \right)$ with respect to $\mathbf{d}_{n+1/2}$ and $\mathbf{\lambda}^c$ yields
\begin{equation}
\mathbf{M}\left[\frac{4}{\Delta t^2}\left(\mathbf{y}-\mathbf{d}_n\right)-\frac{2}{\Delta t}\dot{\mathbf{d}}_n\right]+\mathbf{T}\left(\mathbf{y}\right)-\mathbf{F}_{n+1/2}+\mathbf{F}^c\left(2 \mathbf{y}-\mathbf{d}_n\right) = \mathbf{0}
\end{equation}
\begin{equation}
g_i^c\left(2\mathbf{y}-\mathbf{d}_n\right)= 0, \ \ \ \ \forall \,i \in N_c
\end{equation}
where $\mathbf{y}\equiv\mathbf{d}_{n+1/2}$ is the independent kinematic unknown vector, and
\begin{equation}
\mathbf{F}^c\left(2\mathbf{y}-\mathbf{d}_n\right)\equiv - \sum_{i\in N_c}\lambda_i^c \nabla g_i^c\left(\mathbf{d}_{n+1}\right).
\end{equation}
Note that $\mathbf{T}(\mathbf{y})$ should be read as $\mathbf{T}$ evaluated at $\mathbf{y}$ and so on. From this point forward, we will drop the subscript on the symbol and the gradient of a field should be interpreted to be with respect to the argument of that field, unless otherwise specified.
%
%
\subsection{Energy conservation} \label{sec_con:Econs}
For the energy to be conserved over a time interval $\left[ t_n , t_{n+1} \right]$, the integral of its rate of change over that interval must be equal to zero
\begin{align} \notag
E_{n+1}-E_n &= \int_{t_n}^{t_{n+1}}\frac{dE}{dt} dt  \\
&= \int_{t_n}^{t_{n+1}}\dot{\mathbf{d}}\cdot\left[\mathbf{V}+\mathbf{T}+\mathbf{F}^c-\mathbf{F}\right]dt = 0.
\end{align}
In the context of the midpoint time integration scheme, this condition becomes
\begin{equation}
E_{n+1}-E_n \approx \Delta \mathbf{d}\cdot\left[\mathbf{M}\ddot{\mathbf{d}}_{n+1/2}+\mathbf{T}_{n+1/2}+\mathbf{F}_{n+1/2}-\mathbf{F}^c\right] = 0,
\end{equation}
where $\Delta\mathbf{d} = \mathbf{d}_{n+1}-\mathbf{d}_{n} = 2 \left( \mathbf{y}-\mathbf{d}_{n} \right)$ is the displacement jump over the time step.

As demonstrated by Bachau \cite{bachauMSD00} and Lenz \cite{lensMSD04}, for the total energy of the constrained system to be conserved, the work done by the constraint forces must vanish over the time step. Hence,
\begin{equation}
\mathbf{F}^c\cdot\Delta\mathbf{d}=-\sum_{i\in N_c}\lambda_i^c \nabla g_i^c\left(\mathbf{d}_{n+1}\right)\cdot\Delta\mathbf{d} = \sum_{i\in N_c}\mathbf{F}_i^c\cdot\Delta\mathbf{d}=0.
\end{equation}
First, note that, at the moment of contact, node $p$ is either inside or at the surface of the element. Let $\hat{p}$ be the target location of $p$ on the surface of the element, that is, its location after contact is resolved. Since $\hat{p}$ is a point of the contact element, its deformed coordinates, displacements and incremental displacements over the time step satisfy the equations
\begin{align}
\mathbf{x}^{\hat{p}} &= \sum_\alpha N^\alpha\left(\mathbf{\zeta}^{\hat{p}}\right)\mathbf{x}^\alpha \\
\mathbf{d}^{\hat{p}} &= \sum_\alpha N^\alpha\left(\mathbf{\zeta}^{\hat{p}}\right)\mathbf{d}^\alpha \\
\Delta\mathbf{d}^{\hat{p}} &= \sum_\alpha N^\alpha\left(\mathbf{\zeta}^{\hat{p}}\right)\Delta\mathbf{d}^\alpha.
\end{align}
For each constraint in $i \in N_c$, the gradient of the constraint function with respect to the nodal displacement vector $\mathbf{d}$ is given by the equation (see Appendix \ref{AppA:ConsLin})
\begin{equation}
\nabla g_i^c = \left[\mathbf{D}_p^T-N^\alpha\left(\mathbf{\zeta}^p\right)\mathbf{D}_\alpha^T\right]\mathbf{f}_p^c,
\end{equation}
where $\mathbf{f}_p^c$ contains information about the direction of contact and $\mathbf{D}_\alpha^T$ plays the role of applying the contact force vector at node $\alpha$. Therefore, the work performed by the contact constraint over the time step can be computed as
\begin{equation} \label{eq_con:conWa}
\lambda_i^c\nabla g_i^c\cdot\Delta\mathbf{d}=\lambda_i^c \left[\mathbf{D}_p^T-N^\alpha\left(\mathbf{\zeta}^p\right)\mathbf{D}_\alpha^T\right]\mathbf{f}_p^c\cdot\Delta\mathbf{d}.
\end{equation}
Note that,
\begin{align} \notag
\left[\mathbf{D}_p^T-N^\alpha\left(\mathbf{\zeta}^p\right)\mathbf{D}_\alpha^T\right]\mathbf{f}_p^c\cdot\Delta\mathbf{d} &= \mathbf{f}_p^c\cdot\left[\mathbf{D}_p-N^\alpha\left(\mathbf{\zeta}^p\right)\mathbf{D}_\alpha\right]\Delta\mathbf{d}   \\  
\label{eq_con:conWb}
&=\mathbf{f}_p^c\cdot\left[\Delta\mathbf{d}^p-N^\alpha\left(\mathbf{\zeta}^p\right)\Delta\mathbf{d}^\alpha\right].
\end{align}
Given that $\mathbf{\zeta}^p = \mathbf{\zeta}^{\hat{p}}$ at the solution, equations \eqref{eq_con:conWa} and \eqref{eq_con:conWb} lead to
\begin{align} \notag
\lambda_i^c\nabla g_i^c\cdot\Delta\mathbf{d} &= \lambda_i^c\mathbf{f}_p^c\cdot\left[\Delta\mathbf{d}^p-\Delta\mathbf{d}^{\hat{p}}\right] \\
\label{eq_con:conWc}
&= 0 \ \ \ \ \textrm{if} \ \ \lambda_i^c = 0 \ \  \textrm{or} \ \ \mathbf{f}_p^c\cdot\Delta\mathbf{d}^p = \mathbf{f}_p^c\cdot\Delta\mathbf{d}^{\hat{p}}.
\end{align}
The work of the contact forces clearly vanishes when either the contact force is zero  (just before impact or after separation) or when $\mathbf{f}_p^c \cdot \Delta \mathbf{d}^p = \mathbf{f}_p^c \cdot \Delta \mathbf{d}^{\hat{p}}$ (after impact or before separation). In a continuum setting, this condition translates to $\mathbf{f}_p^c \cdot\dot{\mathbf{d}}^p = \mathbf{f}_p^c \cdot\dot{\mathbf{d}}^{\hat{p}}$, which is the well-known persistency condition predicted by wave propagation. The logical approach to follow to ensure an accurate solution is to integrate these two phases separately. However, in the context of an implicit single-point time integration scheme, such as the midpoint rule used herein, if contact has been detected and $\lambda_i^c \neq 0$, equation \eqref{eq_con:conWc} implies that
\begin{equation}
\mathbf{f}_p^c\cdot\dot{\mathbf{d}}_{n+1/2}^p\,\Delta t = \mathbf{f}_p^c\cdot\dot{\mathbf{d}}_{n+1/2}^{\hat{p}}\,\Delta t
\end{equation}
\begin{equation} \label{eq_con:persistency}
\frac{1}{2}\mathbf{f}_p^c\cdot\left[\dot{\mathbf{d}}_{n}^p+\dot{\mathbf{d}}_{n+1}^p\right] = \frac{1}{2}\mathbf{f}_p^c\cdot\left[\dot{\mathbf{d}}_{n}^{\hat{p}}+\dot{\mathbf{d}}_{n+1}^{\hat{p}}\right].
\end{equation}
Therefore, for energy to be conserved, the persistency condition has to hold in an average sense over the time step (for the midpoint rule). This condition corresponds exactly to the algorithmic gap rate proposed by Laursen et al. \cite{laursenIJNME97}. Note that equation \eqref{eq_con:persistency} can be alternatively written as
\begin{align}
\mathbf{f}_p^c\cdot\left[\mathbf{d}_{n+1}^p-\mathbf{d}_{n}^p\right] &= \mathbf{f}_p^c\cdot\left[\mathbf{d}_{n+1}^{\hat{p}}-\mathbf{d}_{n}^{\hat{p}}\right] \\
\label{eq_con:persistency_x}
\Rightarrow
\mathbf{f}_p^c\cdot\left[\left(\mathbf{d}_{n+1}^p+\mathbf{X}^p\right)-\left(\mathbf{d}_{n}^p+\mathbf{X}^p\right)\right] &= \mathbf{f}_p^c\cdot\left[\left(\mathbf{d}_{n+1}^{\hat{p}}+\mathbf{X}^{\hat{p}}\right)-\left(\mathbf{d}_{n}^{\hat{p}}+\mathbf{X}^{\hat{p}}\right)\right].
\end{align}
Rearranging the terms in equation \eqref{eq_con:persistency_x}, we find
\begin{align}
\mathbf{f}_p^c\cdot\left[\mathbf{x}_{n+1}^p-\mathbf{x}_{n+1}^{\hat{p}}\right] &= \mathbf{f}_p^c\cdot\left[\mathbf{x}_{n}^p-\mathbf{x}_{n}^{\hat{p}}\right] \\
\label{eq_con:persistency_gap}
\Rightarrow
g_i^c\left(\mathbf{d}_{n+1}\right) &= g_i^c\left(\mathbf{d}_{n}\right),
\end{align}
which is an alternate form of the algorithmic persistency condition of equation \eqref{eq_con:persistency}. We distinguish between the following two cases:

(a) \textit{Persistent and static contact:} In this case we have $g_i^c(\mathbf{d}_n)=0$ and $\mathbf{f}_p^c \cdot\dot{\mathbf{d}}_n^p = \mathbf{f}_p^c \cdot\dot{\mathbf{d}}_n^{\hat{p}}$. Thus, if the persistency condition is enforced at the end of the time step $\mathbf{f}_p^c \cdot\dot{\mathbf{d}}_{n+1}^p = \mathbf{f}_p^c \cdot\dot{\mathbf{d}}_{n+1}^{\hat{p}}$, equation \eqref{eq_con:persistency} is satisfied and energy is conserved. The same result can be obtained by simply enforcing the contact constraint at the end of the time step $g_i^c(\mathbf{d}_{n+1})=0$. From equation \eqref{eq_con:persistency_gap}, the average persistency condition is automatically satisfied and energy is conserved. As a result, the persistency condition is also satisfied at the end of the time step $\mathbf{f}_p^c \cdot\dot{\mathbf{d}}_{n+1}^p = \mathbf{f}_p^c \cdot\dot{\mathbf{d}}_{n+1}^{\hat{p}}$.

(b) \textit{Impact:} Since $g_i^c(\mathbf{d}_n)\geq 0$ and $\mathbf{f}_p^c \cdot\dot{\mathbf{d}}_n^p \neq \mathbf{f}_p^c \cdot\dot{\mathbf{d}}_n^{\hat{p}}$ in this case, enforcing the persistency condition $\mathbf{f}_p^c \cdot\dot{\mathbf{d}}_{n+1}^p = \mathbf{f}_p^c \cdot\dot{\mathbf{d}}_{n+1}^{\hat{p}}$ leads to $\mathbf{f}_p^c \cdot\dot{\mathbf{d}}_{n+1/2}^p \neq \mathbf{f}_p^c \cdot\dot{\mathbf{d}}_{n+1/2}^{\hat{p}}$ and therefore energy is not conserved. This issue is well documented in the literature and has been addressed in various ways. Solberg and Papadopoulos \cite{solbergFEAD98} recommend enforcing the persistency condition at step $n$, at the cost of introduction an energy error of order $\mathcal{O}(h)$, where $h$ is the spatial mesh size. Laursen and Chawla \cite{laursenIJNME97} used the algorithmic gap rate (equation \eqref{eq_con:persistency}) to impose the persistency condition in an average sense over the time step. The disadvantage of this method is that it can lead to geometrically inadmissible configurations, as pointed out in \cite{laursenIJNME97}. The general consensus in the literature is that a velocity correction is needed for nonsmooth events such as impact. Hughes et al. \cite{hughesCMAME76} used the wave propagation properties of the medium to calculate these corrections. This approach, however, may not be straightforward in the general case of multi-dimensional nonlinear elasticity. Laursen and Love \cite{laursenIJNME02} proposed discrete velocity jumps at the contact interface that can be computed as a post-processing step. Bachau \cite{bachauMSD00} and Lenz \cite{lensMSD04} suggested discretizing each active contact constraint, such that the following holds
\begin{equation}
\mathbf{F}_i^c\cdot\Delta\mathbf{d} = g_i^c\left(\mathbf{d}_{n+1}\right)-g_i^c\left(\mathbf{d}_{n}\right)=0.
\end{equation}
Since the contact constraints considered herein are written in the reference coordinates of the penetrated element, this approach is problematic, if even possible. Using a similar approach, Hesch and coworkers recently proposed an algorithmic formulation of the contact forces that conserves energy and momentum in the discrete problem \cite{heschIUTAM07,heschIJNME08}.

A complementary way of satisfying energy conservation for impact is by imposing it as an additional constraint via a Lagrange multiplier. Accordingly, the modified discrete Lagrangian to be extremized becomes
\begin{equation}
\hat{\mathcal{L}}\left(\mathbf{d}_{n+1/2},\mathbf{\lambda}^c\right)=\mathcal{L}\left(\mathbf{d}_{n+1/2}\right)-\frac{1}{2}\sum_{i\in N_c}\lambda_i^c g_i^c\left(\mathbf{d}_{n+1}\right)-\frac{1}{2}\lambda^E g^E\left(\mathbf{d}_{n+1},\dot{\mathbf{d}}_{n+1}\right),
\end{equation}
where
\begin{equation}
g^E\left(\mathbf{d}_{n+1},\dot{\mathbf{d}}_{n+1}\right) = \frac{1}{2}\dot{\mathbf{d}}_{n+1}\cdot\mathbf{M}\dot{\mathbf{d}}_{n+1}+\mathit{\Phi}_{n+1}-\frac{1}{2}\dot{\mathbf{d}}_{n}\cdot\mathbf{M}\dot{\mathbf{d}}_{n}-\mathit{\Phi}_{n}-\mathbf{F}_{n+1/2}\cdot\Delta\mathbf{d}
\end{equation}
is the energy conservation constraint. Note that, in $g^E$, $\mathbf{d}_{n}$, $\mathbf{d}_{n+1}$, $\dot{\mathbf{d}}_{n}$ and $\dot{\mathbf{d}}_{n+1}$ can be restricted to the degrees of freedom of the contacting bodies. Defining
\begin{equation} \label{eq_con:FE}
\mathbf{F}^E \equiv - \lambda^E\nabla g^E \left(\mathbf{d}_{n+1},\dot{\mathbf{d}}_{n+1}\right),
\end{equation}
the system of equations to be solved for equilibrium becomes
\begin{equation}\label{eq_con:equil_mid}
\mathbf{M}\ddot{\mathbf{d}}_{n+1/2}+\mathbf{T}_{n+1/2}-\mathbf{F}_{n+1/2}+\mathbf{F}^c+\mathbf{F}^E = \mathbf{0}
\end{equation}
subject to the set of constraints
\begin{equation}
g_i^c\left(\mathbf{d}_{n+1}\right) = 0 \ \ \ \ \forall \,i\in N_c
\end{equation}
and
\begin{equation}
g^E\ ( \mathbf{d}_{n+1},\dot{\mathbf{d}}_{n+1}\ ) = 0.
\end{equation}
Therefore, the role of the energy Lagrange multiplier is to introduce algorithmic forces/accelerations that would produce the velocity corrections necessary to conserve energy during non-smooth events. A single constraint is needed to account for all such events occurring during a given time step $\left[ t_n, t_{n+1} \right]$. Naturally, this extends the number of equations to be solved by one, but the added cost is small compared to the original size of the problem. When the average persistency condition is satisfied, the value of $\lambda^E$ and therefore $\mathbf{F}^E$ results to be zero. This result can be obtained by computing the work done by the forces in equation \eqref{eq_con:equil_mid}
\begin{equation}
\Delta\mathbf{d}\cdot\left[\mathbf{M}\ddot{\mathbf{d}}_{n+1/2}+\mathbf{T}_{n+1/2}-\mathbf{F}_{n+1/2}+\mathbf{F}^c+\mathbf{F}^E\right] = 0.
\end{equation}
It can be shown that the work of the inertia forces over the time step corresponds exactly to the change in kinetic energy. Furthermore, we assume that the work done by the internal and external forces is approximately equal to the change in potential energy (this is an exact equality for linear systems and generally holds up to an error of order $\mathcal{O}(\Delta t)$ for nonlinear systems). As a result,
\begin{equation}
\Delta\mathbf{d}\cdot\left[\mathbf{M}\ddot{\mathbf{d}}_{n+1/2}+\mathbf{T}_{n+1/2}-\mathbf{F}_{n+1/2}\right] \approx g^E = 0,
\end{equation}
which leads to
\begin{equation} \label{eq_con:Econsa}
\Delta\mathbf{d}\cdot\left[\mathbf{F}^c+\mathbf{F}^E\right] = 0.
\end{equation}
Substituting equations \eqref{eq_con:conWc} and \eqref{eq_con:FE} into equation \eqref{eq_con:Econsa} yields
\begin{equation} \label{eq_con:Econsb}
\Delta\mathbf{d}\cdot\left[\mathbf{F}^c+\mathbf{F}^E\right] = -\lambda_i^c\mathbf{f}_p^c\cdot\left[\Delta\mathbf{d}^p-\Delta\mathbf{d}^{\hat{p}}\right] - \lambda^E\nabla g^E\cdot\Delta\mathbf{d} = 0.
\end{equation}
From equation \eqref{eq_con:Econsb}, we can observe that, if the algorithmic persistency condition is satisfied, i.e. if $\mathbf{f}_p^c \cdot\Delta\mathbf{d}^p =\mathbf{f}_p^c \cdot\Delta\mathbf{d}^{\hat{p}}$, then the following holds:
\begin{equation} \label{eq_con:Econsc}
\lambda^E\nabla g^E\cdot\Delta\mathbf{d} = 0.
\end{equation}
Since $ \nabla g^E \cdot \Delta \mathbf{d}\neq 0$ in general, equation \eqref{eq_con:Econsc} yields $\lambda^E=0$. Conversely, if the algorithmic persistency condition is not satisfied, then the work of the (algorithmic) forces introduced by the energy constraint yields the correction needed to counterbalance the work of the contact constraints.
\begin{remark}
Even in the case of impact, if the two contacting bodies are arbitrarily close right before impact such that $g_i^c(\mathbf{d}_n )\approx 0$, the algorithmic persistency condition is satisfied and the total energy is conserved without the need for the additional Lagrange multiplier. Therefore, in the limit of temporal refinement, the Lagrange multiplier approach reverts back to the enforcement of the algorithmic persistency condition.
\end{remark}
\begin{remark}
The Lagrange multiplier approach is applicable to any other time integration procedure. Thus, if the error in integrating the potential energy is relatively large, then a higher-order conservative formulation of the continuum can be used instead of the midpoint rule. It is useful to point out that Hughes et al. \cite{hughesJAM78} implemented a Lagrange multiplier method to achieve conservation of energy for general (unconstrained) nonlinear dynamic systems. Thus, for nonlinear systems in the absence of contact, the energy Lagrange multiplier method coupled with the midpoint rule, as presented herein, is similar to the method of Hughes et al. \cite{hughesJAM78}. In the presence of contact, the Lagrange multiplier serves the purpose of correcting the error in energy due to both the nonlinearity of the system and the non-smooth dynamic contact events.
\end{remark}

In the following section, we describe the implementation of the Lagrange multiplier approach for the conservation of energy using the midpoint rule.
%
%
\subsection{Linearization and solution} \label{sec_con:LinSol}
The discretized equilibrium and constraint equations can be summarized as
\begin{align} \label{eq_con:r}
\mathbf{r}  = &\mathbf{M}\left[\frac{4}{\Delta t^2}\left(\mathbf{y}-\mathbf{d}_n\right) - \frac{2}{\Delta t}\dot{\mathbf{d}}_n\right] + \mathbf{T}\left(\mathbf{y}\right)-\mathbf{F}_{n+1/2}+\mathbf{F}^c\left(2\mathbf{y}-\mathbf{d}_n\right)+\mathbf{F}^E = 0 \\
&g_i^c\left(2\mathbf{y}-\mathbf{d}_n\right) = 0 \ \ \ \ \forall \, i\in N_c  \\
&g^E(2\mathbf{y}-\mathbf{d}_n) = 0.
\end{align}
In these equations,
\begin{align}
\notag
g^E \left(2\mathbf{y}-\mathbf{d}_n\right) &= \frac{1}{2}\left[\frac{4}{\Delta t}\left(\mathbf{y}-\mathbf{d}_n\right)-\dot{\mathbf{d}}_n\right]
\cdot\mathbf{M}\left[\frac{4}{\Delta t}\left(\mathbf{y}-\mathbf{d}_n\right)-\dot{\mathbf{d}}_n\right] - \mathit{\Phi}\left(2\mathbf{y}-\mathbf{d}_n\right) \\
\label{eq_con:gE}
&- \frac{1}{2}\dot{\mathbf{d}}_n\cdot\mathbf{M}\dot{\mathbf{d}}_n - \mathit{\Phi}\left(\mathbf{d}_n\right)-2\mathbf{F}_{n+1/2}\cdot\left(\mathbf{y}-\mathbf{d}_n\right),
\end{align}
$\mathbf{F}^c\left(2\mathbf{y}-\mathbf{d}_n\right)= - \sum \lambda_i^c \nabla g_i^c (2 \mathbf{y}- \mathbf{d}_n)$ as previously defined, and $\mathbf{F}^E = -\lambda^E\nabla g^E$. Since $g^E$ is a discrete (in time) function of $\mathbf{d}$ and $\dot{\mathbf{d}}$, the gradient of the energy conservation constraint should be computed as follows:
\begin{equation} \label{eq_con:gradgE}
\nabla g^E = \partial g^E /\partial\mathbf{d}_{n+1} = \beta\,\partial\left[\frac{1}{2}\dot{\mathbf{d}}_{n+1}\cdot\mathbf{M}\dot{\mathbf{d}}_{n+1}\right]
/\partial\dot{\mathbf{d}}_{n+1}+\partial\mathit{\Phi}\left(\mathbf{d}_{n+1}\right)/\partial\mathbf{d}_{n+1},
\end{equation}
where $\beta\equiv\partial\dot{\mathbf{d}}_{n+1}/\partial{\mathbf{d}}_{n+1}= 2/\Delta t$ for the midpoint rule. System \eqref{eq_con:r}-\eqref{eq_con:gE} can be solved for the kinematic variables $\mathbf{y}$ and the Lagrange multipliers $\mathbf{\lambda}=\left[\mathbf{\lambda}^c, \lambda^E\right]$, using Newton's method. The directional derivatives \begin{footnote}{The directional derivative $D$ of a field $\mathbf{f}(\mathbf{z})$ in the direction of a variation $\mathbf{w}$ is defined as $D \mathbf{f}\left(\mathbf{z}\right)\cdot\mathbf{w}\equiv\frac{d}{d\epsilon}\left[\mathbf{f}\left(\mathbf{z}+\epsilon\mathbf{w}\right)\right]_{\epsilon = 0}=\nabla \mathbf{f}\,\mathbf{w}$} \end{footnote}
 of the system equations with respect to the incremental variables are,
\begin{equation}
D\mathbf{r}\cdot\Delta\mathbf{y} = \tilde{\mathbf{K}}_{eff}\Delta\mathbf{y},
\end{equation}
where
\begin{equation}  \label{eq_con:Keff}
\tilde{\mathbf{K}}_{eff} = \frac{2}{\Delta t^2}\mathbf{M}+\mathbf{K}_t\left(\mathbf{y}\right) - 2\sum_{i\in N_c}\lambda_i^c\mathbf{H}_i^c\left(2\mathbf{y}-\mathbf{d}_n\right)-\lambda^E\mathbf{H}^E.
\end{equation}
In this equation, $\mathbf{H}_i^c$ is the Hessian of the active contact constraint $i$ (see Appendix \ref{AppA:ConsLin}) and $\mathbf{H}^E$ is the Hessian of the energy conservation constraint, calculated by taking the directional derivative of equation \eqref{eq_con:gradgE} in the direction of $\Delta\mathbf{y}$
\begin{align}
\notag
\mathbf{H}^E\Delta\mathbf{y} &= D\left[\frac{2}{\Delta t}\mathbf{M}\left[\frac{4}{\Delta t}\left(\mathbf{y}-\mathbf{d}_n\right)-\dot{\mathbf{d}}_n\right]+\mathbf{T}\left(2\mathbf{y}-\mathbf{d}_n\right) \right]\cdot\Delta\mathbf{y}  \\
&= \left[\frac{8}{\Delta t^2}\mathbf{M}+2\,\mathbf{K}_t\left(2\mathbf{y}-\mathbf{d}_n\right)\right]\Delta \mathbf{y}.
\end{align}
The remaining directional derivatives are
\begin{align}
&D\mathbf{r}\cdot\Delta\lambda_i^c = -\nabla g_i^c\left(2\mathbf{y}-\mathbf{d}_n\right)\Delta\lambda_i^c \ \ \ \ \forall \,i\in N_c \\
&D\mathbf{r}\cdot\Delta\lambda^E = -\left\{\frac{2}{\Delta t}\mathbf{M}\left[\frac{4}{\Delta t}\left(\mathbf{y}-\mathbf{d}_n\right)-\dot{\mathbf{d}}_n\right]+\mathbf{T}\left( 2\mathbf{y}-\mathbf{d}_n\right)-\mathbf{F}_{n+1/2}\right\}\Delta\lambda^E  \\
&D g_i^c\left(2\mathbf{y}-\mathbf{d}_n\right)\cdot \Delta\mathbf{y} = 2 \nabla g_i^c\left(2\mathbf{y}-\mathbf{d}_n\right)\cdot \Delta\mathbf{y} \ \ \ \ \forall \,i\in N_c \\
\notag
&D g^E\left(2\mathbf{y}-\mathbf{d}_n\right)\cdot \Delta\mathbf{y} = 2 \nabla g^E \cdot \Delta\mathbf{y} \\
& \ \ \ \ \ \ = \left\{\frac{4}{\Delta t}\mathbf{M}\left[\frac{4}{\Delta t}\left(\mathbf{y}-\mathbf{d}_n \right)-\dot{\mathbf{d}}_n \right]+2\,\mathbf{T}\left(2\mathbf{y}-\mathbf{d}_n \right)-2\,\mathbf{F}_{n+1/2}\right\}\cdot\Delta\mathbf{y} \\
&D g_i^c\cdot\Delta \lambda_i^c = D g_i^c\cdot\Delta \lambda^E = 0 \ \ \ \ \forall \,i\in N_c \\
&D g^E\cdot\Delta \lambda_i^c = D g^E\cdot\Delta \lambda^E = 0.
\end{align}
Accordingly, the system of linearized equations to be solved at each iteration is
\begin{equation}
\left[ \begin{array}{cc}
  \tilde{\mathbf{K}}_{eff}\left(\mathbf{y}\right) & -\mathbf{J}\left(2\mathbf{y}-\mathbf{d}_n\right)	\\
  2\,\mathbf{J}^T\left(2\mathbf{y}-\mathbf{d}_n\right) & \mathbf{0}
\end{array} \right ]
\left\{ \begin{array}{c} \Delta\mathbf{y}^k \\ \Delta\mathbf{\lambda}^k \end{array} \right\}
= - \left\{ \begin{array}{c} \mathbf{r}\left(\mathbf{y}\right) \\ \mathbf{g}\left(2\mathbf{y}-\mathbf{d}_n\right) \end{array} \right\},
\end{equation}
where $\tilde{\mathbf{K}}_{eff}$ is given by equation \eqref{eq_con:Keff} and,
\begin{align}
\mathbf{J} &= \left[\nabla g_1^c \ \ \nabla g_2^c \ \ ,..., \ \ \nabla g_n^c \ \ \nabla g^E\right] \\
\mathbf{g} &= \left[\,g_1^c \ \ g_2^c \ \ ,..., \ \ g_n^c \ \ g^E\right].
\end{align}
%

%
%
%
%
\section{Numerical examples} \label{sec_con:NumEx} 
This section illustrates the implementation of the suggested approach in the solution of a few examples. In all examples, we assume finite strains/large deformations and a compressible Neo-Hookean material with a strain energy density function given by
\begin{equation}
\psi\left(\mathbf{u}\right) = \frac{\lambda}{2}\left(\mathrm{log}J\right)^2-\mu\,\mathrm{log}J+\frac{\mu}{2}\,\mathrm{tr}\left(\mathbf{C}-\mathbf{I}\right),
\end{equation}
where $J = \textrm{det}(\mathbf{F})$ is the deformation Jacobian, $\mathbf{C}$ is the right Cauchy-Green strain tensor and $\lambda$, $\mu$ are the Lam\'{e} material constants. Consistent units of mass, force, time and length are implied for all numerical quantities used in these examples.
%
%
\subsection{2D example} \label{sec_con:2Dex}
\begin{figure}
	\centering
	\includegraphics[clip]{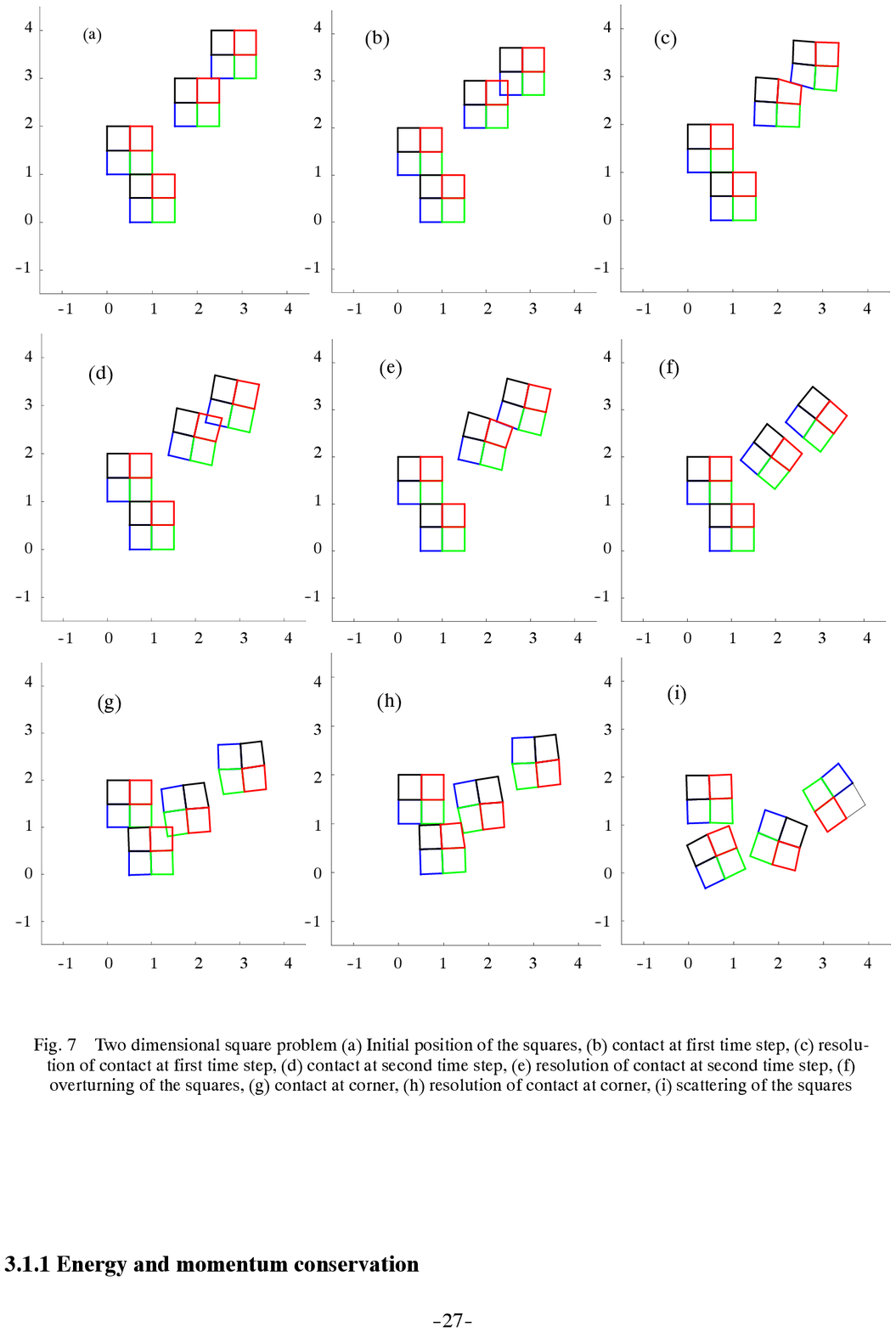}
	\caption{Two dimensional square problem (a) Initial position of the squares, (b) contact at first time step, (c) resolution of contact at first time step, (d) contact at second time step, (e) resolution of contact at second time step, (f) overturning of the squares, (g) contact at corner, (h) resolution of contact at corner, (i) scattering of the squares  \label{fig:con_motion_Q4}}
	\end{figure}

The first example consists of a set of cubes initially set as shown in Figure \ref{fig:con_motion_Q4} (a) \cite{kaneCMAME99}. The cubes are discretized using 4-node quadrilateral elements. The upper-most cube is then given an initial downward velocity of $v_0=-3000$. The cubes have unit side length and the following material properties: $\rho_0 = 1000$, $\lambda = 1.75 \times 10^{11}$ and $\mu = 0.801 \times 10^{11}$ (which is equivalent to the properties $E = 2.151 \times 10^{11}$ and $\nu = 0.343$). The solution is initially carried out over 15 time
steps using time increments of $\Delta t = 7 \times 10^{- 5}$ and a lumped mass matrix.

Contact occurs after the first step as shown in Figure \ref{fig:con_motion_Q4} (b), which clearly induces a large amount of mass penetration. The contact constraint is enforced at this point and the cubes are deformed to preclude mass penetration, as shown in Figure \ref{fig:con_motion_Q4} (c). Subsequently, the upper square keeps moving downwards, leading to contact occurring again at the second time step as shown in Figure \ref{fig:con_motion_Q4} (d). This second contact is resolved leading to the result shown in Figure \ref{fig:con_motion_Q4} (e).

After this point, the bodies start rotating due to the angular momentum transferred by the high impact forces, resulting in the configuration shown in Figure \ref{fig:con_motion_Q4} (f). At this
stage, the two squares separate and the motion happens without any spurious vibrations, due to the conservation of energy during impact. At time $t = 63 \times 10^{- 5}$ another contact event happens at a corner, as shown in Figure \ref{fig:con_motion_Q4} (g), and is successfully resolved in Figure \ref{fig:con_motion_Q4} (h). This particular scenario consists of many repetitive events, all at the corner, and sometimes involving more than two cubes. The ability of the algorithm to treat these multiple events is remarkable.
%
%
%
\subsubsection{Energy and momentum conservation} \label{sec_con:2DexEMcons}
	\begin{figure}
	\centering
	\includegraphics[clip]{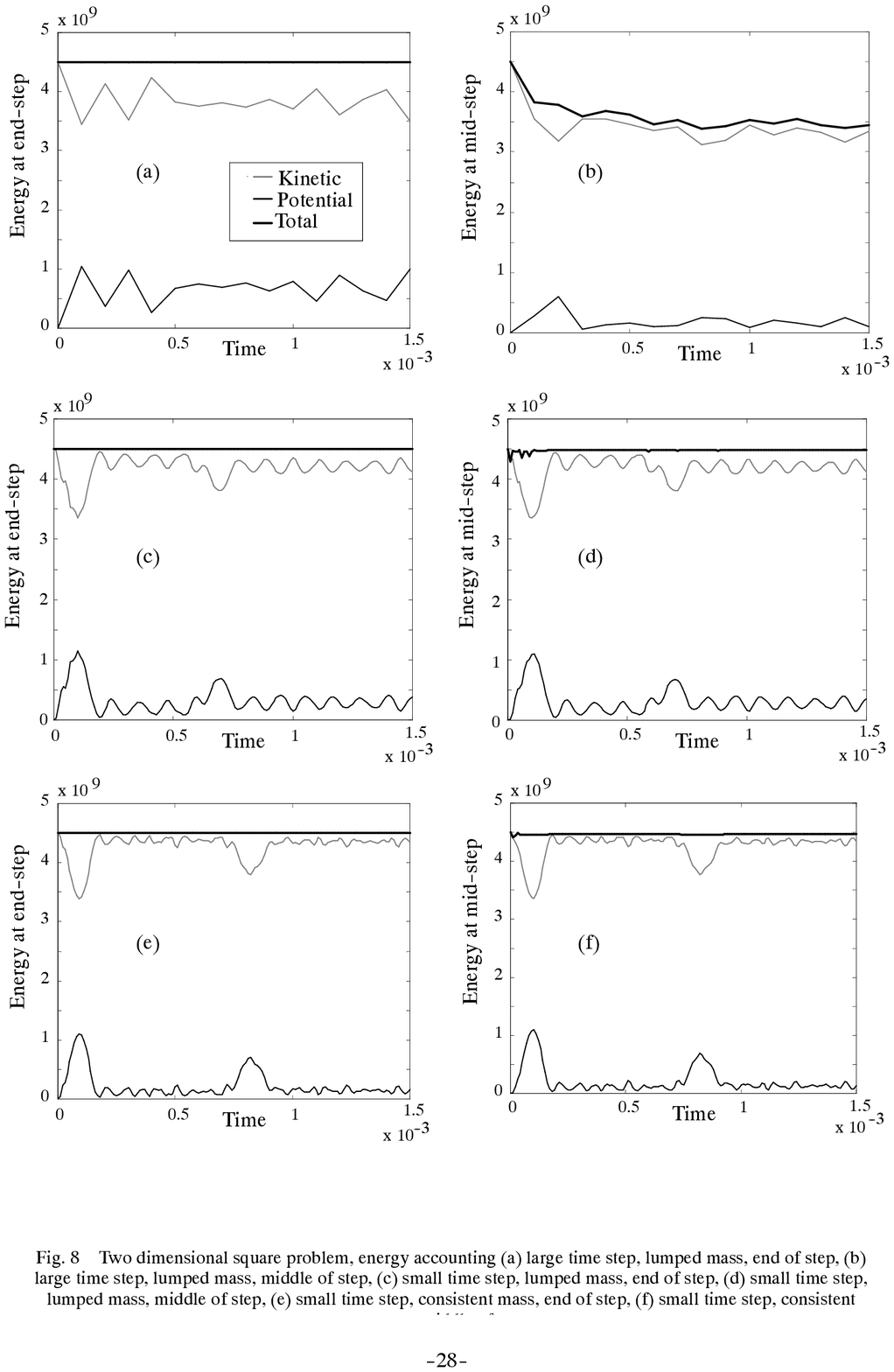}
	\caption{Two dimensional square problem, energy accounting (a) large time step, lumped mass, end of step, (b) large time step, lumped mass, middle of step, (c) small time step, lumped mass, end of step, (d) small time step, lumped mass, middle of step, (e) small time step, consistent mass, end of step, (f) small time step, consistent mass, middle of step  \label{fig:con_E_Q4}}
	\end{figure}

Figure \ref{fig:con_E_Q4} (a) shows the evolution of the system energy with time, where the potential, kinetic and total energy are calculated at the end of each time step. It is clear that, in spite of relatively large changes in both the kinetic and the potential energy, an increase in the potential energy is balanced by an equivalent decrease the kinetic energy, and vice versa, so that the total energy of the system at the end of each time step is exactly conserved. The exchange between the potential and kinetic energy becomes clearer when carrying out the solution using a smaller time step of $\Delta t = 10^{- 5}$, as shown in Figure \ref{fig:con_E_Q4} (c). It can then be observed
that, aside from a few smaller vibrations between individual contact events, the increase in potential energy corresponds to contact events between two or more cubes.
This increase is due to the high deformations induced by contact and leads to an equal decrease in the kinetic energy of the system. Consequently, the contacting cubes slow down while pulling apart from each other, until full separation occurs and the cubes regain full speed while moving as rigid bodies.

Figures \ref{fig:con_E_Q4} (b) and (d) show the energy distributions at midstep, obtained using $\Delta t = 7 \times 10^{- 5}$, and $\Delta t = 10^{- 5}$, respectively. It can be observed from Figure \ref{fig:con_E_Q4} (b), that for relatively large steps, energy conservation is not achieved at midstep. In particular, the solution is highly dissipative. This
phenomenon can be explained by the fact that the Lagrange multiplier guarantees conservation at the end of the time step specifically. The added force vector has
the effect of generating corrective accelerations, and therefore corrective velocities, to balance the energy at the end of the step. Since the energy conservation error (due to the work produced by contact forces) at mid-step is expected to be less than that at the end of the step, these algorithmic corrections introduce some unbalance in total energy at mid-step. In any case, for a smaller time step, the amount of numerical dissipation is reduced and energy is conserved at mid-step throughout the analysis, as shown in Figure \ref{fig:con_E_Q4} (d). Figures \ref{fig:con_E_Q4} (e) and (f) display the results obtained using a consistent mass matrix, in which case the kinetic vs. potential energy exchanges occur mainly around contact events, and the system vibrations between these events are reduced.

\begin{figure}
	\centering
	\includegraphics[clip]{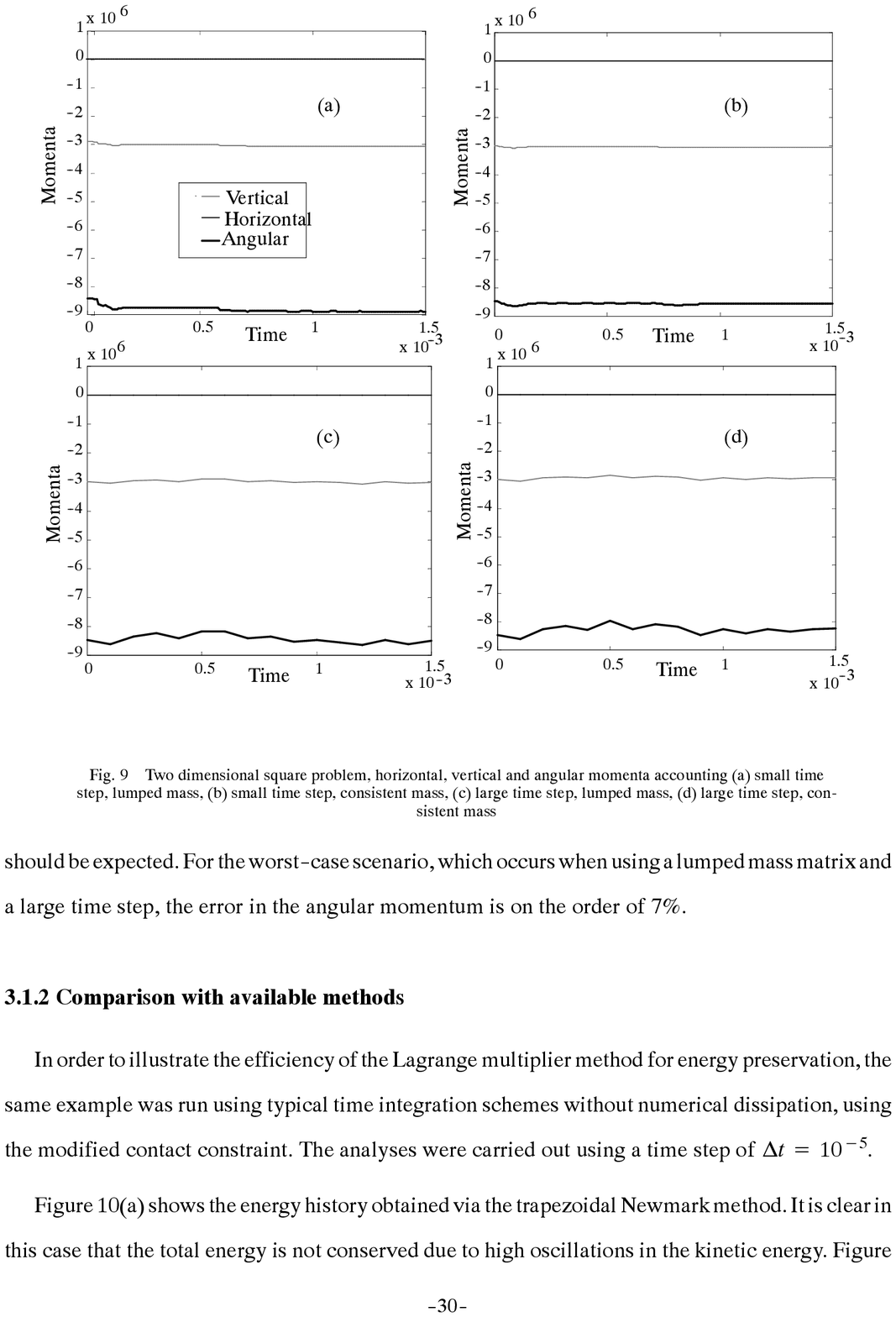}
	\caption{Two dimensional square problem, horizontal, vertical and angular momenta accounting (a) small time step, lumped mass, (b) small time step, consistent mass, (c) large time step, lumped mass, (d) large time step, consistent mass  \label{fig:con_M_Q4}}
	\end{figure}

Figures \ref{fig:con_M_Q4} (a)-(d) depict the evolution of the linear and angular momenta during the motion, for different choices of the time step size and mass matrix formulation. It can
be clearly seen that the momenta, and in particular the angular momentum, are not exactly conserved. The error can be attributed to the added force vector produced by
the energy conservation Lagrange multiplier. This force vector is purely algorithmic and does not represent any physical applied force on the system. The results show better conservation for smaller time step and a consistent mass matrix, as should be expected. For the worst-case scenario the error in the angular momentum is on the order of 4 
%

%
%
	\begin{figure}
	\centering
	\includegraphics[clip]{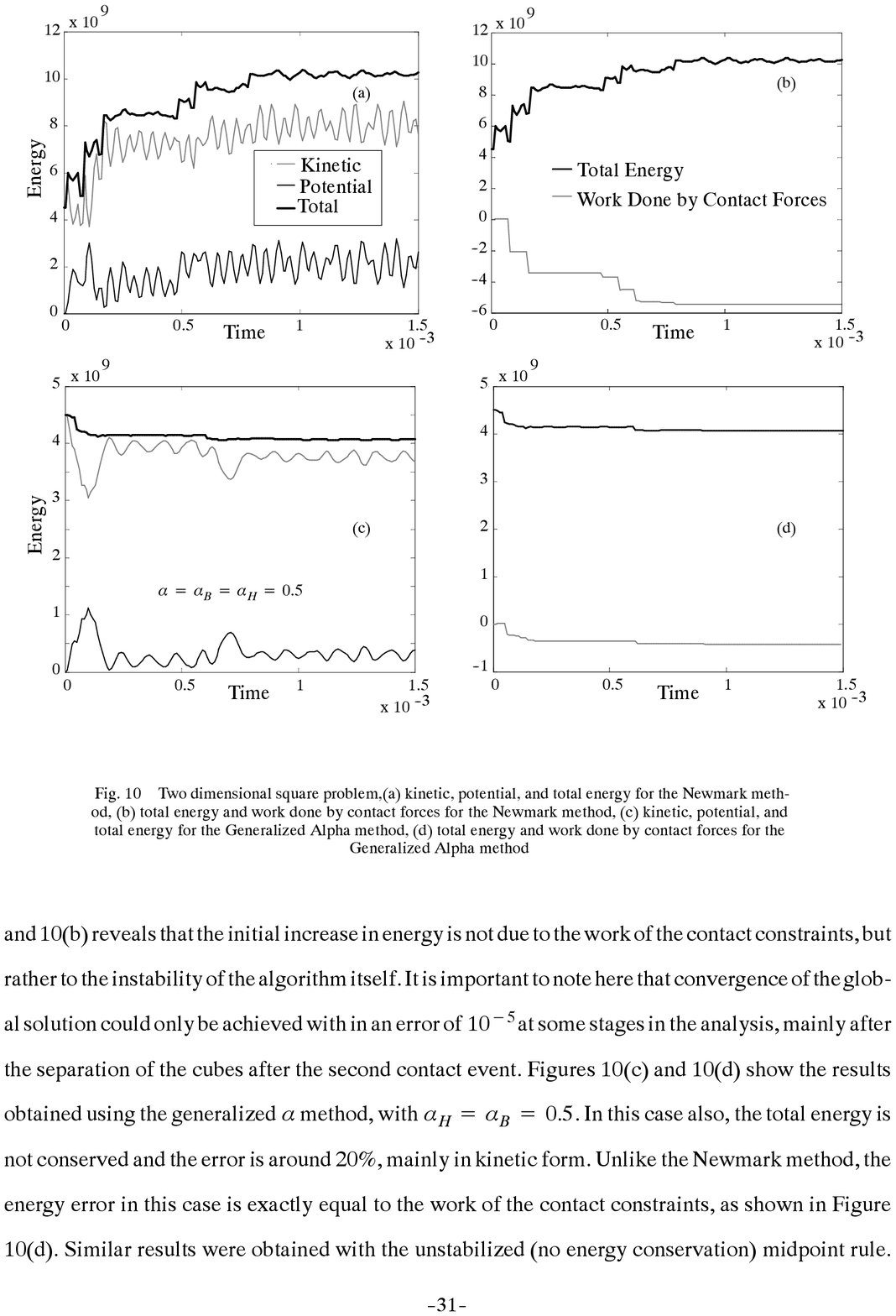}
	\caption{Two dimensional square problem,(a) kinetic, potential, and total energy for the Newmark method, (b) total energy and work done by contact forces for the Newmark method, (c) kinetic, potential, and total energy for the Generalized Alpha method, (d) total energy and work done by contact forces for the Generalized Alpha method  \label{fig:con_comp}}
	\end{figure}
\subsubsection{Comparison with available methods} \label{sec_con:Comp}

In order to illustrate the efficiency of the Lagrange multiplier method for energy preservation, the same example was run using typical time integration schemes without
numerical dissipation, using the modified contact constraint. The analyses were carried out using a time step of $\Delta t = 10^{- 5}$.

Figure \ref{fig:con_comp} (a) shows the energy history obtained via the trapezoidal Newmark method. It is clear in this case that the total energy is not conserved due to high oscillations in the kinetic energy. Figure \ref{fig:con_comp} (b) depicts the history of work performed by the contact constraints. Comparison of Figures \ref{fig:con_comp} (a) and (b) reveals that the initial increase in energy is not due to the work of the contact constraints, but rather to the instability of the algorithm itself. It is important to note here that convergence of the global solution could only be achieved with in an error of $10^{- 5}$ at some stages in the analysis, mainly after the separation of the cubes after the second
contact event. Figures \ref{fig:con_comp} (c) and (d) show the results obtained using the generalized a method, with $\alpha_H = \alpha_B = 0.5$. In this case also, the total energy is not conserved and the error is around 12\%, mainly in kinetic form. Unlike the Newmark method, the energy error in this case is exactly equal to the work of the contact constraints, as shown in Figure \ref{fig:con_comp} (d). Similar results were obtained using the midpoint rule without the energy conservation constraint.

These results justify using the Lagrange multiplier approach to guarantee the conservation of energy, even though some error is introduced in
the conservation of momenta. For a larger time step, the error in momenta conservation in the Lagrange multiplier method is on the order of 4\%, whereas the solution using the aforementioned standard methods displays a large energy error and the solution may diverge altogether. When a smaller time step is used, momenta can be conserved within a 2\% error with the Lagrange multiplier approach, as opposed to a 12\% energy error without it.
%
%
%
	\begin{figure}
	\centering
	\includegraphics[clip]{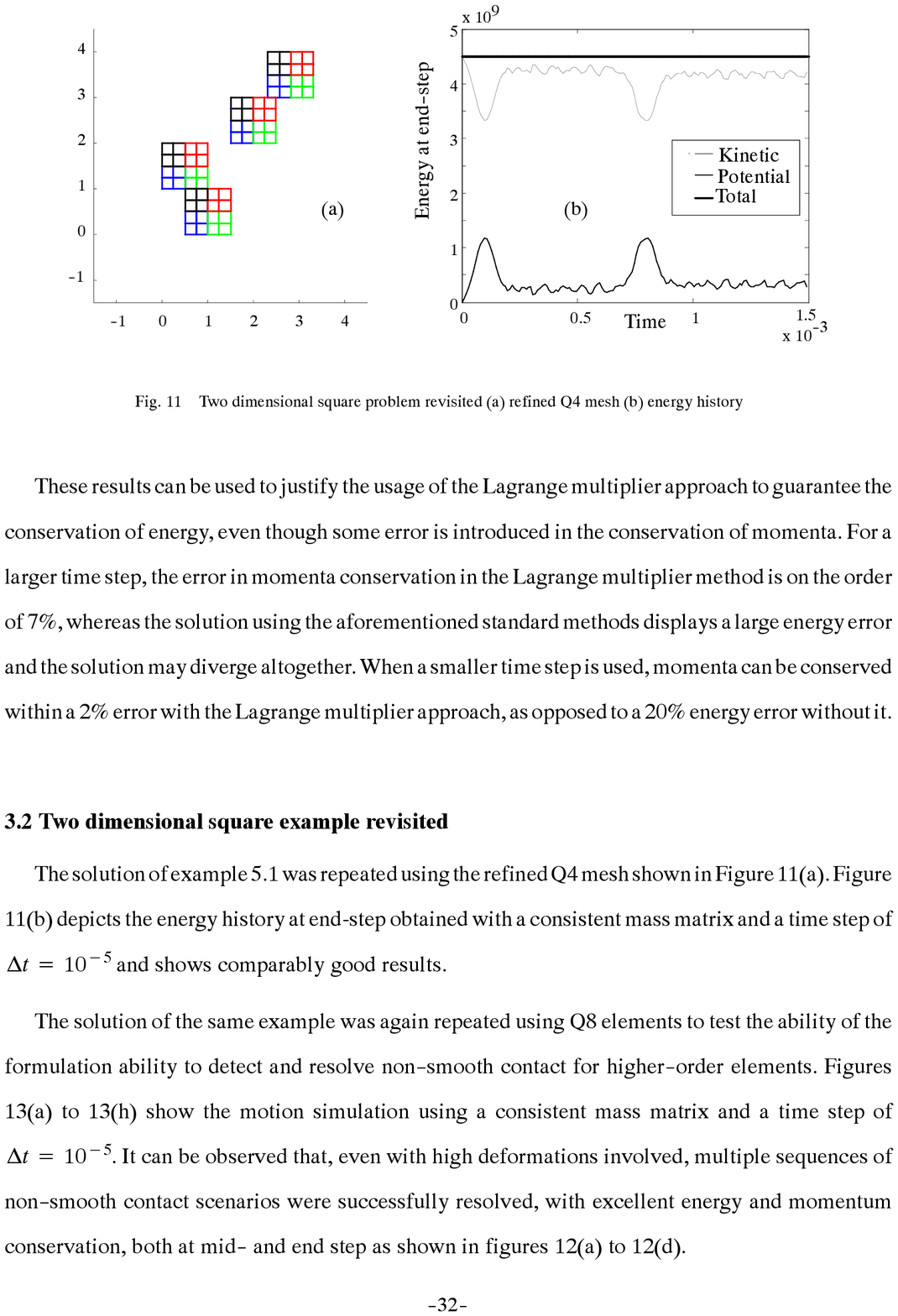}
	\caption{Two dimensional square problem revisited (a) refined Q4 mesh (b) energy history  \label{fig:con_ref_Q4}}
	\end{figure}
\subsection{Two-dimensional square example revisited} \label{sec_con:2Dex_rev}
\begin{figure}
	\centering
	\includegraphics[clip]{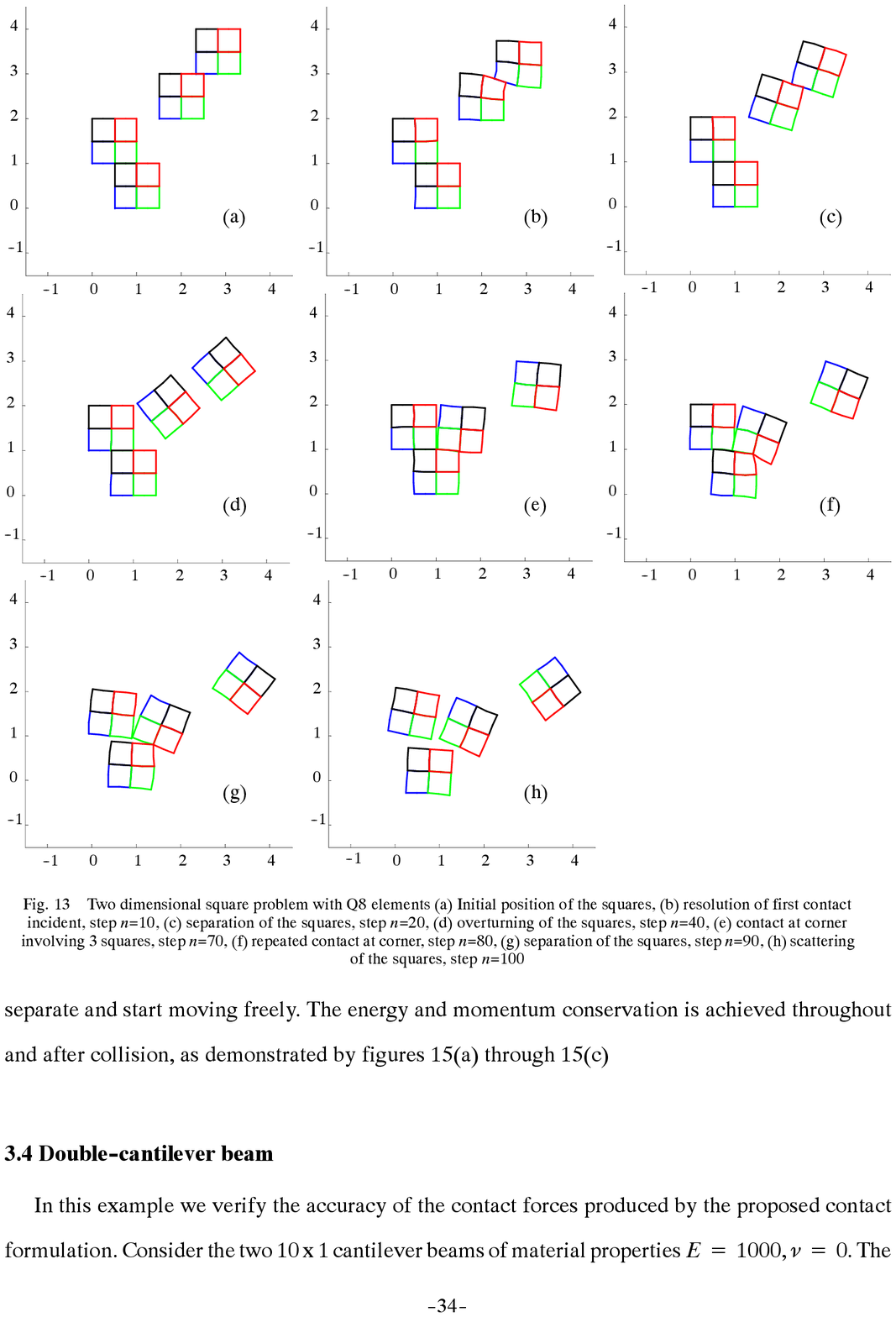}
	\caption{Two dimensional square problem with Q8 elements (a) Initial position of the squares, (b) resolution of first contact incident, step $n$=10, (c) separation of the squares, step $n$=20, (d) overturning of the squares, step $n$=40, (e) contact at corner involving 3 squares, step $n$=70, (f) repeated contact at corner, step $n$=80, (g) separation of the squares, step $n$=90, (h) scattering of the squares, step $n$=100  \label{fig:con_motion_Q8}}
	\end{figure}

\begin{figure}
	\centering
	\includegraphics[clip]{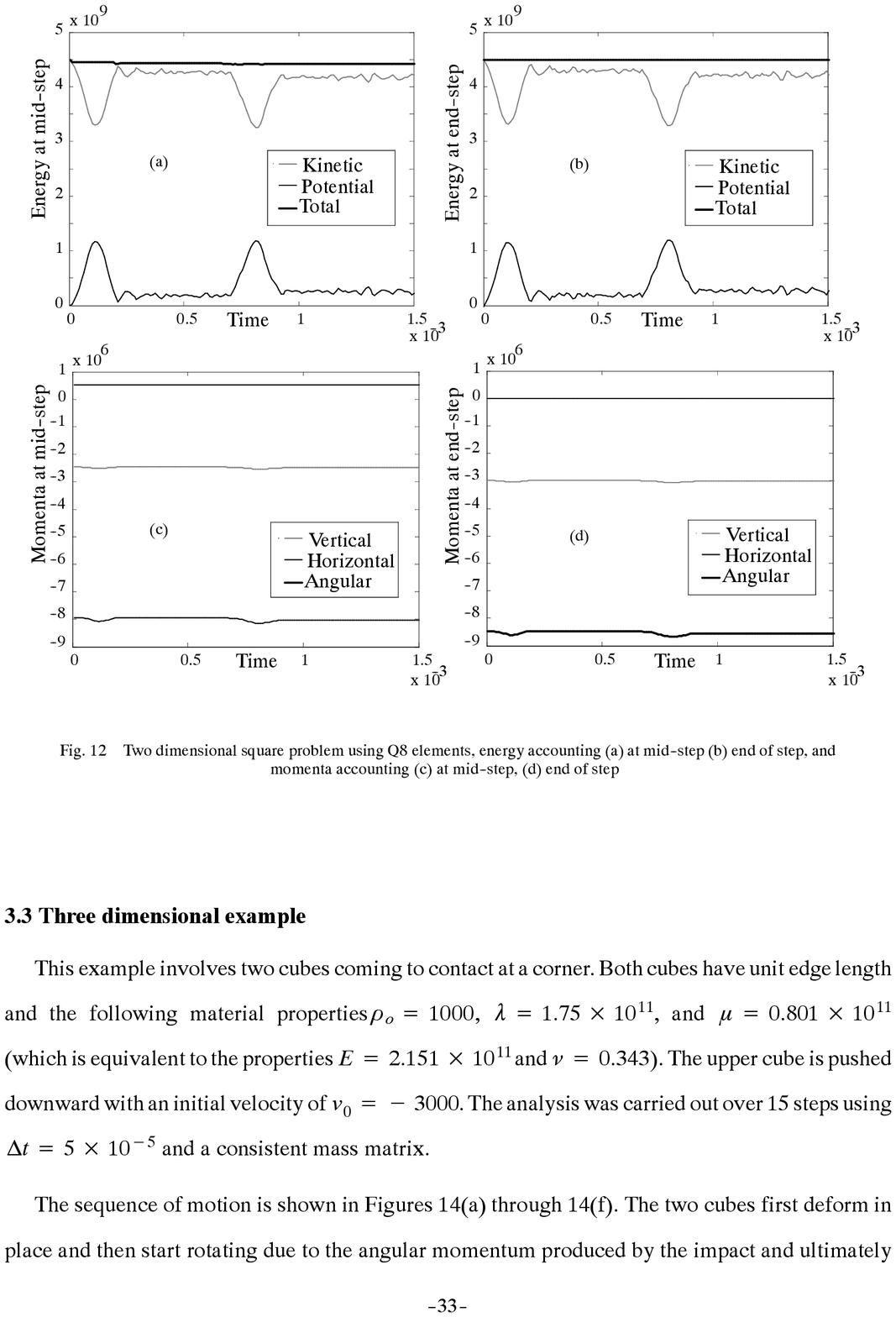}
	\caption{Two dimensional square problem using Q8 elements, energy accounting (a) at mid-step (b) end of step, and momenta accounting (c) at mid-step, (d) end of step  \label{fig:con_EM_Q8}}
	\end{figure}

The solution of example 5.1 was repeated using the refined Q4 mesh shown in Figure \ref{fig:con_ref_Q4} (a). Figure \ref{fig:con_ref_Q4} (b) depicts the energy history at end-step obtained with a consistent mass
matrix and a time step of $\Delta t = 10^{- 5}$ and shows comparably good results. The solution of the same example was again repeated using Q8 elements to test the ability of the formulation ability to detect and resolve non-smooth contact for higher-order elements. Figures \ref{fig:con_motion_Q8} (a)-(h) show the motion simulation using a consistent mass matrix and a time step of $\Delta t = 10^{- 5}$. It can be observed that, even with high deformations involved, multiple sequences of non-smooth contact
scenarios were successfully resolved, with excellent energy and momentum conservation, both at mid- and end-step as shown in Figures \ref{fig:con_EM_Q8} (a)-(d). The solution in this case matches the result obtained using the refined Q4 mesh.
%
%
%
%
	\begin{figure}
	\centering
	\includegraphics[clip]{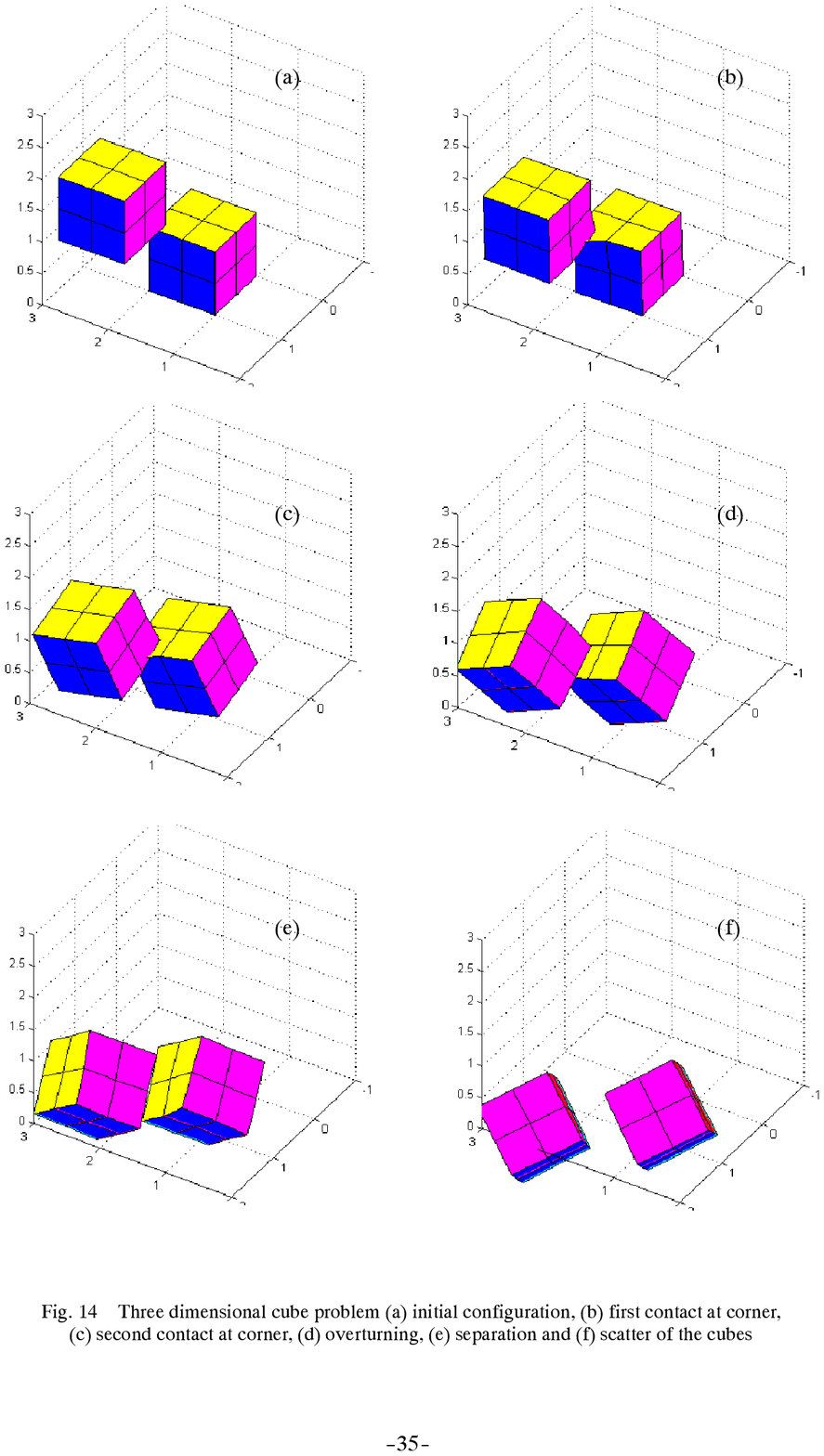}
	\caption{Three dimensional cube problem (a) initial configuration, (b) first contact at corner, (c) second contact at corner, (d) overturning, (e) separation and (f) scatter of the cubes  \label{fig:con_motion_3D}}
	\end{figure}
\subsection{Three-dimensional example} \label{sec_con:3Dex}
This example involves two cubes coming to contact at a corner. Both cubes have unit edge length and the following material properties $\rho_0 = 1000$, $\lambda = 1.75 \times 10^{11}$, and $\mu = 0.801 \times 10^{11}$ (which is equivalent to the properties $E = 2.151 \times 10^{11}$ and $\nu = 0.343$). The upper cube is pushed downward with an initial velocity of $v_0 = - 3000$. The analysis was carried out over 15 steps using $\Delta t = 5 \times 10^{- 5}$ and a consistent mass matrix.

\begin{figure}
	\centering
	\includegraphics[clip]{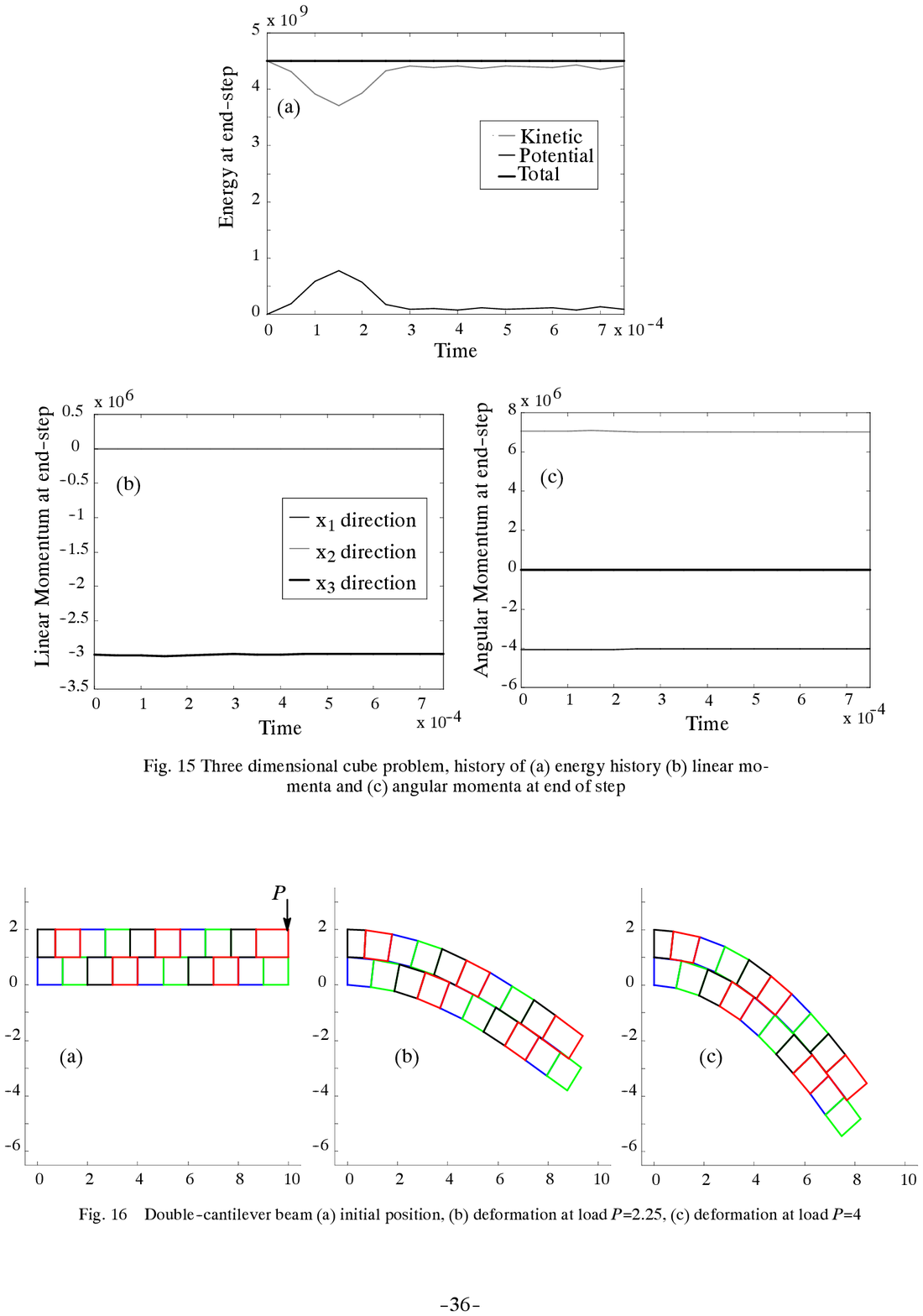}
	\caption{Three dimensional cube problem, history of (a) energy (b) linear momenta and (c) angular momenta at end of step  \label{fig:con_EM_3D}}
	\end{figure}

The sequence of motion is shown in Figures \ref{fig:con_motion_3D} (a) through (f). The two cubes first deform in place and then start rotating due to the angular momentum produced by the impact and ultimately separate and start moving freely. The energy and momentum conservation is achieved throughout and after collision, as demonstrated by Figures \ref{fig:con_EM_3D} (a) through (c).
%
%
%
%
\subsection{Double-cantilever beam} \label{sec_con:DblCantBeam}
In this example we verify the accuracy of the contact forces produced by the proposed contact formulation. Consider the two $10 \times 1$ cantilever beams with material
properties $E = 1000$, $\nu = 0$. The beams are superposed on top of each other and subjected to a quasi-static force $P$ at the right end. We discretize the two beams with non-conforming meshes of QM6 elements to avoid shear locking.

\begin{figure}
	\centering
	\includegraphics[clip]{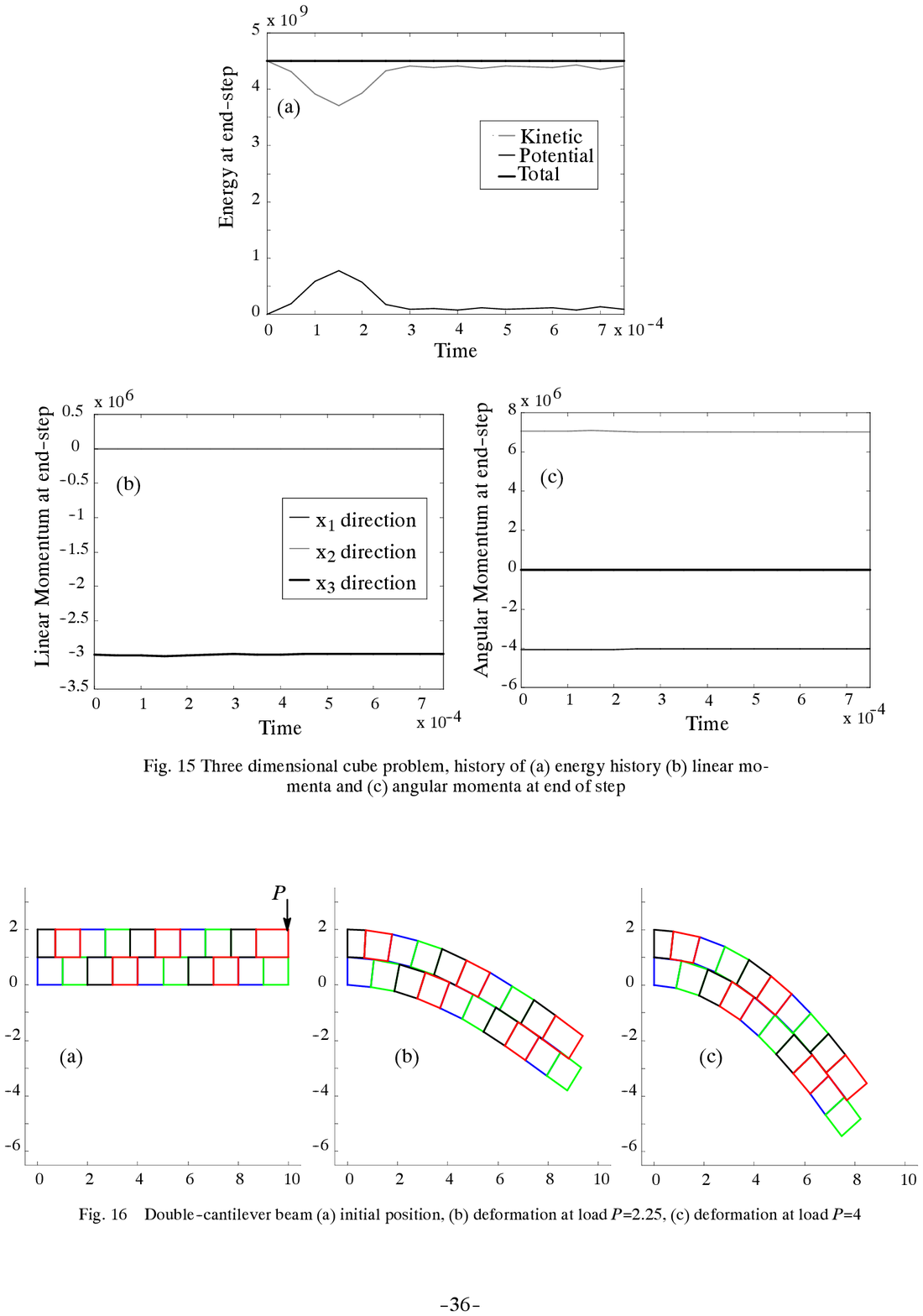}
	\caption{Double-cantilever beam (a) initial position, (b) deformation at load P=2.25, (c) deformation at load P=4  \label{fig:con_cant_def}}
	\end{figure}

The deformation sequence is shown in Figure \ref{fig:con_cant_def}. Figure \ref{fig:con_cant_def} shows the load-deformation curve of the double-cantilever beam compared to that of a single cantilever with $E = 2000$, where the modulus of elasticity was doubled to simulate the presence of two separate beams. The two curves match within a reasonable margin of error; furthermore, in the small deformation range, the tip displacement of the double-cantilever is exactly equal to that obtained using traditional beam theory.
	\begin{figure}
	\centering
	\includegraphics[clip]{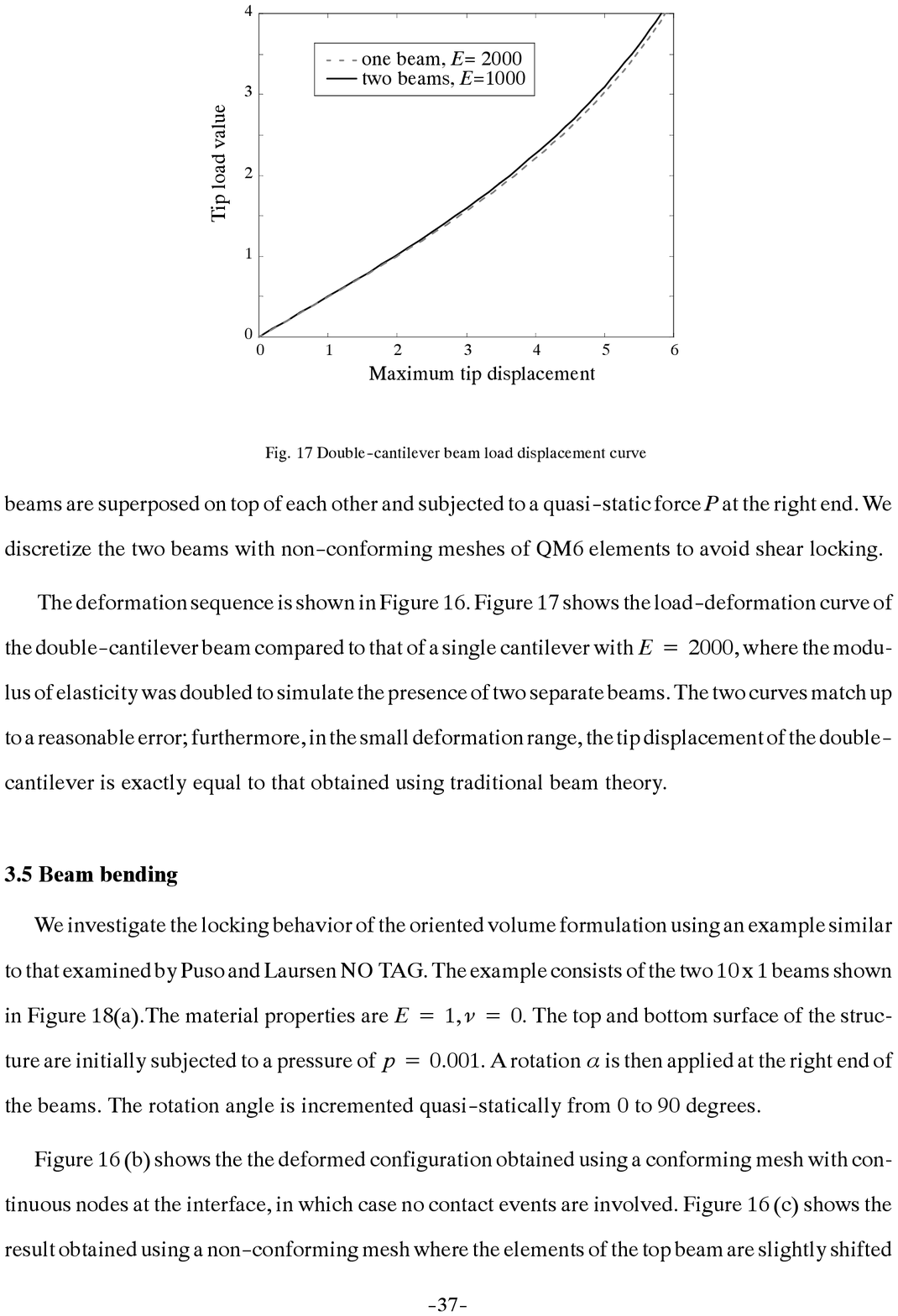}
	\caption{Double-cantilever beam load displacement curve  \label{fig:con_cant_res}}
	\end{figure}
%
%
%
%
%
\section{Conclusions} \label{sec_con:Concl}
We have presented a new approach for the formulation of non-smooth contact. The suggested approach, based on the calculation of an oriented penetration volume, proved to be very efficient in dealing with highly non-smooth contact scenarios such as at corners. In fact, the contact constraint reduces to an oriented gap function, except that the orientation of the gap is determined by the location of the penetrating node inside the penetrated element, as opposed to a dot product with a normal vector, a characteristic that enables it to deal with non-smooth contact locations. Accordingly, the contact constraint retains the advantages of the gap function formulation, such as its simple implementation, while overcoming its main limitation. 

The implementation using the midpoint rule showed the need to conserve total energy during impact. This objective can be achieved provided that the work done by the contact forces, expressed in terms of the gradient of the contact constraints, vanishes exactly over the time step. This condition applies to static and persistent contact, but cannot be met in dynamic impact that results in the separation or scattering of the contacting bodies. In this case, non-conservation of energy manifests in the form of spurious oscillations after contact, since the additional work input due to the contact forces transmits into kinetic energy. A Lagrange multiplier method proved to be a successful approach towards the exact satisfaction of energy conservation, while preserving the linear and angular momenta of the system up to a reasonable error. The role of the Lagrange multiplier is to introduce algorithmic forces/accelerations that would produce the velocity corrections necessary to conserve energy during non-smooth events occurring during a given time step. This approach proved to exactly conserve energy at the end of the time step while introducing some dissipation at mid-step and some error in the conservation of momenta. The amount of energy dissipation and momentum error was shown to decrease
with spatial and temporal mesh refinement. When conservation of energy is achieved through the algorithmic persistency condition, the energy Lagrange multiplier results to be zero. This holds for static/persistent contact and also for impact in the limit of temporal refinement. This result implies that the Lagrange multiplier method preserves the algorithmic consistency of the underlying formulation.

The analysis showed to be consistently stable for both relatively large and relatively small spatial and time discretization, with an improved performance with refinement. The suggested formulation is readily applicable to higher-order 2D 4-node quadrilateral elements and to 3D elements. This formulation would also be easily extended to triangular elements.
%
%
%







%% file: motivation.tex
\chapter {Interface stabilization: motivation}

	\begin{figure}
	\centering
	\includegraphics[clip]{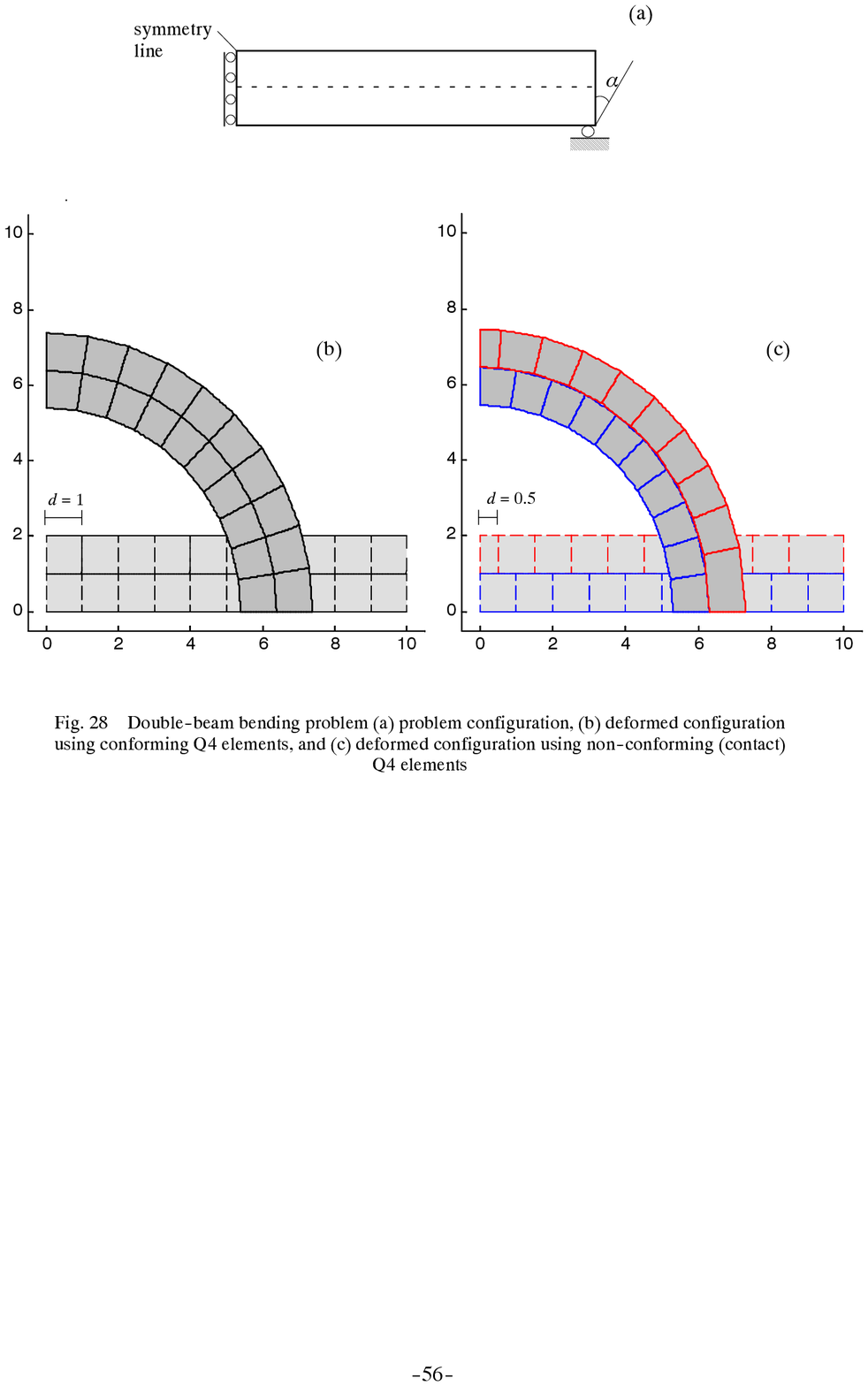}
	\caption{Double-beam bending problem with $p=0.001$ (a) problem configuration, (b) deformed configuration using conforming Q4 elements, and (c) deformed configuration using non-conforming (contact) Q4 elements  \label{fig:mot_beam}}
	\end{figure}
	\begin{figure}
	\centering
	\includegraphics[clip]{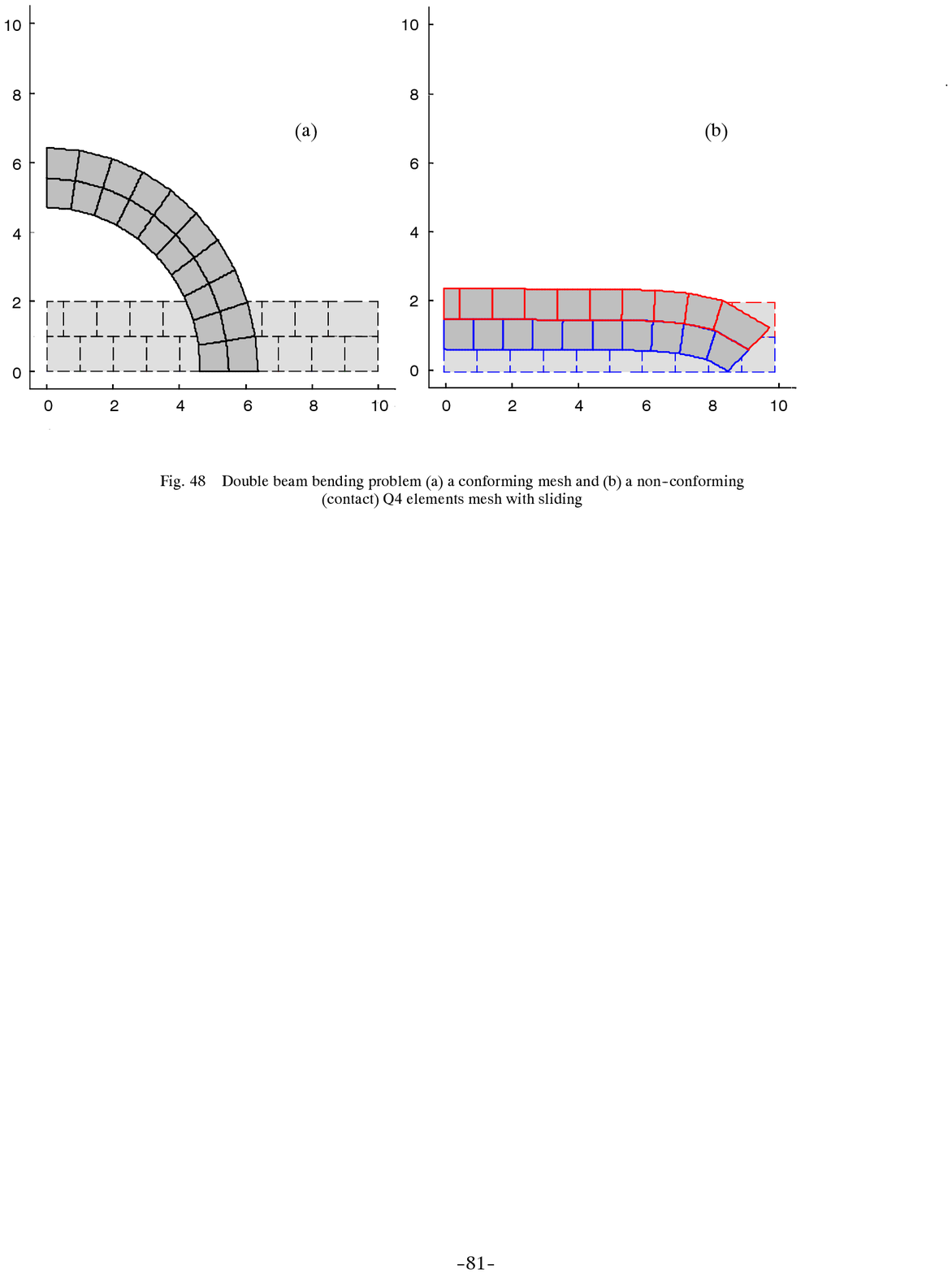}
	\caption{Double-beam bending problem with $p=0.1$ (a) deformed configuration using conforming Q4 elements, and (b) deformed configuration using non-conforming (contact) Q4 elements with sliding \label{fig:mot_sliding}}
	\end{figure}
We investigate the locking behavior of the oriented volume formulation using an example similar to that examined by Puso and Laursen \cite{pusoCMAME04} for the mortar method. The problem consists of two beams discretized using $10 \times 1$ quadrilateral elements, as shown in Figure \ref{fig:mot_beam}(a). The beams are made of a Neo-Hookean material with strain density function given by 
\begin{equation}
\psi\left(\mathbf{u}\right) = \frac{1}{2}\kappa\left[\frac{1}{2}(J^2-1)-\mathrm{log}J\right]-\frac{1}{2}\mu\,\left[ J^{-2/3}\mathrm{tr}\left(\mathbf{C}\right)\right],
\end{equation}
where $\kappa = E/3(1-2\nu)$ and $\mu = E/2(1+\nu)$ are the bulk and shear moduli, respectively \cite{hjelmstad2005}. The material properties are $E = 1$, $\nu = 0$. The top and bottom surface of the structure are initially subjected to a pressure of $p = 0.001$. A rotation $\alpha$ is then applied at the right end of the beams.
The rotation angle is incremented quasi-statically from $0^o$ to $90^o$. Since the quasi-static loading precludes any dynamic effects, the energy Lagrange multiplier was not used in this example.

Figure \ref{fig:mot_beam}(b) shows the deformed configuration obtained using a conforming mesh with continuous nodes at the interface, in which case no contact events are involved. Figure \ref{fig:mot_beam}(c) shows the result obtained using a non-conforming mesh where the elements of the top beam are slightly shifted to the left and the nodes do not match at the interface. This configuration clearly involves multiple contact events and the applied pressure forces all contact constraints to be activated before rotation is applied. The final deformed shape clearly matches that of the conforming mesh.

This case was reported to suffer locking when a traditional two-pass node-to-surface contact formulation is used \cite{pusoCMAME04}. No locking is observed in this case and the analysis was successfully carried out until the end. This result can be attributed to the observation that, when the rotation is applied, some nodes that were in contact in the initial configuration (pressure only) separate and are therefore removed from the contact active set. Note that the pressure value used in this case is smaller than that reported by Puso and Laursen \cite{pusoCMAME04}. 

%
	\begin{figure}
	\centering
	\includegraphics[clip]{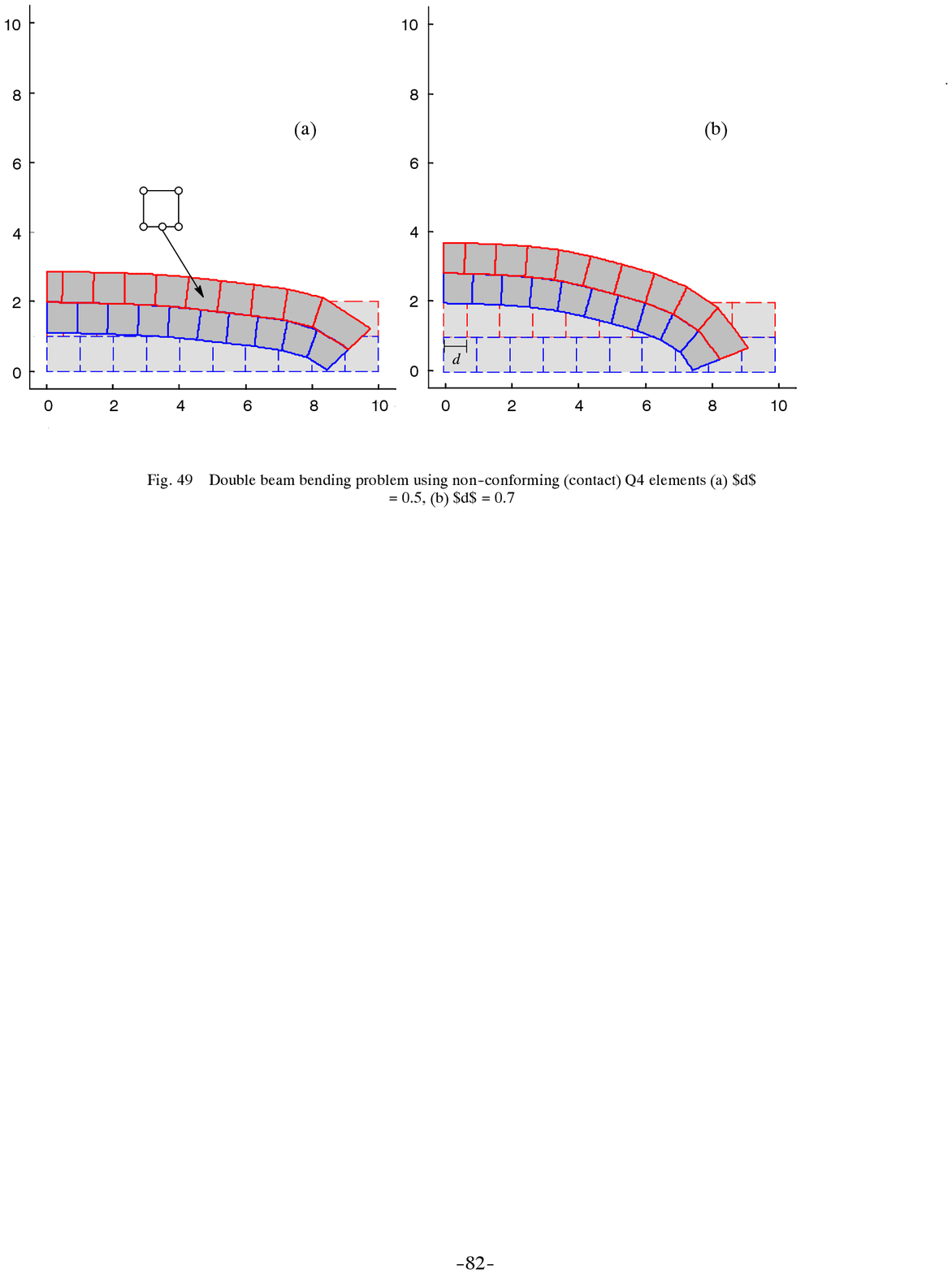}
	\caption{Double beam bending problem using non-conforming (contact) Q4 elements with no sliding (a) d = 0.5, (b) d = 0.7  \label{fig:mot_lock}}
	\end{figure}
Now consider the case $p = 0.1$. The deformed configuration is compared in Figure \ref{fig:mot_sliding} (b) to the result obtained using a conforming mesh shown in Figure \ref{fig:mot_sliding} (a). At this pressure level all contact events at the interface get activated and the formulation shows severe locking. Moreover, the nodes of the top beam slide substantially along the surface of the bottom one, which leads the whole system to become unstable. This behavior is not a feature of the original symmetric problem and is due to the geometric nonlinearity of the problem, as the sum of the stiffnesses of the two beams is substantially lower than that of the original structure. 

To eliminate the instability due to sliding, we tie the middle of the lower beam to the surface of the top one by introducing a node at that location and assembling the top element over 5 nodes (Q5). Figure \ref{fig:mot_lock} shows the results obtained for Q4 meshes with different values of $d$, where $d$ measures the shift of the top beam mesh with respect to the lower one. We can observe that all the nodes along the interface remain in contact with the surface during bending, and locking is observed before the rotation
angle reaches 90 degrees, as shown in Figure \ref{fig:mot_lock}. It is interesting to note that, as the motion progresses, the interface nodes shift until the node-to-surface contact becomes node-to-node. Also, the mesh with a smaller value of $d$ locks earlier.

\begin{figure}
	\centering
	\includegraphics[clip]{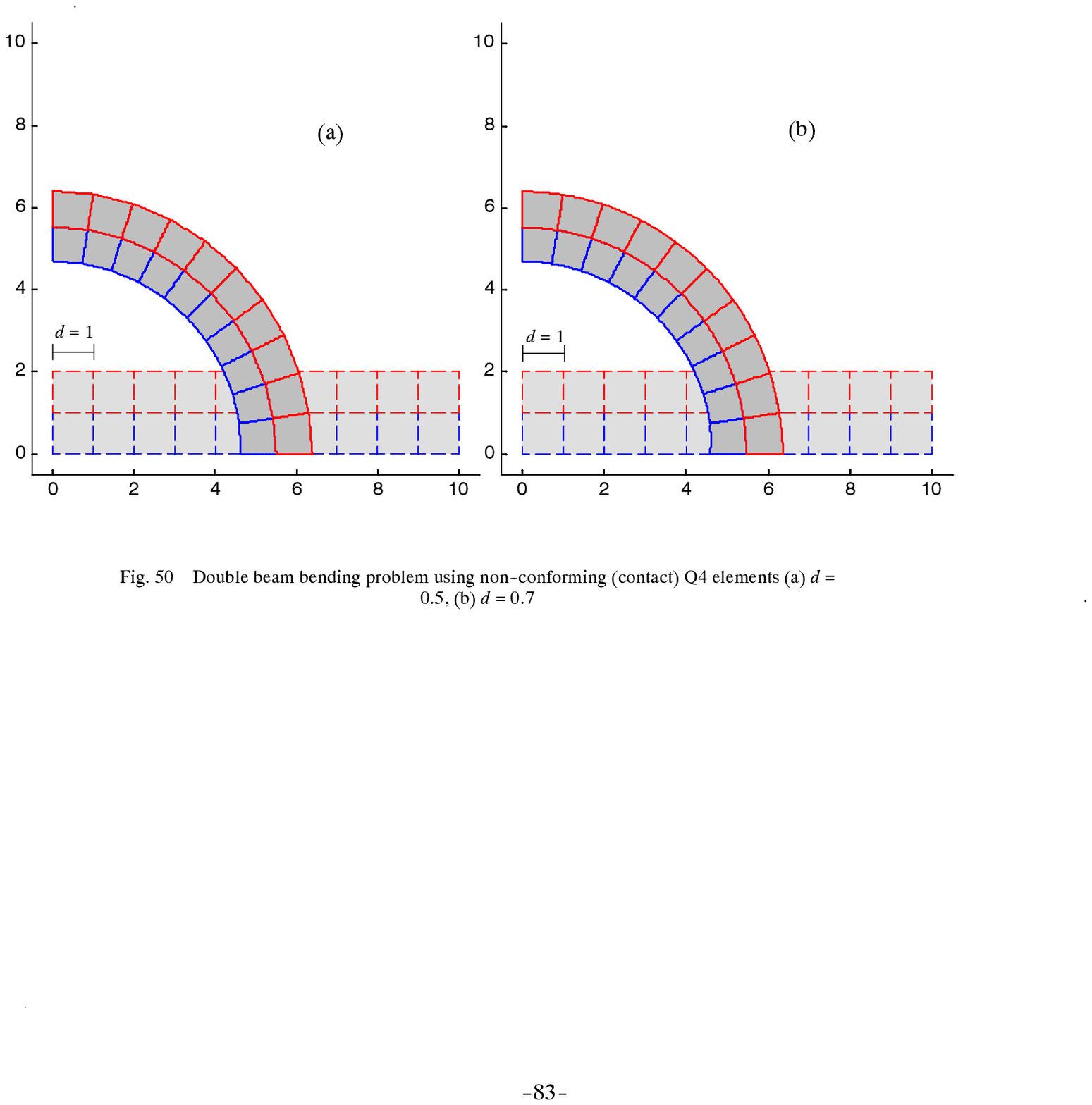}
	\caption{Double-beam bending problem, deformed configuration using (a) matching (contact) Q4 elements along, and (b) matching (contact) Q8 elements  \label{fig:mot_node}}
	\end{figure}

Now let $d=1$. For this special configuration the node-to-surface gap function reduces to a node-to-node constraint. Figure \ref{fig:mot_node} (a) shows the result obtained with a mesh of Q4 elements. Unlike the non-matching mesh of Figure \ref{fig:mot_lock}, no surface locking is observed in this case. A similar result is obtained using a mesh of Q8 elements, as depicted in Figure \ref{fig:mot_node} (b). These results suggest that surface locking can be avoided if contact can be guaranteed to happen between two nodes.

\begin{figure}
	\centering
	\includegraphics[clip]{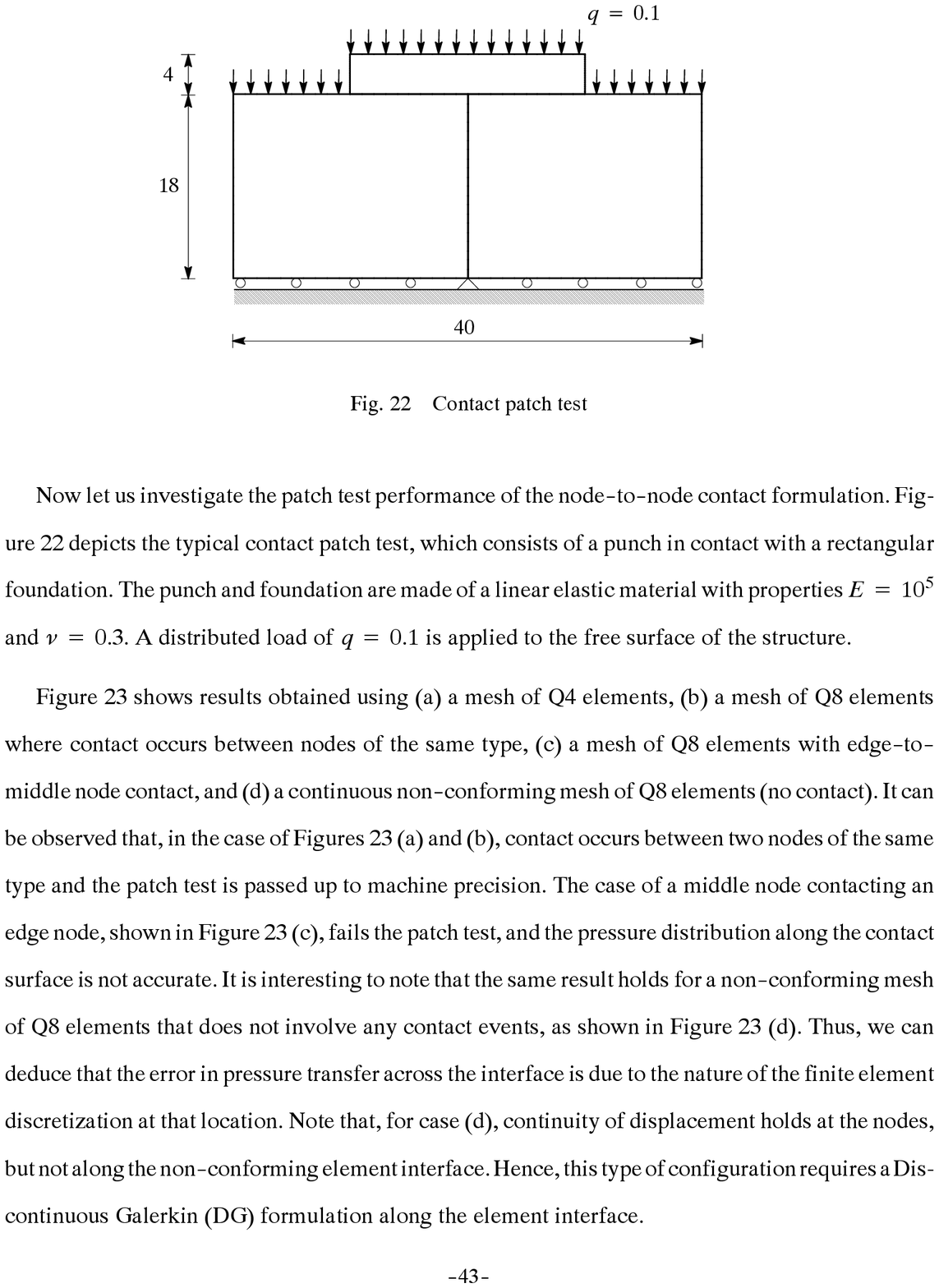}
	\caption{The contact patch test  \label{fig:mot_patch}}
	\end{figure}

Next we investigate the patch test performance of the node-to-node contact formulation. Figure \ref{fig:mot_patch} depicts the typical contact patch test, which consists of a punch in contact with a rectangular foundation. The punch and foundation are made of a linear elastic material with properties $E = 10^5$ and $\nu = 0.3$. A distributed load of $q = 0.1$ is applied to the free surface of the structure.
	\begin{figure}
	\centering
	\includegraphics[clip]{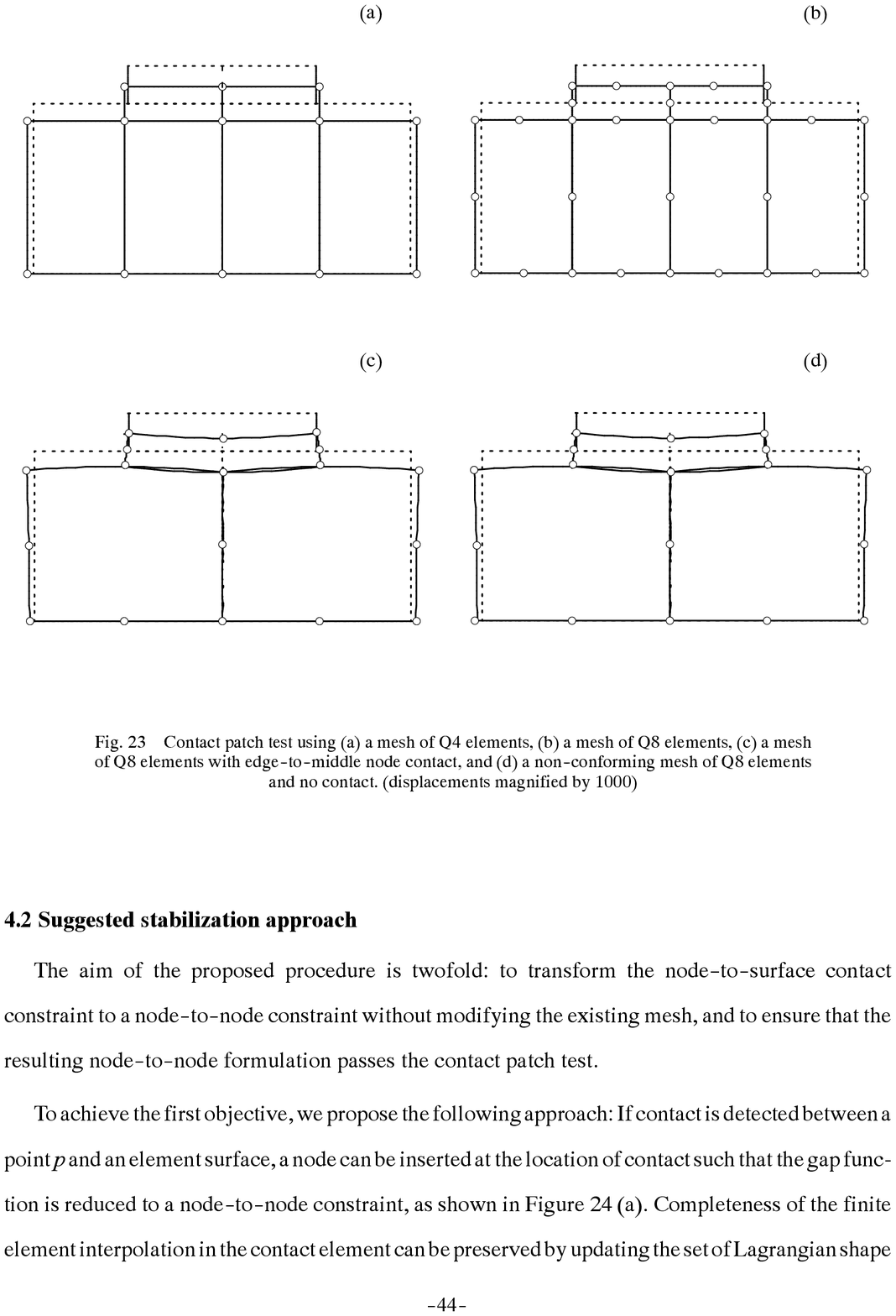}
	\caption{Contact patch test using (a) a mesh of Q4 elements, (b) a mesh of Q8 elements, (c) a mesh of Q8 elements with edge-to-middle node contact, and (d) a non-conforming mesh of Q8 elements and no contact. (displacements magnified by 1000)  \label{fig:mot_patchresults}}
	\end{figure}

Figure \ref{fig:mot_patchresults} shows results obtained using (a) a mesh of Q4 elements, (b) a mesh of Q8 elements where contact occurs between nodes of the same type, (c) a mesh of Q8 elements with edge-to-middle node contact, and (d) a continuous non-conforming mesh of Q8 elements (no contact). It can be observed that, in the case of Figures \ref{fig:mot_patchresults} (a) and (b), contact occurs between two nodes of the same type and the patch test is passed up to machine precision. The case of a middle node contacting an edge node, shown in Figure \ref{fig:mot_patchresults} (c), fails the patch test, and the pressure distribution along the contact surface is not accurate. It is interesting to note that the same result holds for a non-conforming mesh of Q8 elements that does not involve any contact events, as shown in Figure \ref{fig:mot_patchresults} (d). Thus, we can deduce that the error in pressure transfer across the interface is due to the nature of the finite element discretization at that location. Note that, for cases (c) and (d), continuity of displacement holds at the nodes but not along the inter-element interfaces.

These observations provide the background and motivation for the stabilized interface formulation proposed next. As mentioned in Chapter 1 and evidenced by the above example, the issues of locking and patch test performance are common to contact problems and to the coupling of non-conforming meshes. Therefore, to simplify the presentation of the proposed stabilized interface formulation, we start by applying it to the coupling of non-conforming meshes in static hyperelasticity, as described in the following chapter. In Chapter 6 we extend the formulation to the treatment of contact interfaces in a quasi-static framework.

%% file: DG.tex

%
%
%
%
%
%
\chapter {A stabilized interface formulation for the coupling of non-conforming meshes in hyperelasticity}
%
%
%
%
\section{Introduction}   \label{sec_DG:intro}


Non-conforming meshes are typically used to allow for the partial adaptive refinement of a finite
element mesh, such as at locations of stress concentration or steep gradients, without imposing a high
computational cost in locations where such a refinement is not needed. The main challenge in
coupling non-conforming meshes is to capture interface effects accurately. This includes enforcing
continuity of the displacement field while ensuring a complete transfer of tractions across the non-conforming interface. Different methods have been historically proposed to achieve this purpose.

Coupling methods can be classified into \textit{primal} and \textit{dual} methods. The distinction between these two families of methods is in the nature of the interface variables, a primal field such as displacement for the former, and a dual traction field in the latter. For the elliptic problem of hyperelasticity, the mathematical framework for primal methods is a minimization problem, while dual formulations typically result in a saddle-point extremization problem. An overview of the most prominent members of these two classes of methods is presented next. 

\subsection{Dual methods} \label{sec_DG:intro_dual}

Most currently available coupling methods are dual approaches that employ a field of Lagrange multipliers to enforce geometric compatibility at the interface. The choice of the Lagrange multiplier field is not trivial since not all combinations of primal/dual variable discretizations satisfy the \textit{inf-sup} or Ladyzhenskaya-Babu\u{s}ka-Brezzi (LBB) condition, a requirement to achieve spatial convergence of the saddle-point problem \cite{brezziSV01}. In particular, the discrete Multi-Point-Constraint method of enforcing displacement continuity between a set of \textit{slave} nodes at one side of the interface and the \textit{master} surface on the other using a set of discrete Lagrange multipliers is generally not capable of representing a state of constant pressure and therefore fails the well-known patch test \cite{papadopoulosCMAME92}. Conceptually similar to the contact gap function, this approach has been shown to present various numerical difficulties including system ill-conditioning and sensitivity to the choice of the master/slave pairing \cite{ainsworthCMAME01}. The master/slave bias can be eliminated via a two-pass procedure where each of the two surfaces serves as a master to the nodes of the other. The two-pass approach, however, has been shown to be over constrained \cite{papadopoulosCMAME92}. 

The Finite Element Tearing and Interconnecting (FETI) method is a dual domain decomposition approach first introduced by Farhat and Roux \cite{farhatIJNME91}. As the name suggests, the method is based on decomposing the given domain into a set of smaller subdomains that can be analyzed separately. The subdomains are connected using a set of discrete Lagrange multipliers applied to the interface nodes. In a conforming mesh, the Lagrange multipliers represent the forces required to achieve displacement continuity at the nodes, which translates into pointwise continuity along the interface. For non-conforming meshes, a field of Lagrange multipliers is introduced to enforce weak continuity along the non-conforming interface, yielding a set of discrete \textit{resultant} Lagrange multipliers at the nodes. This idea originated in the mortar method \cite{bernardiNPDEA92}. The interpolation of the Lagrange multipliers is based on the discretization of the master surface, thus the weak continuity condition is actually a projection of the slave kinematic variables onto the master surface. The FETI method has been widely used as a parallel iterative algorithm and a lot of research has been devoted to the development of pre-conditioners and iterative solvers for parallel efficiency. More details about the FETI method and its implementation can be found in references \cite{farhatCMAME94,parkIJNME97,farhatCMAME00,farhatIJNME01,farhatIJNME05,farhatIJNME07,bavestrelloCMAME07,bernardiCMAME08,bernardiCMAME09}.

Another popular member in the family of dual coupling methods is the mortar method. Originally proposed by Bernardi et al. \cite{bernardiNPDEA92} for two-dimensional problems, the mortar method was extended to three dimensions by Belgacem and Maday \cite{belgacemMMAN97}. The mortar concept is based on enforcing weak compatibility of the non-conforming interfaces using a Lagrange multiplier field interpolated on the master (or non-mortar) surface. The difference between the mortar and the non-conforming FETI approaches lies in the choice of the Lagrange multiplier field. A direct comparison gives the edge to the mortar method, based on the fact that it produces better-conditioned system matrices and an \textit{inf-sup} condition that is independent of mesh size \cite{lacourBIT97}. A later theoretical development of the mortar formulation can be traced to the work of Wohlmuth \cite{wohlmuthSIAMJNA01} who proposed a dual Lagrange multiplier space for triangular discretizations. Recently, Felmisch et al. applied this method to curved interfaces both in 2D \cite{flemischIJNME05} and 3D \cite{flemischCMAME07}. Bechet et al. \cite{bechetIJNME08} proposed a mortar formulation for the treatment of stiff interface conditions using the extended finite element method.

Jiao and Heath \cite{jiaoIJNME04} introduced the common refinement approach where the dual field is discretized as a linear interpolation on a \textit{commonly refined} mesh to which all interface nodes are mapped. Another dual coupling method of note can be found in the work of Park et al. \cite{parkIJNME00} and Felippa et al. \cite {felippaCMAME01} who proposed a multi-field variational framework that incorporates local and global Lagrange multiplier fields, in addition to the interface displacements. 

One major advantage of dual methods is that continuity of tractions is satisfied pointwise by the dual field and this family of methods passes the patch test by design. However, the stability of the coupling formulation hinges on the choice of the Lagrange multiplier discretization. Not all convenient primal/dual variable interpolations satisfy the LBB condition. Furthermore, the Lagrange multiplier field is based on the master or non-mortar side of the interface and the values at the slave nodes are obtained as a direct interpolation of the nodal master values. Therefore, this family of method is inherently biased by the choice of the master surface. Better results have been observed when the coarser of the two meshes was designated as the master surface, but the choice is not always trivial.

Recent trends in dual coupling methods have focused on developing stabilized methods to relax the LBB restrictions. This idea is well researched in the field of mixed methods for fluid and solid mechanics. A lot of the work on dual stabilized interface formulations has been done on embedded interfaces within the context of the eXtended Finite Element Method. Dolbow and coworkers \cite{jiIJNME04,mouradIJNME07,sandersIJNME08} have studied the stability of Lagrange multiplier formulation and proposed a bubble-enriched formulation that enables the enforcement of Dirichlet boundary and interface continuity conditions. Hansbo et al. \cite{hansboNM05} proposed a stabilized Lagrange multiplier method where the stabilization terms were calculated using Nitsche's method.

In summary, dual coupling methods introduce a field of Lagrange multipliers to enforce weak compatibility across non-conforming interfaces. This approach changes the pure displacement finite element formulation to a hybrid one that is governed by the LBB restriction. Dual interface formulations differ by their discretization of the Lagrange multiplier field. In all cases, however, this discretization is based in the master side of the interface. This family of methods, therefore, is naturally biased by the choice of the master surface. 

\subsection{Primal methods} \label{sec_DG:intro_primal}

In primal coupling methods, the interface is represented by its displacement fields and no dual variables are introduced. Primal methods are therefore not subject to the LBB restrictions, and are better suited to fit in a pure displacement framework. These methods, however, are challenged by the task of enforcing both geometric compatibility and continuity of tractions using a primal variable field. The fact that the discretization of the primary field is predetermined by the nonconforming meshes adds to the complexity of this challenge. These issues have adversely affected the popularity of primal coupling approaches.

The Discontinuous Galerkin framework is perhaps the most widely used primal coupling approach. This approach is natural to the coupling problem since it readily assumes discontinuous discretizations on all inter-element interfaces. DG methods have proved extremely efficient in modeling problems where the discontinuity of the primal fields mirrors physical events such as shocks. Outside these scenarios, however, DG methods can be very computationally extensive.
The DG formulation is based on identifying a set of target continuous fields for the displacement and traction fields on each interface, and mapping the discretized displacement and traction fields on each surface to these target fields, or \textit{numerical fluxes}, in a weak weighted residual form. The definition of these target fields or numerical fluxes is key to the stability of DG formulation. A number of DG approaches have been proposed in the literature. These include the methods of Brezzi et al. \cite{brezzietal}, Bassi et al. \cite{bassietal}, Local Discontinuous Galerkin \cite{cockburnSIAMJNA99,cockburnSIAMJNA01}, NIPG \cite{NIPG}, Babu\u{s}ka-Zl\'amal \cite{babuskaSIAMJNA73}, Brezzi et al. \cite{brezzietal2}, IP, Bassi-Rebay \cite{bassirebay} and Baumann-Oden \cite{baumannoden}. A good survey of these methods can be found in \cite{arnoldSIAMJNA02}. In this survey, the authors propose a unified framework for DG where all afore-mentioned methods fit with different choices for the numerical fluxes. A stability analysis reveals the methods of Baumann-Oden and Bassi-Rebay to be unstable. The methods of Babu\u{s}ka-Zl\'amal and Brezzi et al., although stable, are inconsistent. The remaining methods are consistent and stable due to a mesh-dependent penalty parameter.

The Nitsche method \cite{nitsche71} is a consistent primal formulation that employs a penalty approach to enforcing kinematic conditions. Originally introduced for the treatment of rough Dirichlet boundaries, this method is experiencing a resurgence in recent research. Hansbo and Hansbo \cite{hansboCMAME02} proposed a Nitsche-based formulation for unfitted meshes. Becker et al. \cite{beckerMMNA03} applied the method to the coupling of non-conforming meshes. A direct equivalence between this method and bubble-stabilized dual Lagrange multiplier formulations was established \cite{jiIJNME04,mouradIJNME07,sandersIJNME08}. The Nitsche method has been used as a basis for developing primal stabilized interface formulations for embedded interfaces \cite{dolbowIJNME08},\cite{NIPG}.

The Interface Element Method is a recent primal coupling method introduced by H.G. Kim \cite{kimCMAME02,kimCMAME03}. In this approach, all elements on either side of the non-conforming interface, designated as interface elements are modified to include \textit{pseudo nodes}, such that their surfaces match those of other elements across the interface in a Moving Least Squares sense. This method is very cumbersome and the resulting shape functions are not necessarily local to the interface element. Further developments of this method \cite{choCMAME05,choIJNME06I,choIJNME06II,limIJNME07} focused on the latter issue and lead to a modified element formulation where the weight functions and integration points in the element are designed to ensure locality. This approach, known as the variable-node element method, is very computationally extensive, especially in 3D \cite{kimCMAME08}, and requires a very involved formulation and a special integration rule for the enriched element.

The advantages of primal methods over dual ones are the unbiased treatment of the interface and the absence of the LBB restrictions. One general feature of all stable primal and dual methods is the weak enforcement of the continuity across the interface. This condition is imposed by the Lagrange multiplier field in dual methods, by a penalty-type formulation in primal methods, and by the MLS enrichment in the variable-node element method. The only method that enforces exact compatibility at a set of discrete locations is the MPC method. This method, however, is unstable as discussed in Section \ref{sec_DG:intro_dual}.

Exact geometric compatibility is an issue that is more relevant to contact problems than it is for coupling nonconforming meshes. In the latter, gaps or overlaps can only happen as a result of the interface compatibility error. In contact problems, however, gaps are possible and physically admissible, while overlaps are not. This issue is of a greater importance in non-smooth contact where capturing the exact geometry of the interface is essential to the accuracy of the solution. Therefore, since our goal is first and foremost contact problems, we consider an interface formulation that can enforce geometric compatibility at a set of discrete locations to be of great value. Such formulation would be directly compatible with the non-smooth contact constraint presented in Chapter 2. 

The aim of this chapter is to present a new primal approach for the coupling of non-conforming meshes that would combine the simplicity of the MPC approach with the consistency of other more involved primal methods such as DG. The proposed approach is based on a local enrichment of the non-conforming interface that would enable an unbiased enforcement of the displacement continuity constraints at all interface nodes without introducing additional variables or increasing the size of the problem. We show that a local DG-based interface stabilization procedure guarantees a consistent transfer of the traction field across the non-conforming interface. The result is a robust two-pass formulation that passes the patch test and displays excellent convergence rates with mesh refinement.

The outline of the chapter is as follows. In Section \ref{sec_DG:coupling_form}, we present the mathematical formulation of the equations of motion. Section \ref{sec_DG:patch_test} discusses the necessary requirements for passing the patch test, and Section \ref{sec_DG:mpc} analyzes the performance of the MPC method in light of these requirements. We then describe the suggested approach for the formulation of the contact constraints in Section \ref{sec_DG:prop_form} and present some numerical examples in Section \ref{sec_DG:results}. Section \ref{sec_DG:conclusion} discusses conclusions and future work.


%
%
%
%
\section{Formulation of the coupling problem} \label{sec_DG:coupling_form}
\begin{figure}
	\centering
	\includegraphics[clip]{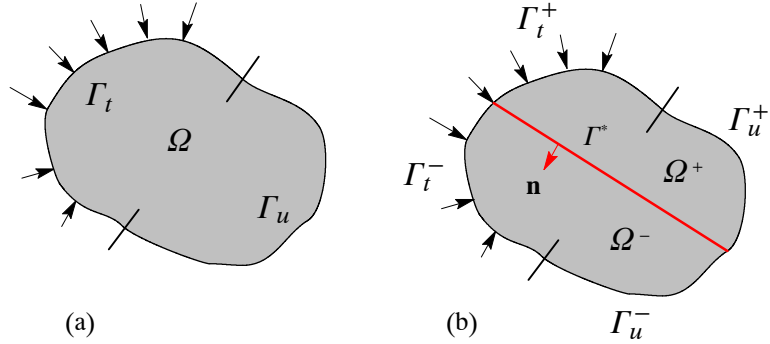}
	\caption{The coupling problem: (a) original domain definition $\Omega$ (b) partitioning into two subdomains $\Omega=\Omega^+\cup\Omega^-$  \label{fig:DG_coupling}}
	\end{figure}

First, let us revisit the formulation of the coupling problem for the case of static finite elasticity (the formulation of the dynamic case can be obtained in a similar manner, but we chose to omit the inertia terms for the sake of simplicity). Consider the open domain $\Omega$ shown in Figure \ref{fig:DG_coupling} (a). Let $\Gamma=\Gamma_t\cup\Gamma_u$ be the boundary of the domain $\Omega$, where $\Gamma_t$ and $\Gamma_u$ denote the Neumann and Dirichlet parts of that boundary, respectively, and $\Gamma=\Gamma_t\cap\Gamma_u=\Phi$. The strong form of the governing equations for the continuum problem can be stated as follows: 
(S) Given the body force vector $\mathbf{b}$, the prescribed tractions $\mathbf{h}$ on $\Gamma_t$ and the prescribed displacement
field $\mathbf{g}$ on $\Gamma_u$ find $\mathbf{u}$ such that
\begin{equation} \label{eq_DG:Sequil}
\mathrm{div}\mathbf{S}+\mathbf{b} = \mathbf{0} \ \textrm{in} \ \Omega; \ \ \ \mathbf{Sn} = \mathbf{h} \ \textrm{on} \ \Gamma_t; \ \ \ \mathbf{u} = \mathbf{g} \ \textrm{on} \ \Gamma_u,
\end{equation}
where '$\mathrm{div}$' is the divergence operator and $\mathbf{S}$ is the Cauchy stress tensor. Defining the spaces $S = \left\{\mathbf{u}\in H^1(\Omega):\mathbf{u}=\mathbf{g} \ \textrm{on} \ \Gamma_u\right\}$ and $V = \left\{\mathbf{w} \in H^1(\Omega):\mathbf{w}=\mathbf{0} \ \textrm{on} \ \Gamma_u\right\}$, the weighted residual form of \eqref{eq_DG:Sequil} can be written as:
(W) Given the body force vector $\mathbf{b}$, the prescribed traction field $\mathbf{h}$ on $\Gamma_t$ and the prescribed displacement
field $\mathbf{g}$ on $\Gamma_u$ find $\mathbf{u}\in S$ such that
\begin{equation} \label{eq_DG:Wequil}
\int_{\Omega}\left[\mathrm{div}\mathbf{S}+\mathbf{b}\right]\cdot\mathbf{w}\,d \Omega + \int_{\Gamma_t}\left[\mathbf{h}-\mathbf{Sn}\right]\cdot\mathbf{w}\,d \Gamma = 0,
\end{equation}
for all $\mathbf{w}\in V$. The \textit{symmetric} form of equation \eqref{eq_DG:Wequil} can be found by applying the divergence theorem to the first term in this equation as follows.
\begin{equation} \label{eq_DG:Wequils}
\int_{\Omega}\left[-\mathbf{S}\cdot\nabla_{\mathbf{x}}\mathbf{w}+\mathbf{b}\cdot\mathbf{w}\right],d \Omega + \int_{\Gamma_t}\mathbf{h} \cdot\mathbf{w}\,d \Gamma = 0,
\end{equation}
where $\mathbf{x}$ is the spatial field variable. To simplify the subsequent derivations, we will drop the subscript in the gradient notation, and assume the gradient to be with respect to the spatial coordinates coordinates $\mathbf{x}$, unless otherwise noted.

Now consider a partition of the domain $\Omega$ into two non-overlapping open subdomains $\Omega^+$ and $\Omega^-$ such that $\Omega =\bar{\Omega}^+ \cup \bar{\Omega}^-$ as shown in Figure \ref{fig:DG_coupling} (b). Let $\Gamma^+$ and $\Gamma^-$ be the (outer) subdomain boundaries and $\Gamma^*$ the smooth interface with normal $\mathbf{n}$ between the two subdomains $\Gamma^* = \bar{\Omega}^+\cap\bar{\Omega}^-$. The choice of the superscript in naming the subdomains is based on the definition of $\mathbf{n}$ as the outward pointing normal to $\Gamma^*$ in $\Omega^+$. Consequently, the outward pointing normal to $\Gamma^*$ in $\Omega^-$ is $-\mathbf{n}$. The Neumann and Dirichlet parts of the subdomain boundary are
\begin{align} \notag
\Gamma_t^{+}=\Gamma_t\cap\Gamma^{+}, \ \ \ \ \Gamma_u^{+}&=\Gamma_u\cap\Gamma^{+}, \ \ \ \ \Gamma_t^{+}\cap\Gamma_u^{+} = \phi  \\
\Gamma_t^{-}=\Gamma_t\cap\Gamma^{-}, \ \ \ \ \Gamma_u^{-}&=\Gamma_u\cap\Gamma^{-}, \ \ \ \ \Gamma_t^{-}\cap\Gamma_u^{-} = \phi.
\end{align}
We denote the displacement fields in $\Omega^+$ and $\Omega^-$ as $\mathbf{u}^+ = \mathbf{u}|_{\Omega^+}$ and $\mathbf{u}^-= \mathbf{u}|_{\Omega^-}$, respectively. Similarly, $\mathbf{S}^+ = \mathbf{S}|_{\Omega^+}$ and $\mathbf{S}^-= \mathbf{S}|_{\Omega^-}$. The subdivision is completely artificial, and is not based on any physical discontinuity in the original problem. Thus, to preserve the consistency of the formulation, continuity of the displacement and traction fields has to be maintained along the interface. The governing equations for the coupled problem are
\begin{align} \label{eq_DG:Sequil_+}
\mathrm{div}\mathbf{S}^{+}+\mathbf{b} = \mathbf{0} \ \textrm{in} \ \Omega^+; \ \ \ \mathbf{S}^{+}\mathbf{n} = \mathbf{h} \ \textrm{on} \ \Gamma_t^{+}; \ \ \ \mathbf{u}^+ = \mathbf{g} \ \textrm{on} \ \Gamma_u^{+} \\
\label{eq_DG:Sequil_-}
\mathrm{div}\mathbf{S}^{-}+\mathbf{b} = \mathbf{0} \ \textrm{in} \ \Omega^-; \ \ \ \mathbf{S}^{-}\mathbf{n} = \mathbf{h} \ \textrm{on} \ \Gamma_t^{-}; \ \ \ \mathbf{u}^- = \mathbf{g} \ \textrm{on} \ \Gamma_u^{-} \\
\label{eq_DG:Sequil_S*}
\mathbf{S}^+\mathbf{n} = \mathbf{S}^-\mathbf{n} \ \textrm{on} \ \Gamma^{*}; \\ \label{eq_DG:Sequil_u*}
\mathbf{u}^{+} = \mathbf{u}^{-} \ \textrm{on} \ \Gamma^{*}.
\end{align}
We define
\begin{align}
\mathbf{w}^+ \in V^+ = \left\{\mathbf{w}^+ \in H^1(\Omega^+):\mathbf{w}^+=\mathbf{0} \ \textrm{on} \ \Gamma_u^+\right\} \\
\mathbf{w}^- \in V^- = \left\{\mathbf{w}^- \in H^1(\Omega^-):\mathbf{w^-}=\mathbf{0} \ \textrm{on} \ \Gamma_u^-\right\} \\
\mathbf{w}^*\in V^* = \left\{\mathbf{w}^* \in H^{1/2}(\Gamma^*)\right\}
\end{align}
to be the variational displacement fields in $\Omega^+$ , $\Omega^-$ and along $\Gamma^*$, respectively. The weighted residual form of equations \eqref{eq_DG:Sequil_-}-\eqref{eq_DG:Sequil_S*} can be written as follows:
\begin{align} \notag
\int_{\Omega^+}\left[\mathrm{div}\mathbf{S}^++\mathbf{b}\right]\cdot\mathbf{w}^+\,d \Omega &+ \int_{\Omega^-}\left[\mathrm{div}\mathbf{S}^-+\mathbf{b}\right]\cdot\mathbf{w}^-\,d \Omega 
-\int_{\Gamma^*} \left[\mathbf{S}^+-\mathbf{S}^-\right]\mathbf{n}\cdot \mathbf{w}^*\,d\Gamma  \\  \label{eq_DG:WCequil} +
\int_{\Gamma_t^+}\left[\mathbf{h}-\mathbf{S}^+\mathbf{n}\right]\cdot\mathbf{w}^+\,d \Gamma &+ \int_{\Gamma_t^-}\left[\mathbf{h}+\mathbf{S}^-\mathbf{n}\right]\cdot\mathbf{w}^-\,d \Gamma = 0,
\end{align}
for all $\mathbf{w}^+ \in V^+,\mathbf{w}^{-}\in V^-$, and $\mathbf{w}^{*}\in V^*$. In this primal formulation of the coupled problem, the displacement continuity condition \eqref{eq_DG:Sequil_u*} does not appear in the weak form and is treated as a \textit{Dirichlet-type} condition for the continuous problem. The variational displacement fields can be arbitrary and need not obey the continuity constraint that applies to the real field. 

Applying the divergence theorem to the terms $\int_{\Omega^+}\mathrm{div}\mathbf{S}^+ \cdot\mathbf{w}^+\,d \Omega$ and $\int_{\Omega^-}\mathrm{div}\mathbf{S}^- \cdot\mathbf{w}^-\,d \Omega$ in equation \eqref{eq_DG:WCequil}, and imposing homogeneous boundary conditions on the Dirichlet boundary yields
\begin{align} \notag
\int_{\Omega^+}\left[-\mathbf{S}^+\cdot\nabla\mathbf{w}^{+}+\mathbf{b}\cdot\mathbf{w}^{+}\right]\,d \Omega 
+\int_{\Omega^-}\left[-\mathbf{S}^-\cdot\nabla\mathbf{w}^{-}+\mathbf{b}\cdot\mathbf{w}^{-}\right]\,d \Omega  \\ \notag
+\int_{\Gamma^*}\mathbf{S}^+\mathbf{n}\cdot\mathbf{w}^+\,d\Gamma - \int_{\Gamma^*}\mathbf{S}^-\mathbf{n}\cdot\mathbf{w}^-\,d\Gamma 
-\int_{\Gamma^*} \left[\mathbf{S}^+-\mathbf{S}^-\right]\mathbf{n}\cdot \mathbf{w}^* \,d\Gamma \\ \label{eq_DG:WCequil_div}
+\int_{\Gamma_t^+}\left[\mathbf{h}-\mathbf{S}^+\mathbf{n}\right]\cdot\mathbf{w}^{+}+\int_{\Gamma_t^-}\left[\mathbf{h}+\mathbf{S}^-\mathbf{n}\right]\cdot\mathbf{w}^{-}\,d\Gamma=0.
\end{align}
Simplifying and rearranging the terms in equation \eqref{eq_DG:WCequil_div} we get
\begin{align} \notag
\int_{\Omega^+}\left[-\mathbf{S}^+\cdot\nabla\mathbf{w}^{+}+\mathbf{b}\cdot\mathbf{w}^{+}\right]\,d \Omega 
+\int_{\Omega^-}\left[-\mathbf{S}^-\cdot\nabla\mathbf{w}^{-}+\mathbf{b}\cdot\mathbf{w}^{-}\right]\,d \Omega  \\ \label{eq_DG:WCequil_div2}
+\int_{\Gamma_t^+}\mathbf{h}\cdot\mathbf{w}^{+}\,d\Gamma+ \int_{\Gamma_t^-}\mathbf{h}\cdot\mathbf{w}^{-}\,d\Gamma + I_A -I_B= 0. 
\end{align}
where
\begin{align} \notag
I_A &= \int_{\Gamma^*}\mathbf{S}^+\mathbf{n}\cdot\mathbf{w}^+\,d\Gamma - \int_{\Gamma^*}\mathbf{S}^-\mathbf{n}\cdot\mathbf{w}^-\,d\Gamma  \\ \label{eq_DG:IAIB}
I_B &= \int_{\Gamma^*} \left[\mathbf{S}^+-\mathbf{S}^-\right]\mathbf{n}\cdot \mathbf{w}^*\,d\Gamma
\end{align}
The term $I_A$ is often referred to as the interface \textit{jump}, while $I_B$ is the weak statement of traction equilibrium along the interface. 
From equations \eqref{eq_DG:WCequil_div2} and \eqref{eq_DG:IAIB}, we can observe that for the choice $\mathbf{w}^+ = \mathbf{w}^-=\mathbf{w}^*$, the interface term $I_A - I_B$ vanishes and 
equation \eqref{eq_DG:WCequil_div2} simplifies as
\begin{align} \notag
\int_{\Omega^+}\left[-\mathbf{S}^+\cdot\nabla\mathbf{w}^+ +\mathbf{b}\cdot\mathbf{w}^{+}\right]\,d \Omega 
+\int_{\Omega^-}\left[-\mathbf{S}^-\cdot\nabla\mathbf{w}^- +\mathbf{b}\cdot\mathbf{w}^{-}\right]\,d \Omega  \\ \label{eq_DG:Gequil_div2}
+\int_{\Gamma_t^+}\mathbf{h}\cdot\mathbf{w}^{+}\,d\Gamma+ \int_{\Gamma_t^-}\mathbf{h}\cdot\mathbf{w}^{-}\,d\Gamma = 0. 
\end{align}
Therefore the solution to the coupled problem can be obtained by solving the two sub-problems in $\Omega^+$, $\Omega^-$, without the need for explicitly enforcing the continuity of tractions. The continuity of the displacement field, however, still needs to be addressed at the interface. This result can be extrapolated to the case of multiple domains or finite elements, and is crucial to understanding the performance of finite element methods. In a conforming Bubnov-Galerkin finite element discretization, the real and variational displacement fields are point-wise continuous by design, which simultaneously ensures geometric compatibility and a complete transfer of forces across the interface. This result can be easily obtained for the case of a simple patch test, as shown in the following section. 

\begin{remark}  \label{rem_DG:conf}
Equation \eqref{eq_DG:Gequil_div2} can be interpreted to be the equivalent of the symmetric weighted residual equilibrium form of equation \eqref{eq_DG:Wequils} for the coupled problem. The equivalence, however, holds only for a conforming choice of the variational field $\mathbf{w}^+ = \mathbf{w}^-=\mathbf{w}^*$ for which the interface jump term vanishes exactly. 
\end{remark}

\begin{remark}  \label{rem_DG:nonconf}
For a nonconforming variational field, the weighted residual of equation \eqref{eq_DG:Gequil_div2} does not take into account the jump resulting from the discontinuity of tractions across the interface, and does not lead to a solution where the equilibrium of tractions on the interface is satisfied. To amend these issues, the terms $I_A$ and $I_B$ have to be included in the variational formulation.
\end{remark}

%
%
%
%
%
%
%
\section{The patch test} \label{sec_DG:patch_test}
	\begin{figure}
	\centering
	\includegraphics[clip]{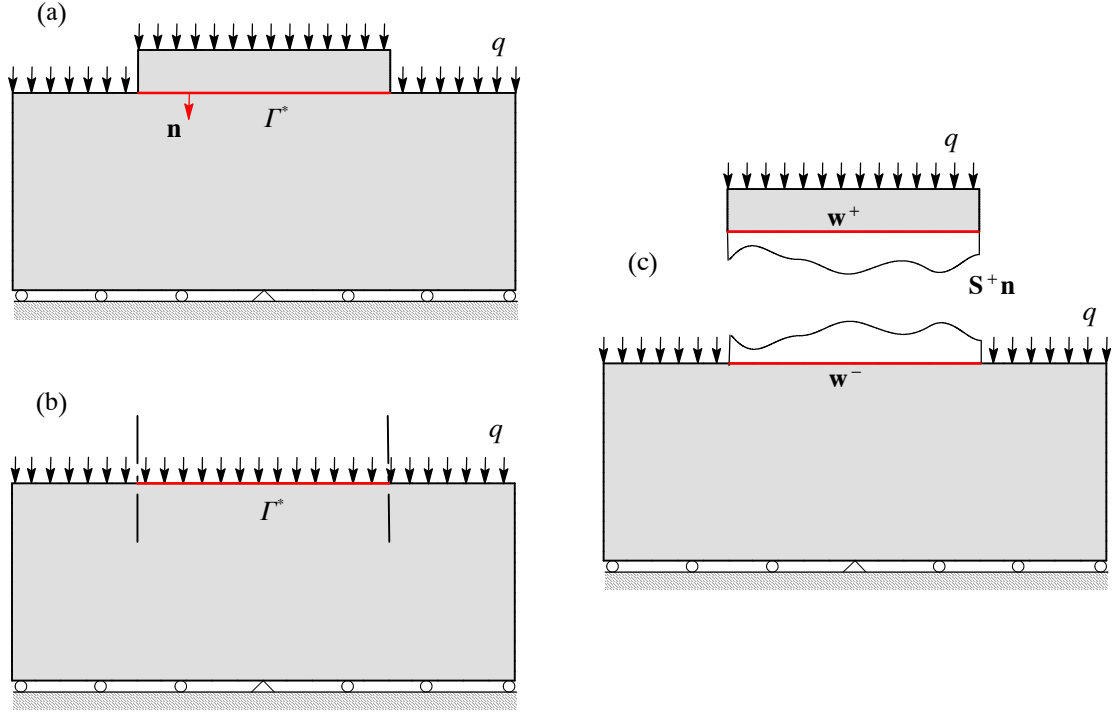}
	\caption{The patch test for (a) the punch-foundation system and (b) the foundation only, and (c) free-body diagram of the punch-foundation system  \label{fig:DG_patch}}
	\end{figure}
The patch test shown in Figure \ref{fig:DG_patch}(a) consists of a punch connected with a rectangular foundation through an interface $\Gamma^*$. The punch and foundation are made of a linear elastic material with properties $E$ and $\nu$. A constant traction vector $\mathbf{q}=-q\,\mathbf{e}_2$ is applied to the free surface of the structure. Let $\mathbf{u}^+$, $\mathbf{w}^+$ be the displacement and variational displacement fields in the punch, and $\mathbf{u}^-$, $\mathbf{w}^-$ their counterparts in the foundation. The exact solution for this problem dictates that the traction field on the interface $\Gamma^*$ between the punch and the foundation be equal to the applied traction. Therefore, the displacement and stress fields in the foundation obtained using this configuration should be equal to those obtained when the load is applied directly to the surface of the foundation without the presence of the punch, as shown in Figure \ref{fig:DG_patch}(b).

Consider the free body diagram shown in Figure \ref{fig:DG_patch}(c). The state of constant stress in the punch implies that the traction field on its lower surface is equal to the applied load on its upper surface $\mathbf{S}^+\mathbf{n} = \mathbf{q}$. Therefore, the resultant traction field applied by the punch to the foundation through the interface is
\begin{equation}
\int_{\Gamma^*}\mathbf{S}^+\mathbf{n}\cdot\mathbf{w}^+\,d\Gamma = \int_{\Gamma^*}\mathbf{q}\cdot\mathbf{w}^+\,d\Gamma
\end{equation}

For the case shown in Figure \ref{fig:DG_patch}(b), the total equivalent traction applied to the top surface of the foundation is $\int_{\Gamma^{*}}\mathbf{q}\cdot\mathbf{w}^-\,d\Gamma$. Therefore, for the two loading scenarios to be equivalent, the following has to be satisfied
\begin{equation}
\int_{\Gamma^{*}}\mathbf{q}\cdot\mathbf{w}^+\,d\Gamma = \int_{\Gamma^{*}}\mathbf{q}\cdot\mathbf{w}^-\,d\Gamma,
\end{equation}
\begin{equation} \label{eq_DG:patch_cond}
\mathbf{q}\cdot\int_{\Gamma^{*}}\left[\mathbf{w}^+-\mathbf{w}^-\right]\,d\Gamma = 0.
\end{equation}
It should be clear from equation \eqref{eq_DG:patch_cond} that when the variational displacement is point-wise continuous $\mathbf{w}^+=\mathbf{w}^-$, the resultant traction vector applied by the punch to the foundation is exactly equal to the resultant of the external pressure. Therefore, the transfer of tractions across the interface is complete and the the patch test is passed. As discussed in the previous section, this condition is satisfied by design in a conforming mesh in the context of the Bubnov-Galerkin method. In a non-conforming mesh, the choice of the variational fields has to satisfy the weak continuity condition of \eqref{eq_DG:patch_cond}. This result explains the stability of the dual and primal methods discussed in Section \ref{sec_DG:intro} and is key to understanding the patch test performance of the MPC method, as we show in the following section. 
%
%
%

\section{The Multi-Point-Constraints method} \label{sec_DG:mpc}
	\begin{figure}
	\centering
	\includegraphics[clip]{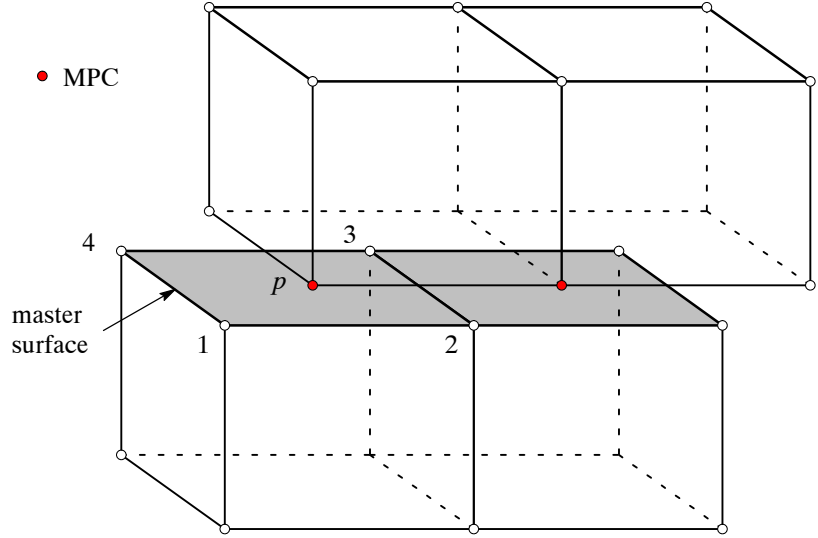}
	\caption{The Multi-Point Constraint method  \label{fig:DG_MPC}}
	\end{figure}
The Multi-Point-Constraints method seeks to enforce strong compatibility between a set of nodes on one side of the non-conforming interface and the opposing element surface, such that each node remains glued to the corresponding element surface during deformation. This condition precludes the possibility of the node either separating from or penetrating the element thereby  creating physically inadmissible gaps or overlaps. This method is most commonly used as a single-pass approach where one side of the interface is designated as the \textit{master} surface to which the \textit{slave} nodes of the opposing side are tied. The constraint is expressed as a linear combination involving the displacements of the slave node and those of the master element nodes, hence its designation as \textit{Multi-Point Constraint}. For the case shown in Figure \ref{fig:DG_MPC}, the constraint takes the shape
\begin{equation}
\mathbf{u}^p = \sum_{i=1}^4 N^i(p)\mathbf{u}^i.
\end{equation}
where $N^i(p), i=1,...,4$ are the shape functions associated with the surface nodes $1,2,3,4$, evaluated at $p$. Given its straightforward implementation, this method has been extensively used in sub-structuring and domain decomposition problems. As a result of applying this constraint, the degrees of freedom of the slave node are eliminated in favor of those of the master element. Therefore, the single-pass MPC method is inherently biased by the choice of the master/slave pairing.  Furthermore, as reported in \cite{papadopoulosCMAME92}, the single-pass MPC method fails the patch test. The master/slave bias can be eliminated by applying the method in a two-pass fashion where each side of the interface serves as master to the nodes of the other. When implemented as a two-pass approach, the MPC method passes the patch test in the case of bilinear elements, but fails it when higher-order elements are used. Furthermore, the two-pass MPC method is known to exhibit locking or artificial stiffening of the interface due to the excessive elimination of interface degrees of freedom. Through the following discussion, we show that these results are directly correlated with the interpolation of the variational displacement field along the interface. We consider the simple patch test presented in Section \ref{sec_DG:patch_test} and examine the displacement fields in the punch and foundation after enforcing the MPCs. We address the cases of planar bilinear and quadratic elements, and restrict our attention to the two-dimensional domain for simplicity.

\subsection{Bilinear elements}
	\begin{figure}
	\centering
	\includegraphics[clip]{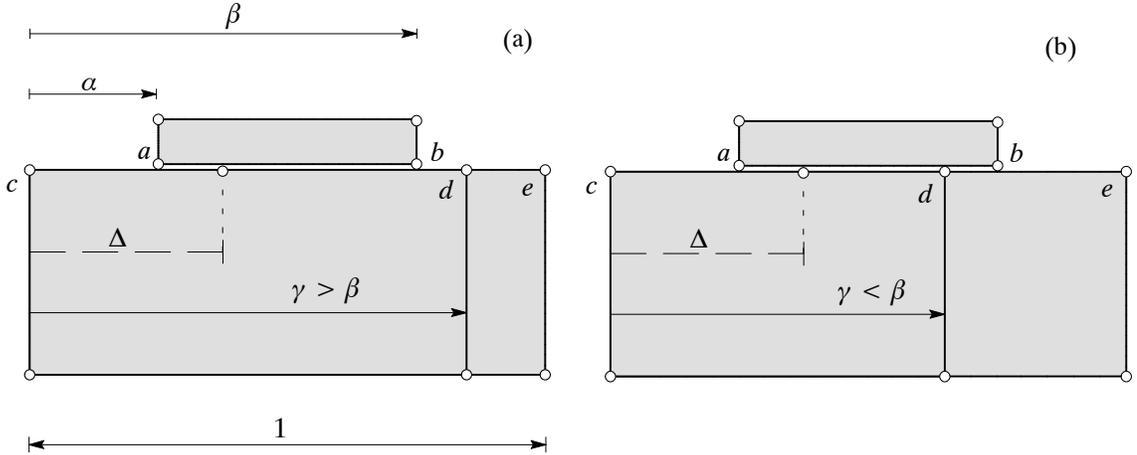}
	\caption{Non-conforming discretization using Q4 elements with (a) two non-matching nodes and (b) three non-matching nodes  \label{fig:DG_MPC_Q4}}
	\end{figure}
First we consider the mesh shown in Figure \ref{fig:DG_MPC_Q4}(a), where the foundation is discretized using two Q4 elements whereas the punch is represented by a single Q4 element that connects to the foundation through the interface $ab$. The foundation is of unit width and $\alpha$, $\beta$ and $\gamma$ denote the horizontal coordinates of nodes $a$, $b$, and $c$, respectively, measured with respect to the left-end of the structure. Let $\mathbf{u}^+$, $\mathbf{u}^-$ be the displacement fields in the punch and foundation, as defined previously. The displacement at a given location $\Delta >\alpha$ on the upper (punch) and lower (foundation) sides of the interface can be expressed as follows,
\begin{align} \label{eq_DG:UDelta+_Q4}
\mathbf{u}^{+\Delta} &= (1-\frac{\Delta - \alpha}{\beta-\alpha}) \,\mathbf{u}^a + \frac{\Delta - \alpha}{\beta-\alpha}\,\mathbf{u}^b \\ \notag
\mathbf{u}^{-\Delta} &= (1-H_{\Delta,\gamma})\left[(1-\frac{\Delta}{\gamma}) \,\mathbf{u}^c + \frac{\Delta}{\gamma}\,\mathbf{u}^d\right] \\ \label{eq_DG:UDelta-_Q4}
 &+ H_{\Delta,\gamma}\left[(1-\frac{\Delta - \gamma}{1-\gamma}) \,\mathbf{u}^d + \frac{\Delta - \gamma}{1-\gamma} \,\mathbf{u}^e\right],
\end{align}
where $H_{\Delta,\gamma}$ is the Heaviside step function $H_{\Delta,\gamma} = 1$ if $\Delta > \gamma$ and 0 otherwise.

We treat $\gamma$ as a variable that measures the relative positioning of the foundation mesh divider with respect to the punch, and we consider the following cases. 
\subsubsection{Case I: $\gamma > \beta$} \label{sec_DG:Q4caseI}
This case corresponds to the configuration shown in Figure \ref{fig:DG_MPC_Q4} (a), where both punch nodes $a$ and $b$ meet the left-side foundation element and no nodes from the foundation mesh are located on the interface. This leaves only the possibility of a single-pass approach to implement the MPCs that tie the punch nodes to the foundation surface as follows
\begin{align} \notag
\mathbf{u}^a = (1-\frac{\alpha}{\gamma}) \,\mathbf{u}^c + \frac{\alpha}{\gamma}\,\mathbf{u}^d \\ \label{eq_DG:SPass_uaub_Q4}
\mathbf{u}^b = (1-\frac{\beta}{\gamma}) \,\mathbf{u}^c + \frac{\beta}{\gamma} \,\mathbf{u}^d.
\end{align}
Given these constraints, we rewrite the displacement field in the punch $\mathbf{u}^{+\Delta}$ and compare it to the displacement $\mathbf{u}^{-\Delta}$ of the foundation. Substituting equation \eqref{eq_DG:SPass_uaub_Q4} into \eqref{eq_DG:UDelta-_Q4} we can write the displacement field in the punch as a function of the foundation nodes displacements as,
\begin{align}\notag
\mathbf{u}^{+\Delta} &= (1-\frac{\Delta - \alpha}{\beta-\alpha}) \,\mathbf{u}^a + \frac{\Delta - \alpha}{\beta-\alpha}\,\mathbf{u}^b   \\ \label{eq_DG:UDelta+}
 &= (1-\frac{\Delta - \alpha}{\beta-\alpha})\left[(1-\frac{\alpha}{\gamma}) \,\mathbf{u}^c + \frac{\alpha}{\gamma}\,\mathbf{u}^d \right] + \frac{\Delta - \alpha}{\beta-\alpha} \left[ (1-\frac{\beta}{\gamma}) \,\mathbf{u}^c + \frac{\beta}{\gamma} \,\mathbf{u}^d \right] 
\end{align}
The terms in \eqref{eq_DG:UDelta+} can be condensed by noting that
\begin{align} \label{eq_DG:Delta_simple}
\left(1-\frac{\Delta - \alpha}{\beta-\alpha}\right) \frac{\alpha}{\gamma} + \left(\frac{\Delta - \alpha}{\beta-\alpha}\right)\frac{\beta}{\gamma}  
= \frac{\alpha}{\gamma}+ \left( \frac{\Delta - \alpha}{\beta-\alpha}\right) \left(\frac{\beta -\alpha}{\gamma}\right) = \frac{\Delta}{\gamma}.
\end{align}
Therefore the displacement field in the punch can be simply expressed as a linear combination of $\mathbf{u}_c$ and $\mathbf{u}_d$ as follows
\begin{equation}
\mathbf{u}^{+\Delta} = (1-\frac{\Delta}{\gamma}) \,\mathbf{u}^c + \frac{\Delta}{\gamma} \,\mathbf{u}^d.
\end{equation}
This expression is identical to the displacement field in the foundation $\mathbf{u}^{-\Delta}$ for $\Delta < \gamma$ (equation \eqref{eq_DG:UDelta-_Q4}). Therefore, the displacement field, and consequently the variational field, is point-wise continuous along the interface and the formulation passes the patch test. The deformation of the interface is defined by the degrees of freedom of the foundation. 
\subsubsection{Case II: $\alpha < \gamma < \beta$}
The foundation mesh in this case is positioned such that its node $d$ falls on the interface, as shown in Figure \ref{fig:DG_MPC_Q4} (b). We compare the two approaches for enforcing the MPCs.
\\
\\
\textit{ a) Single-pass approach}
\\
The choice of the master surface in this case is not trivial. To facilitate the comparison of the results with the previous case we choose to repeat the designation of the foundation as master surface, and enforce the MP constraints,
\begin{align} \notag
\mathbf{u}^a &= (1-\frac{\alpha}{\gamma}) \,\mathbf{u}^c + \frac{\alpha}{\gamma}\,\mathbf{u}^d  \\
\mathbf{u}^b &= (1-\frac{\beta - \gamma}{1-\gamma}) \,\mathbf{u}^d + \frac{\beta - \gamma}{1-\gamma}\,\mathbf{u}^e.
\end{align}
Following a procedure similar to Case I, we incorporate the above-constraints into the displacement field expression in the punch $\mathbf{u}^{+\Delta}$ (equation \ref{eq_DG:UDelta+_Q4}) leading to,
\begin{align}\label{eq_DG:UDelta+_1P}
\mathbf{u}^{+\Delta} &= (1-\frac{\Delta - \alpha}{\beta-\alpha}) \left[ (1-\frac{\alpha}{\gamma}) \,\mathbf{u}^c + \frac{\alpha}{\gamma}\,\mathbf{u}^d \right ] + \frac{\Delta - \alpha}{\beta-\alpha}\left[ (1-\frac{\beta - \gamma}{1-\gamma}) \,\mathbf{u}^d + \frac{\beta - \gamma}{1-\gamma}\,\mathbf{u}^e \right ]  
\end{align}
To simplify equation \eqref{eq_DG:UDelta+_1P} we introduce the following variables,
\begin {equation} \label{eq_DG:lambda}
\lambda^a=\frac{\gamma-\alpha}{\gamma},\ \ \ \lambda^b=\frac{\beta-\gamma}{1-\gamma},\ \ \ \lambda^\Delta=\frac{\Delta-\alpha}{\beta-\alpha}
\end{equation}
The variables $\lambda^a$ and $\lambda^b$ measure the relative positioning of the interface node $d$ with respect to the punch nodes $a$ and $b$. Note that $\lambda^a = 0$ when $\gamma = \alpha$ and $\lambda^b = 0$ when $\gamma = \beta$, and that these two variables cannot be zero simultaneously. $\lambda^\Delta$ measures the relative position of $\Delta$ with respect to the punch nodes $a$ and $b$. Substituting back in equation \eqref{eq_DG:UDelta+_1P} we find
\begin{align} \notag
\mathbf{u}^{+\Delta} 
 &= (1-\lambda^\Delta) (\lambda^a) \,\mathbf{u}^c +  \lambda^\Delta\lambda^b \,\mathbf{u}^e \\ \label{eq_DG:UDelta+_1P_2}
 &+ \left[ (1-\lambda^\Delta) (1-\lambda^a) + \lambda^\Delta (1-\lambda^b)\right ] \,\mathbf{u}^d. 
\end{align}
Equation \eqref{eq_DG:UDelta+_1P_2} describes the displacement field in the punch in terms of the interface degrees of freedom $\mathbf{u}^c$, $\mathbf{u}^d$ and $\mathbf{u}^e$. This result obviously does not match with the foundation interface displacement \eqref{eq_DG:UDelta-_Q4} since the latter would include either $\mathbf{u}^c$ and $\mathbf{u}^d$, or $\mathbf{u}^d$ and $\mathbf{u}^e$, depending on the location $\Delta$. Therefore, the displacement field is discontinuous across the interface, which translates into an incomplete transfer of the traction field. The only possibility for pointwise continuity is when either $\lambda^a = 0$ or $\lambda^b = 0$, in which case the current configuration reverts to the case described in Case I.
\\
\\
\textit{ b) Two-pass approach}
\\
In this case the MPCs are imposed at nodes $a$, $b$ and $d$ simultaneously as follows
\begin{align} \notag
\mathbf{u}^a &= \lambda^a \,\mathbf{u}^c + (1-\lambda^a)\,\mathbf{u}^d  \\ \label{eq_DG:DPass_uaub_Q4}
\mathbf{u}^b &= (1-\lambda^b) \,\mathbf{u}^d + \lambda^b\,\mathbf{u}^e  \\ \notag 
\mathbf{u}^d &= (1-\frac{\gamma-\alpha}{\beta-\alpha}) \,\mathbf{u}^a 
+ \frac{\gamma-\alpha}{\beta-\alpha} \,\mathbf{u}^b \\ \label{eq_DG:DPass_ud_Q4}
&= (1-\lambda^d) \,\mathbf{u}^a + \lambda^d\,\mathbf{u}^b,
\end{align}
where $\lambda^d=\displaystyle\frac{\gamma-\alpha}{\beta-\alpha}$ and $\lambda^a$, $\lambda^b$ are defined in \eqref{eq_DG:lambda}. Note that the expressions for $\mathbf{u}^a$ and $\mathbf{u}^b$ include $\mathbf{u}^d$, which, in turn, is a function of $\mathbf{u}^a$ and $\mathbf{u}^b$. Substituting equation \eqref{eq_DG:DPass_ud_Q4} into equation \eqref{eq_DG:DPass_uaub_Q4}, we find
\begin{align} \notag
\mathbf{u}^d &= (1-\lambda^d) \left[ \lambda^a \,\mathbf{u}^c + (1-\lambda^a)\,\mathbf{u}^d \right] + \lambda^d\ \left[ (1-\lambda^b) \,\mathbf{u}^d + \lambda^b\,\mathbf{u}^e \right] \\ \label{eq_DG:DPass_ud1_Q4}
\mathbf{u}^d \,&\left[1-(1-\lambda^d)\,(1-\lambda^a)-\lambda^d\,(1-\lambda^b) \right]
=\lambda^a\,(1-\lambda^d) \,\mathbf{u}^c + \lambda^b\,\lambda^d \,\mathbf{u}^e.
\end{align}
After some manipulation of the left-hand side, equation \eqref{eq_DG:DPass_ud1_Q4} can be written as,
\begin{align} \label{eq_DG:DPass_ud2_Q4}
\left[ (1-\lambda^d)\lambda^a+\lambda^d\,\lambda^b\right]\mathbf{u}^d = (1-\lambda^d)\lambda^a \,\mathbf{u}^c + \lambda^d\,\lambda^b \,\mathbf{u}^e.
\end{align}
Equation \eqref{eq_DG:DPass_ud2_Q4} implies that, unless the term $\left[(1-\lambda^d)\lambda^a+\lambda^d\,\lambda^b\right]$ is equal to zero, the displacement vector $\mathbf{u}^d$ is a linear interpolation of $\mathbf{u}^c$ and $\mathbf{u}^e$. Therefore, in addition to eliminating the degrees of freedom of the punch for those of the foundation, some of the foundation degrees of freedom are lost as well as a result of applying the MPCs in a two pass approach. Therefore, the flexibility of the foundation at the interface is greatly reduced, a phenomenon called \textit{surface locking}.
The only possibility for $\mathbf{u}^d$ to be independent from $\mathbf{u}^c$ and $\mathbf{u}^e$ occurs when
\begin{equation} \label{eq_DG:nolocking} \left[(1-\lambda^d)\lambda^a+\lambda^d\,\lambda^b\right] =0.
\end{equation}
Since all terms on the left-hand side of this equality are positive, equation \eqref{eq_DG:nolocking} requires the following to be true
\begin{align} \label{eq_DG:nolocking2}
(1-\lambda^d)\lambda^a = 0, \ \ \ \lambda^d\,\lambda^b =0.
\end{align}
As pointed out from equation \eqref{eq_DG:lambda}, $\lambda^a$ and $\lambda^b$ cannot be zero simultaneously. Furthermore, when $\lambda^a = 0$, $\lambda^d = 0$. Also, when $\lambda^b = 0$, $\lambda^d = 1$. Therefore, the conditions \eqref{eq_DG:nolocking2} can only be met when 
\begin{equation} \label{eq_DG:nolocking3}
\lambda^a = 0 \ \ \ \textrm{or} \ \ \ \lambda^b = 0.
\end {equation}
Both cases imply that the displacement continuity condition at $d$ becomes a node-to-node constraint. This type of constraint is very convenient since it does not require a master/slave definition; it is two-pass by nature. Furthermore, it can be handled automatically by the assembly procedure without the need for Lagrange multipliers by treating either one of the two nodes as a duplicate that can be eliminated and replaced by the other in the element formulation.  

For ${\lambda^a \neq 0}$ and ${\lambda^b \neq 0}$, we can simplify equation \eqref{eq_DG:DPass_ud2_Q4} further by noting the following:
\begin{align} \notag
\frac{1-\lambda^d}{\lambda^d} &= \left( 1-\frac{\gamma - \alpha}{\beta-\alpha} \right)\frac{\beta-\alpha}{\gamma-\alpha} = \frac{\beta-\gamma}{\gamma-\alpha} \\ \notag
\frac{(1-\lambda^d)\lambda^a}{\lambda^d\,\lambda^b} &=\frac{\gamma - \alpha}{\gamma} \,\frac{\beta-\gamma}{\gamma-\alpha} \, \frac{1-\gamma}{\beta-\gamma}  \\
&= \frac{1-\gamma}{\gamma} \equiv \rho,
\end{align}
which reduces equation \eqref{eq_DG:DPass_ud2_Q4} to the simple expression
\begin{align}
\mathbf{u}^d&= \frac{\rho}{\rho+1}\mathbf{u}^c + \frac{1}{\rho+1}\mathbf{u}^e \\
&= (1-\gamma)\mathbf{u}^c + \gamma \mathbf{u}^e.
\end{align}
Thus, the two elements constituting the foundation become effectively a single bilinear element. This case becomes  similar to Case I, where we have shown that displacement continuity holds point-wise along the interface, and the patch test is passed. This result explains the reported fact that the two-pass MPC approach passes the patch test. The price, however, is the over-constraining of the interface which leads to surface locking. The locking is not reduced by mesh refinement since for each added node on the interface, an additional constraint is activated and the net effect yields no additional flexibility on the interface. As a result, convergence is not achieved in the limit of mesh refinement. 

\subsection {Quadratic elements}
\begin{figure}
	\centering
	\includegraphics[clip]{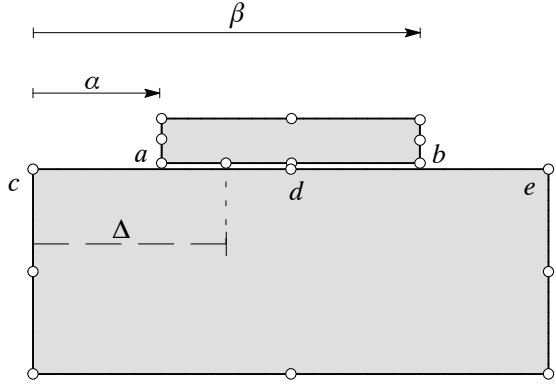}
	\caption{Non-conforming discretization using Q8 elements  \label{fig:DG_MPC_Q8}}
	\end{figure}

To discuss this case we restrict our attention to the simple case shown in Figure \ref{fig:DG_MPC_Q8}, where both the punch and the foundation are represented by a single Q8 element, and the mid-edge node of the punch element matches exactly that of the foundation element. We make this choice to show that, even in the simplest configuration, the drawbacks of the MPC formulation are evident, and we can expect the effect of the issues discussed herein to be more substantial in more general configurations.
The displacement interpolation along the surface of the punch can be expressed as follows
\begin{align} \notag
\mathbf{u}^{+\Delta} &= N_{Q8}^a(\Delta)\,\mathbf{u}^a+ N_{Q8}^b(\Delta)\,\mathbf{u}^b+ N_{Q8}^d(\Delta)\,\mathbf{u}^d \\ \label{eq_DG:UDelta+_Q8}
& + \left[N_{Q4}^a(\Delta) - \frac{1}{2} f^+(\Delta)\right]\mathbf{u}^a +
\left[N_{Q4}^e(\Delta) - \frac{1}{2} f^+(\Delta)\right]\mathbf{u}^e
+ f^+(\Delta)\,\mathbf{u}^d,
\end{align}
where $N_{Q4}^q$ and $N_{Q8}^q$ are the Q4 and Q8 shape function associated with corner node $q$, respectively, and $f^+=N_{Q8}^d$ is the shape function associated with the mid-edge node. Simplifying,
\begin{align} \notag
\mathbf{u}^{+\Delta} &= \mathbf{u}^{+\Delta}_{Q4} + f^+(\Delta) \left[\mathbf{u}^d - \frac{1}{2}(\mathbf{u}^a+\mathbf{u}^b)\right] \\
&= (1-\frac{\Delta - \alpha}{\beta-\alpha})\,\mathbf{u}^a + \frac{\Delta - \alpha}{\beta-\alpha}\,\mathbf{u}^b + f^+(\Delta) \left[\mathbf{u}^d - \frac{1}{2}(\mathbf{u}^a+\mathbf{u}^b)\right],
\end{align}
where $\mathbf{u}^{+\Delta}_{Q4}$ is the bilinear component of the displacement, as given in equation \eqref{eq_DG:UDelta-_Q4}. Similarly, in the foundation,
\begin{align} \notag
\mathbf{u}^{-\Delta} &= \mathbf{u}^{-\Delta}_{Q4} + f^-(\Delta) \left[\mathbf{u}^d - \frac{1}{2}(\mathbf{u}^c+\mathbf{u}^e)\right] \\ \label {eq_DG:UDelta-_Q8}
&= (1-\Delta)\,\mathbf{u}^c +\Delta \mathbf{u}^e +f^-(\Delta)\left[\mathbf{u}^d - \frac{1}{2}(\mathbf{u}^c+\mathbf{u}^e)\right].
\end{align}
It is useful to point out that, even though both $f^+$ and $f^-$ are associated with $d$, they are not the same function since $f^+$ goes to zero at $a$ and $b$, while $f^-$ is zero at $c$ and $e$.
Enforcing the MPCs at $a$ and $b$ leads to
\begin{align} \notag
\mathbf{u}^a &= (1-\alpha)\,\mathbf{u}^c+\alpha\,\mathbf{u}^e +f^-(a)\left[\mathbf{u}^d - \frac{1}{2}(\mathbf{u}^c+\mathbf{u}^e)\right] \\ \label{eq_DG:ua_Q8}
&= (1-\alpha)\,\mathbf{u}^c+\alpha\,\mathbf{u}^e +\Phi^a \\ \notag
\mathbf{u}^b &= (1-\beta)\,\mathbf{u}^c+\beta\,\mathbf{u}^e +f^-(b)\left[\mathbf{u}^d - \frac{1}{2}(\mathbf{u}^c+\mathbf{u}^e)\right] \\  \label{eq_DG:ub_Q8}
&= (1-\beta)\,\mathbf{u}^c+\beta\,\mathbf{u}^e +\Phi^b,
\end{align}
where $\Phi^a = f^-(a)\left[\mathbf{u}^d - \frac{1}{2}(\mathbf{u}^c+\mathbf{u}^e)\right]$ and $\Phi^b = f^-(b)\left[\mathbf{u}^d - \frac{1}{2}(\mathbf{u}^c+\mathbf{u}^e)\right]$. The substitution of equations \eqref{eq_DG:ua_Q8} and \eqref{eq_DG:ub_Q8} into equation \eqref{eq_DG:UDelta+_Q8} yields
\begin{align} \notag
\mathbf{u}^{+\Delta} &= (1-\frac{\Delta - \alpha}{\beta-\alpha})
\left[ (1-\alpha)\,\mathbf{u}^c+\alpha\,\mathbf{u}^e +\Phi^a \right ] \\ 
&+ \frac{\Delta - \alpha}{\beta-\alpha}
\left[ (1-\beta)\,\mathbf{u}^c+\beta\,\mathbf{u}^e +\Phi^b \right] +f^+(\Delta) \left[\mathbf{u}^d - \frac{1}{2}(\mathbf{u}^a+\mathbf{u}^b)\right].
\end{align}
The first two terms in this equation constitute the bilinear component of the displacement field in the punch. From equations \eqref{eq_DG:Delta_simple} and \eqref{eq_DG:UDelta+} for the case of bilinear elements with $\gamma=1$, the following holds
\begin{equation} \label{eq_DG:UQ4inQ8}
(1-\frac{\Delta - \alpha}{\beta-\alpha})\left[(1-\alpha)\,\mathbf{u}^c+\alpha\,\mathbf{u}^e \right ]+ \frac{\Delta - \alpha}{\beta-\alpha} \left[ (1-\beta)\,\mathbf{u}^c+\beta\,\mathbf{u}^e  \right] = (1-\Delta)\,\mathbf{u}^c +\Delta \mathbf{u}^e.
\end{equation}
Given equation \eqref{eq_DG:UQ4inQ8} and the definitions of $\Phi^a$ and $\Phi^b$, we can write the displacement field in the punch for the case of quadratic elements as follows:
\begin{align}\notag
\mathbf{u}^{+\Delta} &= (1-\Delta)\,\mathbf{u}^c+\Delta\,\mathbf{u}^e
+ f^+(\Delta) \left[\mathbf{u}^d - \frac{1}{2}(\mathbf{u}^a+\mathbf{u}^b)\right]
\\ \label{eq_DG:UDelta+_Q8_2} 
& + \left[\mathbf{u}^d - \frac{1}{2}(\mathbf{u}^c+\mathbf{u}^e)\right]\, 
\left[(1-\frac{\Delta - \alpha}{\beta-\alpha})\,f^-(a)+ \frac{\Delta - \alpha}{\beta-\alpha}\,f^-(b)\right].
\end{align}
Comparing equation \eqref{eq_DG:UDelta+_Q8_2} to equation \eqref{eq_DG:UDelta-_Q8} reveals that the displacement fields in the punch and the foundation are not identical, and pointwise displacement continuity is not generally achievable, even in this relatively simple configuration. The cases where the two fields match exactly are limited to the following scenarios

\begin{enumerate}
\item If $\mathbf{u}^d = \frac{1}{2}(\mathbf{u}^a+\mathbf{u}^b)$ and $\mathbf{u}^d = \frac{1}{2}(\mathbf{u}^c+\mathbf{u}^e)$, then $\mathbf{u}^{+\Delta} = \mathbf{u}^{-\Delta}$. In this case, the interpolation of the displacement fields in both the punch and the foundation becomes linear, and the point-wise continuity is achieved. The problem reverts to Case I as discussed above.

\item Another possibility for point-wise continuity is when $f^-(a) = f^-(b) = 0$ and $f^-(\Delta) = f^+(\Delta)$ for any $\Delta$. The only feasible scenario for these conditions to be met is when nodes $a$ and $b$ in the punch match with nodes $c$ and $e$ in the foundation, respectively. In this case, all MPCs become node-to-node constraints and the mesh transforms to a conforming one where the displacement is point-wise continuous by design. It is important to note here that, for quadratic elements, a mesh where all nodes match at the inter-element interfaces is not necessarily a conforming one. For the node-to-node MPCs to result in a conforming mesh, the matching nodes have to be of the same type (mid-node to mid-node and edge node to edge node). Otherwise, the displacement fields between the nodes on either side of the interface do not match due to the difference in interpolating functions associated with each node. For example, if $\mathbf{u}^+ = \sum_i N^i \mathbf{u}^+_i$ and $\mathbf{u}^- = \sum_i M^i \mathbf{u}^-_i$, for the displacement fields to be identical $\mathbf{u}^+ = \mathbf{u}^-$, we need to have $\mathbf{u}_i^+ = \mathbf{u}_i^-$ and $N^i = M^i$ for any node $i$. This obviously is not an issue for bilinear elements since all nodes are edge ones.
\end{enumerate}

In the above section (\ref{sec_DG:mpc}) we have established that continuity of the variational field, at least in a weak sense, is a sufficient condition for the complete transfer of interface tractions. In the context of the MPC, point-wise continuity can be achieved for bilinear elements when the MPCs are enforced at all nodes at the interface (a two-pass approach).  The downside of such approach, however, is the interface locking resulting from the excessive enforcement of the compatibility constraints, which leads to a slow or lack of convergence with mesh refinement. The cure to avoiding locking is to use node-to-node compatibility constraints. This approach has the added benefit of being inherently two-pass and does not require the calculation of Lagrange multipliers as continuity can be assumed by the assembly procedure. This method, however is not guaranteed to pass the patch test for elements of quadratic or higher order. These observations are key to the development of the proposed interface formulation as described below.
%
%
%
%
\section{Proposed interface formulation} \label{sec_DG:prop_form}
The aim of the proposed approach is twofold: to enable a two-pass strategy for the enforcement of geometric compatibility along the interface, and to ensure that the resulting formulation passes the patch test. To achieve the first objective we propose a local enrichment of the interface that would transform the continuity constraints to a set of node-to-node constraints that can be enforced at all nodes along the interface. We then employ a discontinuous Galerkin-based stabilization procedure to ensure a complete transfer of the traction field across the interface. The enrichment and stabilization procedures are described in what follows. 
%
%
\subsection{Enrichment procedure} \label{sec_DG:enr_procedure}
	\begin{figure}
	\centering
	\includegraphics[clip]{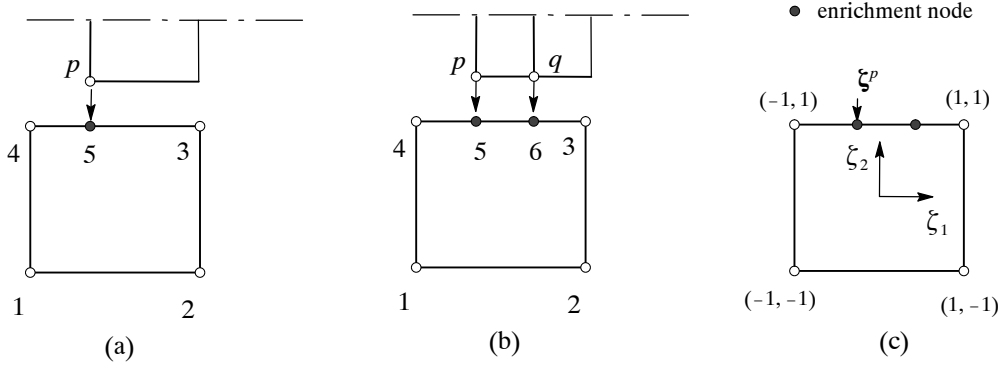}
	\caption{Local enrichment of the interface element for the cases of (a) a single and (b) multiple nodes, and (c) added node reference in the parent domain  \label{fig:DG_enr}}
	\end{figure}
The purpose of the enrichment is to transform the displacement continuity condition
at each interface node to a simple node-to-node constraint without modifying the existing mesh. To achieve this purpose, we adopt
the following approach: At every location where a node $p$ meets an element surface, a node can be
inserted such that the displacement continuity at that location is reduced to a node-to-node constraint, as shown in Figure \ref{fig:DG_enr}(a). Completeness of the finite element interpolation in the enriched element can
be preserved by updating the set of Lagrangian shape functions to account for the additional node. For
example, consider the case illustrated in Figure \ref{fig:DG_enr}(a). Given $\mathbf\zeta^{p}$, the reference of point $p$ in the element
parent coordinates $\mathbf\zeta$ shown in Figure \ref{fig:DG_enr}(c), the updated element shape functions are,
		\begin{align} 
		\tilde{N}^{5}(\mathbf{\zeta},\mathbf{\zeta}^{p}) &= \textstyle\frac{1}{2}(\zeta_{1}+1)
		\displaystyle\frac{(\zeta_{2}+1)(\zeta_{2}-1)}{(\zeta_{2}^{p}+1)(\zeta_{2}^{p}-1)} 
		\\  \label{eq_DG:N_enr}
		\tilde{N}^{\alpha} &= N_{Q4}^{\alpha}-N_{Q4}^{\alpha}(\mathbf{\zeta}^{p})\tilde{N}^{5} 
		\end{align}
for $\alpha =1,...,4$ where $N_{Q4}^{\alpha}$
are the shape functions of a Q4 element. For an element $m$ defined by a set of nodes with spatial coordinates $\mathbf{x}^{\alpha}$ and corresponding shape functions $N^{\alpha}$, where $\alpha = 1,...,n$,
The enriched finite element formulation of the spatial variable $\mathbf{x}$ in the element is given by
		\begin{equation} \label{eq_DG:xEnrIso}
		\mathbf{x}(\mathbf{\zeta},\mathbf{\zeta}^{p}) = \tilde{N}^{\alpha}(\mathbf{\zeta},\mathbf{\zeta}^{p}) \mathbf{x}^{\alpha} 
		+ \tilde{N}^{p}(\mathbf{\zeta},\mathbf{\zeta}^{p})\mathbf{x}^{p} , \ \ \  \alpha = 1,\cdots,n,
		\end{equation}
where $\tilde{N}^{p}$ is the shape function associated with the additional node and the $\tilde{N}^{\alpha}$ are the modified (enriched) element shape functions defined as follows
		\begin{equation} \label{eq_DG:NEnr}
		\tilde{N}^{\alpha}(\mathbf{\zeta},\mathbf{\zeta}^{p}) = 
		N^{\alpha}(\mathbf{\zeta})-N^{\alpha}(\mathbf\zeta^{p})\tilde{N}^{p}(\mathbf{\zeta},\mathbf{\zeta}^{p}),\ \ \ \alpha = 1,\cdots,n.
		\end{equation}
The spatial coordinate of the additional node corresponds to that of node $p$. This procedure can be repeated
with multiple nodes for a given element, as shown in Figure \ref{fig:DG_enr}(b). The case of multiple enrichments can be handled using a recursive approach in which one enrichment is performed at a time. Each added node has the effect of
increasing the order of interpolation along the interface, thus providing the surface with an additional
degree of freedom in each spatial direction. This approach can therefore be safely implemented as a
two-pass method. The enrichment applies only to the shape functions that are non-zero at the interface ($N^3$ and $N^4$ in the given case), and the order of interpolation along all other interfaces remains the same, thereby preserving the conformity of the mesh at these interfaces. The enriched formulation is local to the interface element, and has no effect on the
finite element discretization in the remainder of the mesh. Moreover, the nodal displacement continuity
constraints can be accommodated automatically by the assembly procedure without the need for additional
variables or Lagrange multipliers. 

In order to preserve the isoparametric formulation in the enriched element, the set of shape functions $\tilde{N}^{\alpha} , \tilde{N}^{p}$ can be used to define the interpolation of the material and displacement fields in $m$, such that
		\begin{align} 
		\label{eq_DG:XEnrIso}
		\mathbf{X}(\mathbf{\zeta}) &= \tilde{N}^{\alpha}(\mathbf{\zeta},\mathbf{\zeta}^{p}) \mathbf{X}^{\alpha} 
		+ \tilde{N}^{p}(\mathbf{\zeta},\mathbf{\zeta}^{p})\mathbf{X}^{p}  \ \ \  \alpha = 1,\cdots,n.
		\end{align}
In these equations, $\mathbf{X}^{p}$ denotes the material point corresponding to $\mathbf{x}^{p}$. We assume that the boundary of the material domain remains fixed during the enrichment procedure and therefore the mapping between the reference and material domains remains unaffected by the enrichment. As a result,
		\begin{align} 
		\notag
		N^{\alpha} \mathbf{X}^{\alpha} &= \tilde{N}^{\alpha}(\mathbf{\zeta},\mathbf{\zeta}^{p}) \mathbf{X}^{\alpha} 
		+ \tilde{N}^{p}(\mathbf{\zeta},\mathbf{\zeta}^{p})\mathbf{X}^{p}  \\ \label{eq_DG:XEnrIso2}
		&= \left[ N^\alpha(\mathbf{\zeta}) - N^\alpha(\mathbf{\zeta}^p)\tilde{N}^{p} \right] \mathbf{X}^\alpha + \tilde{N}^{p}\,\mathbf{X}^p  \ \ \  \alpha = 1,\cdots,n.
		\end{align}
To wit,
		\begin{equation} 
		\tilde{N}^{p}\left[\mathbf{X}^p - N^\alpha(\mathbf{\zeta}^p)\mathbf{x}^\alpha\right] = \mathbf{0}.
		\end{equation}
Therefore, the material point corresponding to the added node $\mathbf{X}^{p}$ satisfies the following equation:
		\begin{equation} \label{eq_DG:Xp}
		\mathbf{X}^{p} = N^{\alpha}(\mathbf\zeta^{p})\mathbf{X}^{\alpha},
		\end{equation}
Finally, from equations \eqref{eq_DG:xEnrIso} and \eqref{eq_DG:XEnrIso}, the displacement field in the enriched element can be written as
		\begin{align} 
		\label{eq_DG:uEnrIso}	
		\mathbf{u}(\mathbf{\zeta},\mathbf{\zeta}^{p}) &= \tilde{N}^{\alpha}(\mathbf{\zeta},\mathbf{\zeta}^{p}) \mathbf{d}^{\alpha} 
		+ \tilde{N}^{p}(\mathbf{\zeta},\mathbf{\zeta}^{p})\mathbf{d}^{p} , \ \ \  \alpha = 1,\cdots,n	
		\end{align}

To complete the element formulation, we need to compute the reference $\mathbf\zeta^{p}$ of the enrichment node with respect to the element parent coordinates, as is described in the following section.
%
%
\subsubsection{Computing the enrichment location $\mathbf\zeta^{p}$}
Let $p$ be the node meeting the surface  $\zeta_{j}=c$ of a given element $m$. Let the set of shape functions $N^{\alpha}$ define the isoparametric interpolation of the material and spatial variables in $m$, such that
		\begin{equation} 
		\mathbf{x}(\mathbf\zeta) = N^{\alpha}(\mathbf\zeta)\mathbf{x}^{\alpha}, \ 
		\mathbf{X}(\mathbf\zeta) = N^{\alpha}(\mathbf\zeta)\mathbf{X}^{\alpha},
				   \ \ \ \alpha = 1,\cdots,n
		\end{equation}
To compute the reference of point $p$ in element $m$, we use the same technique discussed in section \ref{sec_con:Zcoord}, by mapping $p$ to the reference coordinates of the (unenriched) element. Therefore, the enrichment reference $\mathbf{\zeta}^{p}$ corresponds to the solution of the system of nonlinear equations
\begin{equation}  \label{eq_DG:enr_loc}
\mathbf{x}^p = N^\alpha\left(\mathbf{\zeta}^p\right)\mathbf{x}^\alpha
\end{equation}
for $\mathbf{\zeta}^p$. To restrict the added node to the surface of the enriched element and minimize overlap, we impose the condition $\zeta_{j}^p=c$. Therefore, Equation \eqref{eq_DG:enr_loc} is in fact a system of \textit{ndof}-1 equations to be solved for $\zeta_{i}^{p}$ where \textit{ndof} denotes the dimensionality of the problem (2D vs 3D) and $i \neq j$.

In the case of multiple enrichments, the reference of the $n$-th enrichment is computed by mapping the corresponding node $p_n$ to the element defined by its original nodes and $n-1$ enrichments. The computational cost associated with this procedure is directly proportional to the order of interpolation in the element. Therefore, the solution can prove costly for multiple enrichments, particularly in the presence of large deformations. To simplify the solution procedure, we make use of the following observation: 

Given that the initial discretization is a Lagrangian mesh, the material point associated with node $p$ remains the same through deformation. In the physical domain, the material point associated with node $p$ is attached to the surface of the element, or, more precisely, to the material point in the element opposing $p$ across the interface (the mesh does not exist in the physical domain). Similarly, the mapping of the element boundary to the material domain remains unchanged during deformation and the mesh/material domain mapping is not affected by the enrichment, as described by equation \eqref{eq_DG:XEnrIso2}. Therefore, it is safe to assume that the material point associated with the enrichment does not change. As a result, the solution to equation \eqref{eq_DG:enr_loc} can be executed in the material domain, even in the case of multiple enrichments. The material point corresponding to each enrichment $q$ to the element is therefore computed by solving the system
\begin{equation}  \label{eq_DG:enr_loc_X}
\mathbf{X}^q = N^\alpha\left(\mathbf{\zeta}^q\right)\mathbf{X}^\alpha
\end{equation}
for each enrichment $\mathbf{\zeta}^q$. The computational cost associated with the solution is minimal and the procedure is very efficient for the case of multiple enrichments. 
%
%
%
%
\subsection{Stabilized interface formulation}  
The enrichment procedure discussed above ensures continuity of the displacement field at all nodal locations along the interface. This however, does not generally imply full geometric compatibility as the displacement field remains discontinuous between the nodes. This effect is amplified by the increased interpolation order dictated by the enrichment procedure. As discussed in Section \ref{sec_DG:patch_test}, within the context of the standard Galerkin method, the discontinuity in the kinematic field leads to an incomplete transfer of the traction field across the interface. This suggests treating the interface using the discontinuous Galerkin method.

In a typical DG formulation, the displacement field is interpolated independently in each element and continuity is enforced in an average sense along the element interfaces. Therefore, DG methods are very computationally involved. Besides, Arnold et al. \cite{arnoldSIAMJNA02} have shown that not all DG formulations are stable.

A full DG method would not be suitable for the case at hand because of two reasons. Firstly, the displacement discontinuity is restricted to the non-matching interface. The remainder of the mesh can still be treated using the continuous Galerkin formulation. Secondly, it does not take advantage of the continuity of displacement at the interface nodes, a fact that could greatly reduce computational cost. Therefore, a special DG formulation can be designed to better accommodate these particular features of the problem. This formulation is described in the following.
	\begin{figure}
	\centering
	\includegraphics[clip]{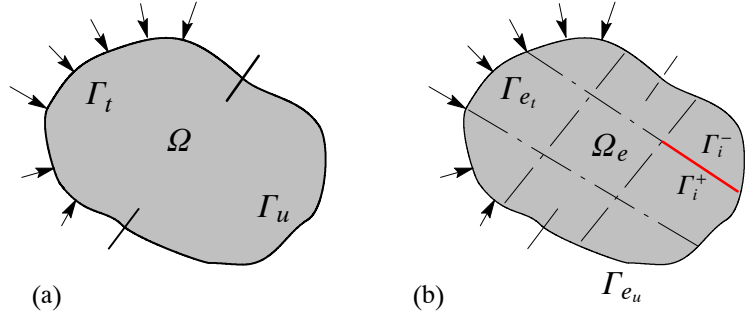}
	\caption{(a) Original configuration and (b) finite element discretization of a domain $\Omega$  \label{fig:DG_FEM}}
	\end{figure}
Consider a discretization of the domain $\Omega$ into a set of finite elements $\Omega_e$ such that $\Omega = \sum_e\Omega_e$, as shown in Figure \ref{fig:DG_FEM}. Defining the finite element subspaces $S^h=\left\{\mathbf{u}^h \in H^1(\Omega_e):\mathbf{u}^h=\mathbf{g} \ \textrm{on} \ {\Gamma_e}_{u}\right\}$ and $V^h=\left\{\mathbf{w}^h \in H^1(\Omega_e):\mathbf{w}^h=\mathbf{0} \ \textrm{on} \ {\Gamma_e}_{u}\right\}$, we can rewrite equation \eqref{eq_DG:Wequil} as follows
\begin{equation} \label{eq_DG:Gequil}
\sum_e\int_{\Omega_e}\left[\mathrm{div}\mathbf{S}+\mathbf{b}\right]\cdot\mathbf{w}^{h}\,d \Omega + \sum_e\int_{\Gamma_t\cap\Gamma_e}\left[\mathbf{h}-\mathbf{Sn}\right]\cdot\mathbf{w}^{h}\,d \Gamma = 0 \ \ \ \forall \ \mathbf{w}^h \in V^h,
\end{equation}
where $\Gamma_e$ is the element boundary. Applying the divergence theorem to the first term in equation \eqref{eq_DG:Gequil} yields
\begin{equation} \label{eq_DG:Gequil_div}
\sum_e\int_{\Omega_e}\left[-\mathbf{S}\cdot\nabla\mathbf{w}^{h}+\mathbf{b}\cdot\mathbf{w}^{h}\right]\,d \Omega +
\sum_e\int_{\Gamma_e}\mathbf{Sn}\cdot\mathbf{w}^h\,d\Gamma+\sum_e\int_{\Gamma_t\cap\Gamma_e}\left[\mathbf{h}-\mathbf{Sn}\right]\cdot\mathbf{w}^{h}\,d\Gamma=0.
\end{equation}
Let $\Gamma_i$ be the set of inter-element interfaces. The boundary of a given element $\Gamma_e$ can be decomposed into $\Gamma_e={\Gamma_e}_{t}\cup{\Gamma_e}_{u}\cup{\Gamma_e}_{i}$ where
\begin{equation}
{\Gamma_e}_{t} = {\Gamma_e}\cap{\Gamma_t}, \ \ \ {\Gamma_e}_{u} = {\Gamma_e}\cap{\Gamma_u}, \ \ \ {\Gamma_e}_{i} = {\Gamma_e}\cap{\Gamma_i}
\end{equation}
and
\begin{equation}
{\Gamma_e}_{t}\cap{\Gamma_e}_{i}={\Gamma_e}_{t}\cap{\Gamma_e}_{u}={\Gamma_e}_{u}\cap{\Gamma_e}_{i}= \phi.
\end{equation}
Therefore, equation \eqref{eq_DG:Gequil_div} reduces to
\begin{equation} \label{eq_DG:Gequil_div_i}
-\sum_e\int_{\Omega_e}\mathbf{S}\cdot\nabla\mathbf{w}^{h}\,d \Omega+\sum_e\int_{\Omega_e}\mathbf{b}\cdot\mathbf{w}^{h}\,d \Omega +
\sum_e\int_{{\Gamma_e}_t}\mathbf{h}\cdot\mathbf{w}^{h}\,d\Gamma+\sum_e\int_{{\Gamma_e}_i}\mathbf{Sn}\cdot\mathbf{w}^h\,d\Gamma=0,
\end{equation}
where $\sum_{e}\int_{\Gamma_e}\mathbf{Sn}\cdot\mathbf{w}^h d \Gamma$ is the element interface term. Let $\mathbf{n}^+$, $\mathbf{w}^+$, and $\mathbf{n}^-$, $\mathbf{w}^-$, be the normal vector and variational displacement vector on each side of interface $i$. The interface term in equation \eqref{eq_DG:Gequil_div_i} can be rewritten as follows
\begin{equation} \label{eq_DG:jump_i}
\sum_e\int_{{\Gamma_e}_i}\mathbf{Sn}\cdot\mathbf{w}^h\,d\Gamma = \sum_i\int_{\Gamma_i^+}\mathbf{Sn}^+\cdot\mathbf{w}^{h+}\,d\Gamma+\sum_i\int_{\Gamma_i^-}\mathbf{Sn}^-\cdot\mathbf{w}^{h-}\,d\Gamma.
\end{equation}
Therefore, the Galerkin form of the governing equations becomes
(G) Given the body force vector $\mathbf{b}$, the prescribed traction field $\mathbf{h}$ on $\Gamma_t$ and the prescribed displacement
field $\mathbf{g}$ on $\Gamma_u$ find $\mathbf{u}^h\in S^h =\left\{\mathbf{u}^h\in H^1(\Omega_e):\mathbf{u}^h=\mathbf{g} \ \textrm{on} \ {\Gamma_e}_u\right\}$ such that
\begin{align} \notag
-\sum_e\int_{\Omega_e}\mathbf{S}\cdot\nabla\mathbf{w}^{h}\,d \Omega&+\sum_e\int_{\Omega_e}\mathbf{b}\cdot\mathbf{w}^{h}\,d \Omega +
\sum_e\int_{{\Gamma_e}_t}\mathbf{h}\cdot\mathbf{w}^{h}\,d\Gamma \\ \label{eq_DG:Gequil_div_i2}
&+\sum_i\int_{\Gamma_i^+}\mathbf{Sn}^+\cdot\mathbf{w}^{h+}\,d\Gamma+\sum_i\int_{\Gamma_i^-}\mathbf{Sn}^-\cdot\mathbf{w}^{h-}\,d\Gamma=0,
\end{align}
for all $\mathbf{w}^h\in V^h =\left\{\mathbf{w}^h\in H^1(\Omega_e):\mathbf{w}^h=\mathbf{0} \ \textrm{on} \ {\Gamma_e}_u\right\}$. 

The continuous Galerkin formulation can be derived as a special case of equation \eqref{eq_DG:Gequil_div_i2}. If $\mathbf{w}^h$ is continuous across the element boundaries, the interface term becomes
\begin{equation} \label{eq_DG:Gequil_continuous}
\sum_i\int_{\Gamma_i^+}\mathbf{Sn}^+\cdot\mathbf{w}^{h+}\,d\Gamma+\sum_i\int_{\Gamma_i^-}\mathbf{Sn}^-\cdot\mathbf{w}^{h-}\,d\Gamma =
\sum_i\int_{\Gamma_i}\left[\mathbf{S}^+ - \mathbf{S}^-\right]\mathbf{n}\cdot\mathbf{w}^h\,d\Gamma,
\end{equation}
which weakly enforces the equilibrium of tractions along the element interfaces. If $\mathbf{w}^h$ is not continuous along the element interfaces, equilibrium of tractions does not necessarily hold and the formulation is unstable.

To stabilize the solution, we subtract the weighted residual of the interface tractions to equation \eqref{eq_DG:Gequil_div_i2}. To wit,
\begin{align} \notag
&-\sum_e\int_{\Omega_e}\mathbf{S}\cdot\nabla\mathbf{w}^{h}\,d \Omega+\sum_e\int_{\Omega_e}\mathbf{b}\cdot\mathbf{w}^{h}\,d \Omega +
\sum_e\int_{{\Gamma_e}_t}\mathbf{h}\cdot\mathbf{w}^{h}\,d\Gamma +\sum_i\int_{\Gamma_i^+}\mathbf{Sn}^+\cdot\mathbf{w}^{h+}\,d\Gamma \\ \label{eq_DG:Gequil_stabilized}
&+\sum_i\int_{\Gamma_i^-}\mathbf{Sn}^-\cdot\mathbf{w}^{h-}\,d\Gamma-\sum_i\int_{\Gamma_i}\left[\mathbf{Sn}^+ + \mathbf{Sn}^-\right]\cdot\bar{\mathbf{w}}^h\,d\Gamma=0 \ \ \ \forall \ \mathbf{w}^h \in V^h,
\end{align}
where $\bar{\mathbf{w}}^h=\left(\mathbf{w}^{h+} + \mathbf{w}^{h-} \right)/2$ is the average of the variational displacement along the interface. This choice guarantees an unbiased method. Simplifying,
\begin{align} \notag
&-\sum_e\int_{\Omega_e}\mathbf{S}\cdot\nabla\mathbf{w}^{h}\,d \Omega+\sum_e\int_{\Omega_e}\mathbf{b}\cdot\mathbf{w}^{h}\,d \Omega +
\sum_e\int_{{\Gamma_e}_t}\mathbf{h}\cdot\mathbf{w}^{h}\,d\Gamma  \\ \label{eq_DG:Gequil_stabilized_simple}
&+\sum_i\int_{\Gamma_i^+}\mathbf{Sn}^+\cdot\left[\mathbf{w}^{h+}-\bar{\mathbf{w}}^h\right]\,d\Gamma +\sum_i\int_{\Gamma_i^-}\mathbf{Sn}^-\cdot\left[\mathbf{w}^{h-}-\bar{\mathbf{w}}^h\right]\,d\Gamma=0 \ \ \ \forall \ \mathbf{w}^h \in V^h,
\end{align}
or, equivalently,
\begin{align} \notag
&-\sum_e\int_{\Omega_e}\mathbf{S}\cdot\nabla\mathbf{w}^{h}\,d \Omega+\sum_e\int_{\Omega_e}\mathbf{b}\cdot\mathbf{w}^{h}\,d \Omega +
\sum_e\int_{{\Gamma_e}_t}\mathbf{h}\cdot\mathbf{w}^{h}\,d\Gamma  \\  \notag
&+\frac{1}{2}\sum_i\int_{\Gamma_i^+}\mathbf{Sn}^+\cdot\left[\mathbf{w}^{h+}-\mathbf{w}^{h-}\right]\,d\Gamma \\ \label{eq_DG:Gequil_stabilized_simple2} &+\frac{1}{2}\sum_i\int_{\Gamma_i^-}\mathbf{Sn}^-\cdot\left[\mathbf{w}^{h-}-\mathbf{w}^{h+}\right]\,d\Gamma=0 \ \ \ \forall \ \mathbf{w}^h \in V^h.
\end{align}
It should be clear that when $\mathbf{w}^{h+}=\mathbf{w}^{h-}=\bar{\mathbf{w}}^h$ the interface term goes to zero and the formulation reverts back to the standard continuous Galerkin method. This feature reinforces the status of the proposed method as a consistent framework that includes the continuous Galerkin method as a subset. As evidenced by equation \eqref{eq_DG:Gequil_stabilized_simple2}, when the variational displacement field is point-wise continuous along the interface, as is the case in conforming meshes, equilibrium of interface tractions is automatically satisfied and the standard (unstabilized) Galerkin formulation passes the patch test. 

The differences between the proposed formulation and standard DG formulations can be summarized
as follows:

\begin{enumerate}
\item In a full DG formulation, compatibility of displacement has to be enforced in a weak sense
along the element interfaces. This condition is replaced here by the strong enforcement of
displacement continuity at the nodes. This approach leads to a big reduction in the storage
and computation cost usually associated with DG.
\item DG formulations typically replace the interface tractions $\mathbf{Sn}$ with \textit{numerical fluxes} that
often include a user-defined stabilization parameter. The stabilization term proposed herein
comes from the enforcement of equilibrium between elements and does not involve user input.
\item The formulation is consistent and includes the continuous Galerkin method as a subset.
Therefore, it can be easily integrated in a standard Finite Element code.
\item The only disadvantage of this approach is that the consistent tangent matrix is not symmetric.
This feature is currently under investigation for possible improvement.
\end{enumerate}

\section{Comments on implementation} \label{sec_DG:implementation}
\begin{itemize}
\item {\textbf{Alternative large-deformation formulation}: We presented the development of the interface formulation in an updated Lagrangian framework where all fields are defined in the deformed configuration of the body. The equivalent of the stabilized interface formulation \eqref{eq_DG:Gequil_stabilized_simple2} in a total Lagrangian framework can be written as follows
\begin{align} \notag
&-\sum_e\int_{V_e}\mathbf{P}\cdot\nabla_{\mathbf{X}}\mathbf{w}^{h}\,dV+\sum_e\int_{V_e}\mathbf{b}_0\cdot\mathbf{w}^{h}\,d \Omega +
\sum_e\int_{{\Gamma_e}_t}\mathbf{h}_0\cdot\mathbf{w}^{h}\,d\Gamma  \\  \notag
&+\frac{1}{2}\sum_i\int_{A_{i}^+}\mathbf{PN}^+\cdot\left[\mathbf{w}^{h+}-\mathbf{w}^{h-}\right]\,dA \\ \label{eq_DG:Gequil_stabilized_PF} &+\frac{1}{2}\sum_i\int_{A_{i}^-}\mathbf{PN}^-\cdot\left[\mathbf{w}^{h-}-\mathbf{w}^{h+}\right]\,dA=0 \ \ \ \forall \ \mathbf{w}^h \in V^h,
\end{align}
where $A_{i}$ is the inter-element interface in the undeformed configuration.}
\item {\textbf{3D interfaces}: The implementation of equation \eqref{eq_DG:Gequil_stabilized_simple2} requires defining the inter-element interfaces, which involves finding the intersections between the surfaces of each two elements that meet at the interface. This process is relatively simple in 2D since the interfaces are one-dimensional. In 3D the interfaces are 2D surfaces of arbitrary shape, as shown in Figure \ref{fig:DG_imp} (a). Integrating the interface terms for such a case requires developing a special integration rule that can handle arbitrary shapes. This issue is universal to coupling algorithms and has been addressed in various ways in the literature. The most common procedure \cite{pusoCMAME04, kimCMAME08} is to subdivide the surface into a set of triangles for which the Gauss integration rule can be used. In this dissertation, we have restricted our attention to planar problems to simplify the coding process. The extension to 3D would only require coding the interface-detection algorithm and does not involve any change in the formulation. We plan to address three-dimensional problems in the future.}
	\begin{figure}
	\centering
	\includegraphics[clip]{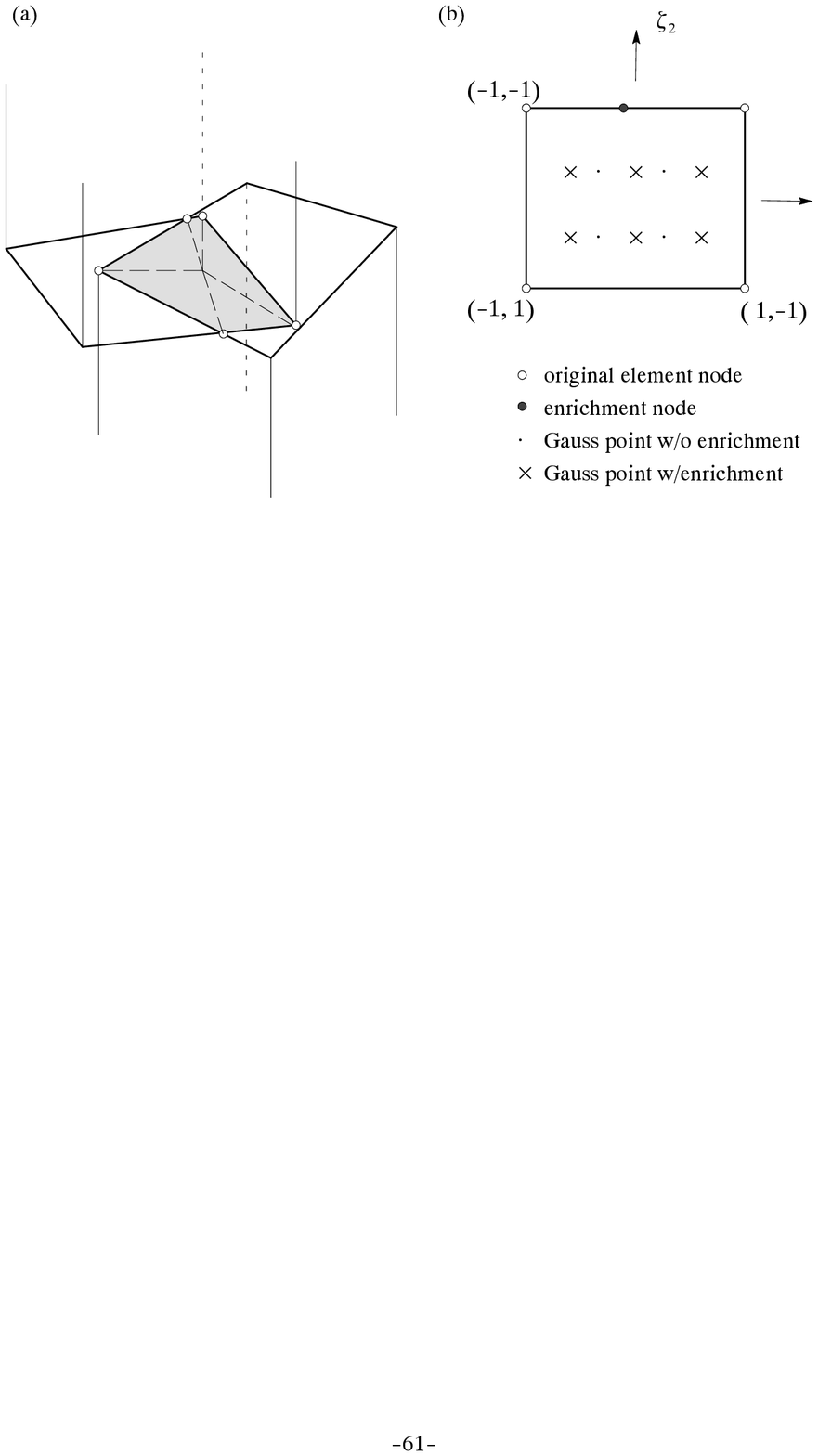}
	\caption{(a) intersecting surfaces in 3D, (b) Gauss integration of the enriched element  \label{fig:DG_imp}}
	\end{figure}
\item{\textbf{Numerical integration}: Consider the enriched element shown in Figure \ref{fig:DG_imp} (b). The enrichment of the top surface $\zeta_2 = 1$ introduces a quadratic term in $\zeta_1$ in the element shape functions (equation \eqref{eq_DG:N_enr}) associated with the nodes located on this interface, while the order of interpolation with respect to $\zeta_2$ remains the same. Therefore, for the element to be integrated properly, the order of the integration rule has to be increased in the direction of $\zeta_1$. When using a Gauss quadrature, one can simply add one integration point in one direction while keeping the number of integration points in the other the same, as shown in Figure \ref{fig:DG_imp} (b). This approach works well for hyperelasticity. For history-dependent materials a progressive integration rule such as the Gauss-Kronrod quadrature can be used alternatively to preserve the computational history at the integration points before enrichment, where necessary.}

\end{itemize}
%
%
%
%
\section{Numerical results} \label{sec_DG:results}
\subsection{Patch test}
The punch and foundation depicted in Figure 2 (a) are made of a linear elastic material with properties $E = 10^5$ and $\nu = 0.3$. A distributed load of $q = 0.1$ is applied to the free surface of the structure. The punch and foundation are discretized using Q4 elements.

\begin{figure}
	\centering
	\includegraphics[clip]{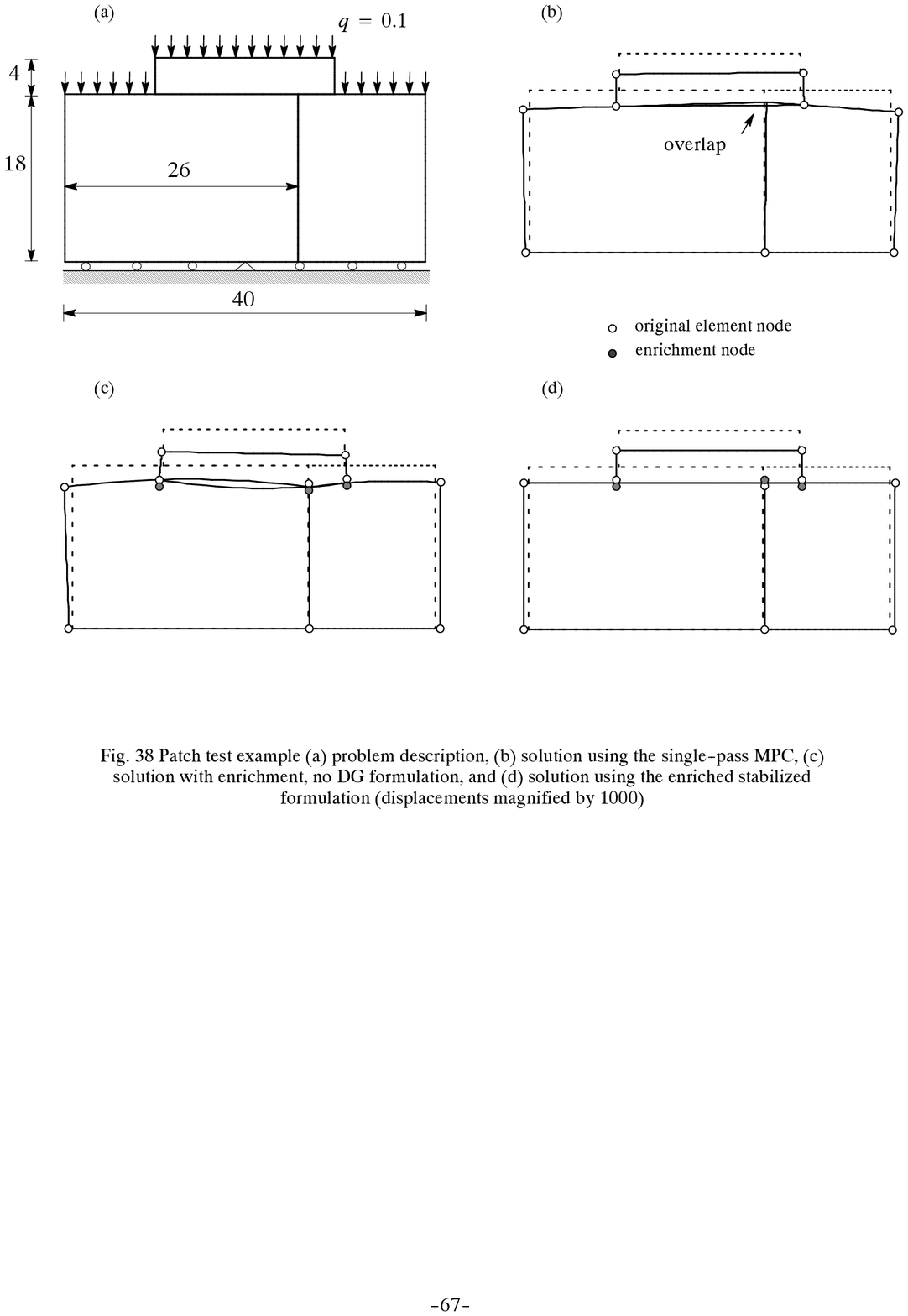}
	\caption{Patch test example (a) problem description, (b) solution using the single-pass MPC, (c) solution with enrichment, no DG formulation, and (d) solution using the enriched stabilized formulation (displacements magnified by 1000)  \label{fig:DG_patchex1}}
	\end{figure}

Figure \ref{fig:DG_patchex1} shows the result of the MPC method, with the foundation designated as the master surface. The solution is clearly not accurate, and involves a foundation node overlapping with the punch. 
Implementing the node-enrichment procedure discussed above (see Figure \ref{fig:DG_enr}) enables us to use a two-pass approach and enforce continuity at all three nodes at the interface. The solutions obtained with and without the DG interface stabilization are shown in Figures \ref{fig:DG_patchex1}(b) and Figure \ref{fig:DG_patchex1}(c), respectively. It is clear from Figures \ref{fig:DG_patchex1}(b) that
the pressure distribution along the interface is not accurate and the patch test is not passed. The deformed
configuration shown in Figures \ref{fig:DG_patchex1}(c) confirms that the stabilized DG formulation is successful
in projecting the pressure field accurately across the interface and the proposed formulation passes the patch test up to machine precision.
	\begin{figure}
	\centering
	\includegraphics[clip]{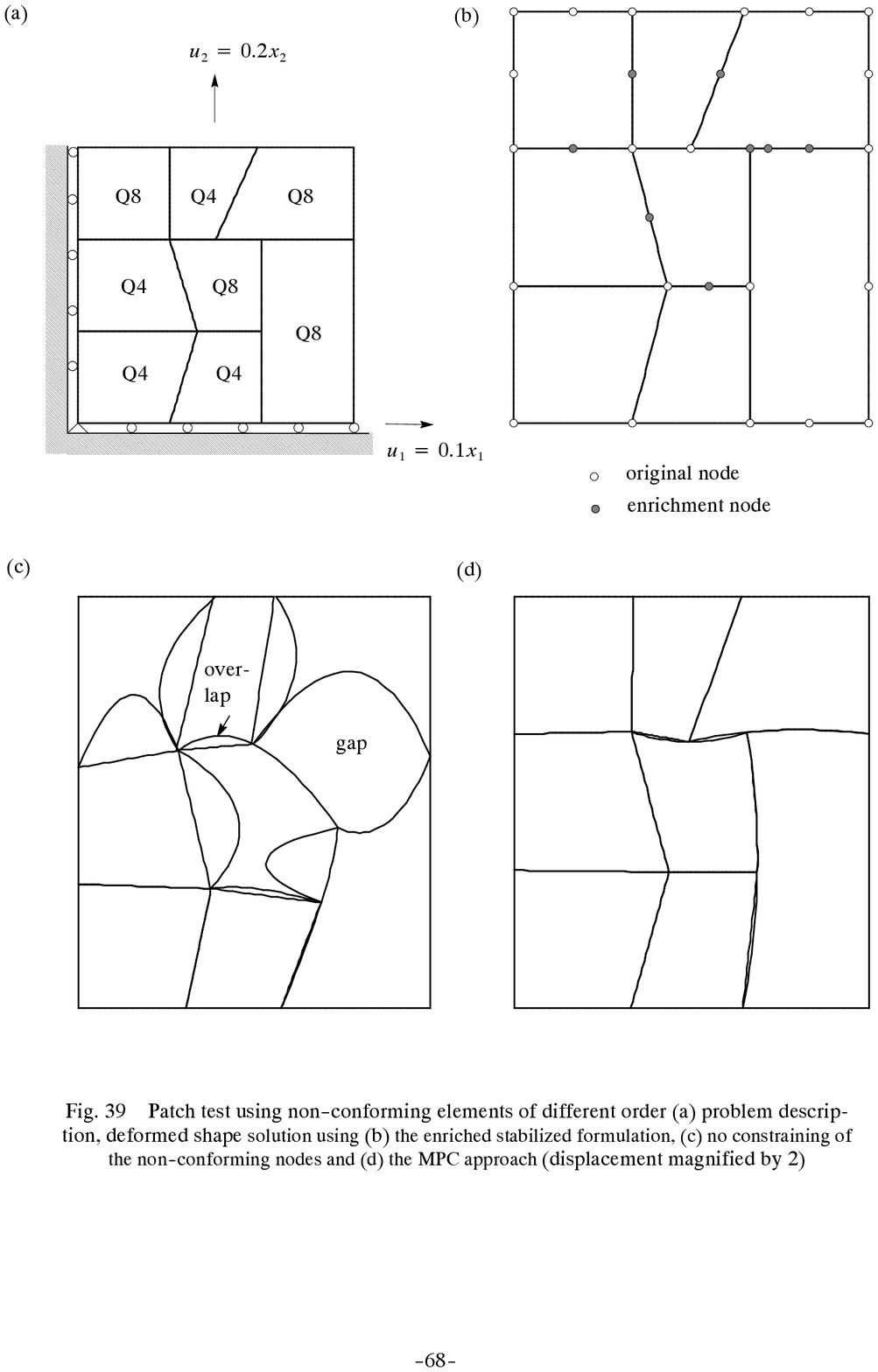}
	\caption{Patch test using non-conforming elements of different order (a) problem description, deformed shape solution using (b) the enriched stabilized formulation, (c) no constraining of the non-conforming nodes and (d) the MPC approach. (The displacement were magnified by 2 in (b),(c) and (d)))  \label{fig:DG_patchex2}}
	\end{figure}

Figure \ref{fig:DG_patchex2} (a) shows a different patch test where a square domain of unit side length is subjected to a prescribed displacement field at its right and top edges. The bottom and left boundaries are connected to roller supports and the material properties are $E = 10^4$ and $\nu = 0.3$. The domain is discretized using different types of elements under plane strain conditions, thereby creating a variety of non-conforming element interfaces. To verify the robustness of the stabilization procedure, the mesh is designed to include elements of non-constant Jacobian. 
Due to the different interpolation order on the interfaces between Q8 and Q4 elements, this configuration contains a number of \textit{floating nodes} that enforce the quadratic interpolation on the Q8 side of the interface.

Figure \ref{fig:DG_patchex2} (c) shows the solution obtained without treating the non-conforming interfaces. The presence of floating nodes leads to large gaps and overlaps in the deformed configuration. Figure \ref{fig:DG_patchex2} (d) shows the result of enforcing the MPCs at the floating nodes. The gaps and overlaps are greatly reduced but the patch test is still failed. The deformed configuration using the proposed formulation is shown in Figure \ref{fig:DG_patchex2} (b) and confirms that the formulation passes the patch test up to machine precision. The same result was obtained using the DG formulation only, i.e. without implementing the enrichment procedure and continuity conditions on the floating nodes. This result confirms the accuracy of the projection of the traction field across the interface by the DG formulation.
%
	\begin{figure}
	\centering
	\includegraphics[clip]{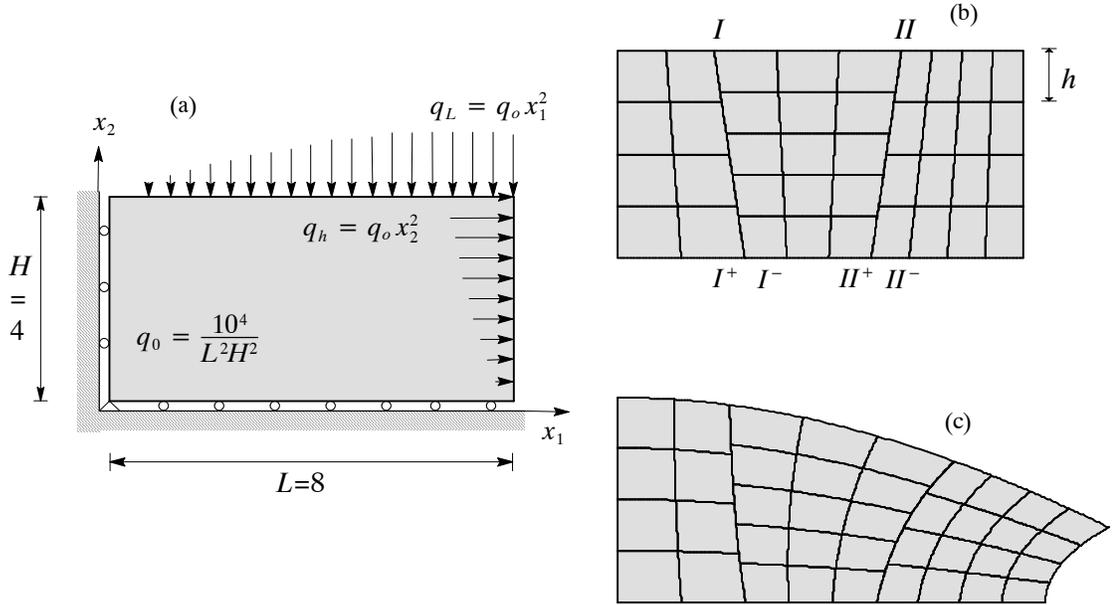}
	\caption{Example with multiple non-conforming interfaces (a) problem description, (b) mesh configuration (c) deformed shape (displacements magnified by 1000)  \label{fig:DG_ex_mesh}}
	\end{figure}

\subsection{Convergence analysis}

Figure \ref{fig:DG_ex_mesh}(a) shows a rectangular domain subjected to a quadratically varying traction field along its
right and top edges. The bottom and left boundaries are connected to roller supports. The domain is made of a linear elastic material with properties $E=10^6$ and $\nu = 0.25$. We discretize the domain using three different layers of elements,
thereby creating two non-conforming interfaces, \textit{I} and \textit{II}, as shown in Figure \ref{fig:DG_ex_mesh}(b). To verify the robustness
of the stabilization procedure, the mesh is designed to include elements of non-constant Jacobian.
The exact solution for this problem is 
\begin{align} \label{eq_DG:exu1}
u_1 &= q_0/E \left[ \nu x_1^3/3+x_1 x_2^2\right] \ \ \ u_2 = -q_0/E \left[ x_2 x_1^2 + \nu x_2^3/3\right] \\
\sigma_{11} &= q_0\,x_2^2  \ \ \ \sigma_{22} = - q_0\,x_1^2  \ \ \  \sigma_{12} = 0.
\end{align}
We designed this problem specifically without shear stress in order to eliminate the effects of shear locking in the Q4 elements.
	\begin{figure}
	\centering
	\includegraphics[clip]{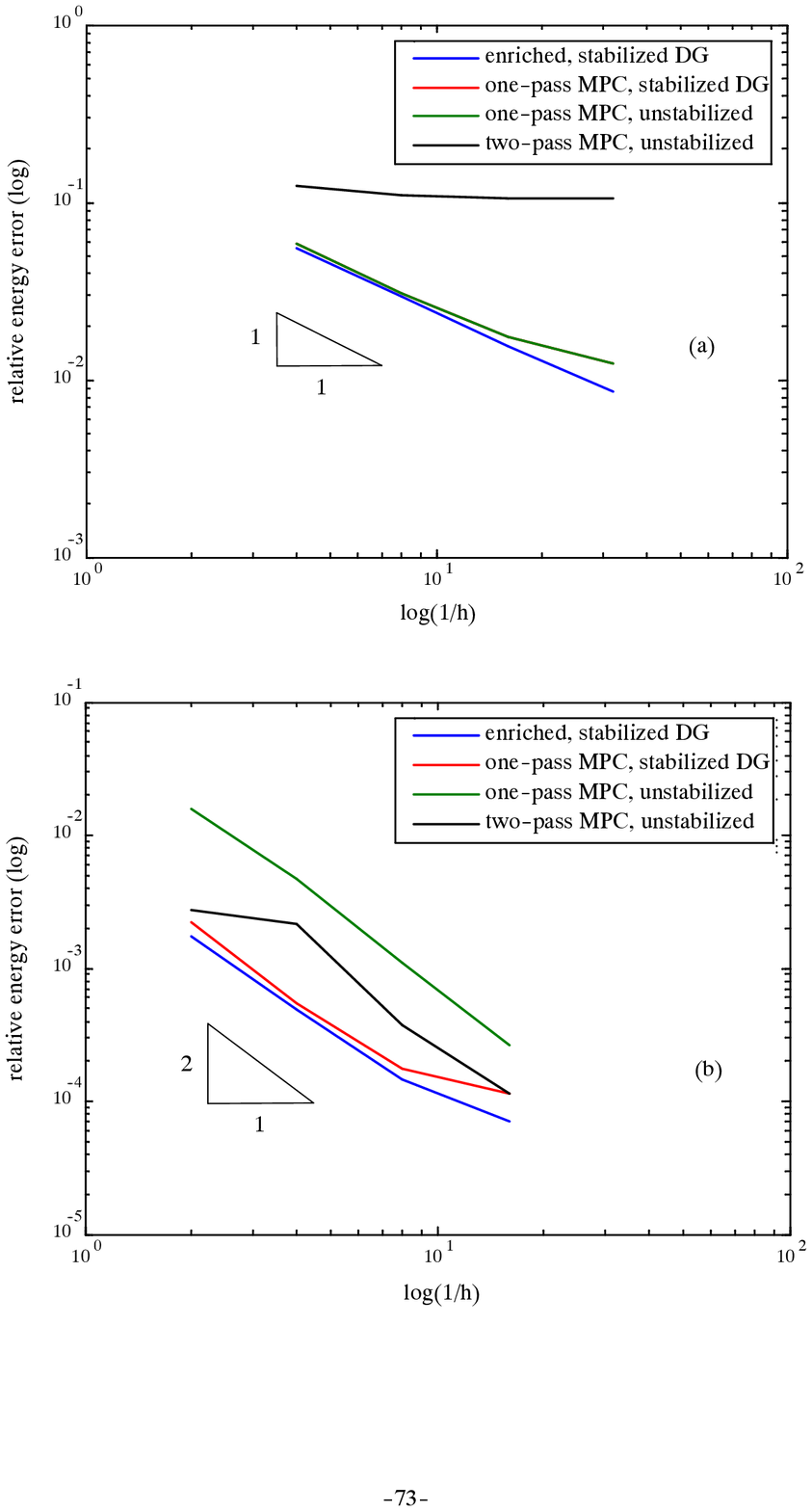}
	\caption{Convergence results using a mesh of (a) Q4 elements and (b) Q8 elements \label{fig:DG_ex_conv}}
	\end{figure}
The deformed configuration is shown in Figure \ref{fig:DG_ex_mesh} (c) and matches the exact solution.

Figures \ref{fig:DG_ex_conv} (a) and (b) shows the spatial convergence in relative energy norm of the proposed formulation using bilinear (Q4)
and quadratic (Q8) elements, respectively. In these figures, we vary the element size $h$ in the third (right-most) layer only while holding the
size of all other elements to be constant. Thus, the energy norms plotted in Figure \ref{fig:DG_ex_conv} were computed over the part of the
mesh that was under refinement only. We compare the convergence rates obtained using 1) the proposed enriched stabilized formulation 2) the single-pass MPC combined with DG stabilization of the interface 3) the single-pass MPC without stabilization, and 4) the two-pass MPC. In the single-pass MPC, we designate the left-hand (coarser) side of the interface to be the master surface. 

The results for the enriched formulation with interface stabilization show a monotonic decrease in
the relative energy norm with convergences rates close to the optimal theoretical values for the underlying mesh.
For very small value of $h$, the error from the unrefined portion of the mesh becomes dominant,
leading to a slightly less-than-optimal convergence rate. 
The performance of the enriched formulation
is clearly superior to that of the MPC approach. Without interface stabilization, the single-pass MPC solution error is 10 times higher than that of the proposed formulation for the case of quadratic elements. Adopting a two-pass approach reduces the energy error for quadratic elements, but the convergence rate in Q4 elements deteriorates. This result is due to the locking caused by enforcing multiple discrete MPCs. The locking is less obvious for Q8 elements, however the system matrices become ill conditioned. 

Applying the DG stabilization in conjunction with the single-pass MPC improves the convergence of the Q8 mesh. This result does not apply to the Q4 mesh, however, and the same result is obtained with and without the DG terms on the interface. This is due to the fact that refining one side of the interface only leads to a Q4 mesh where most elements are connected to the interface in the configuration discussed in Case I of Section \ref{sec_DG:Q4caseI}. In this case enforcing the MPCs leads to point-wise continuity along the interface and the DG terms go to zero, as expected.

\begin{figure}
	\centering
	\includegraphics[clip]{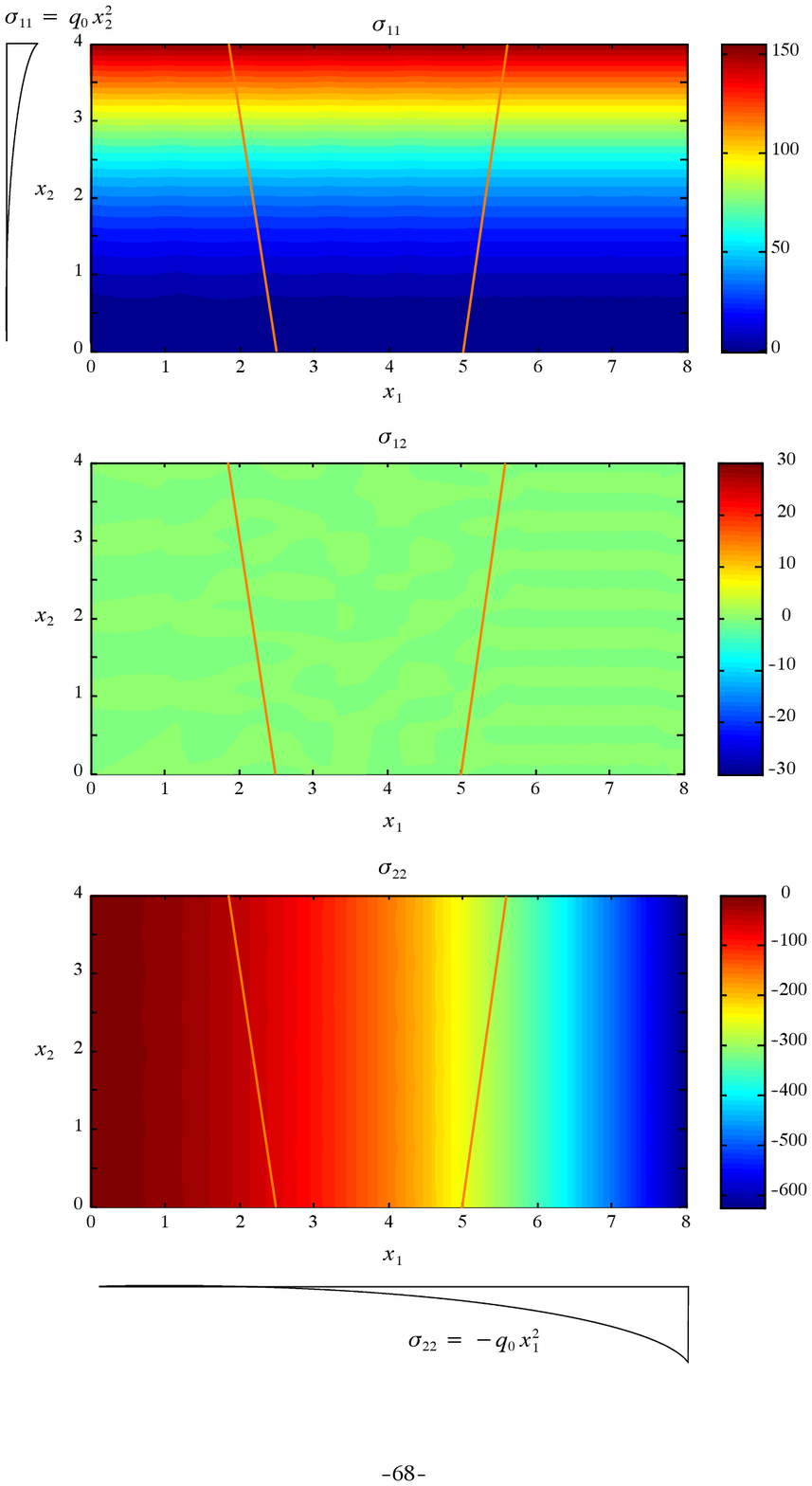}
	\caption{Stress fields \label{fig:DG_ex_stress}}
	\end{figure}

Figure \ref{fig:DG_ex_stress} shows the stress distribution in the domain, as obtained using the proposed formulation. The results match very closely with the exact solution and no perturbations are observed at the interfaces. 
	\begin{figure}
	\centering
	\includegraphics[clip]{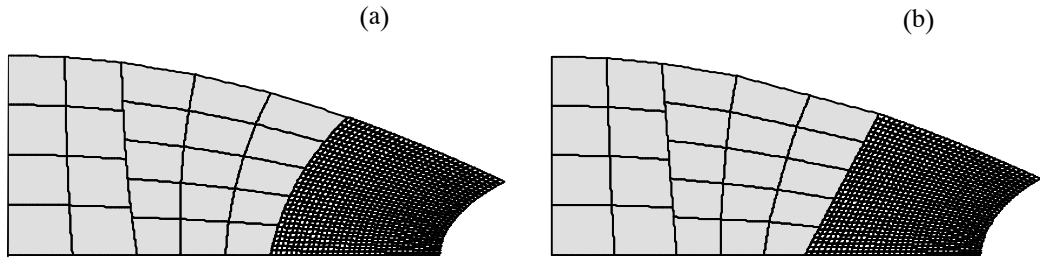}
	\caption{The deformed shape obtained using (a) the proposed interface formulation and (b) the two-pass MPC method showing severe locking  \label{fig:DG_ex_lock_def}}
	\end{figure}
	\begin{figure}
	\centering
	\includegraphics[clip]{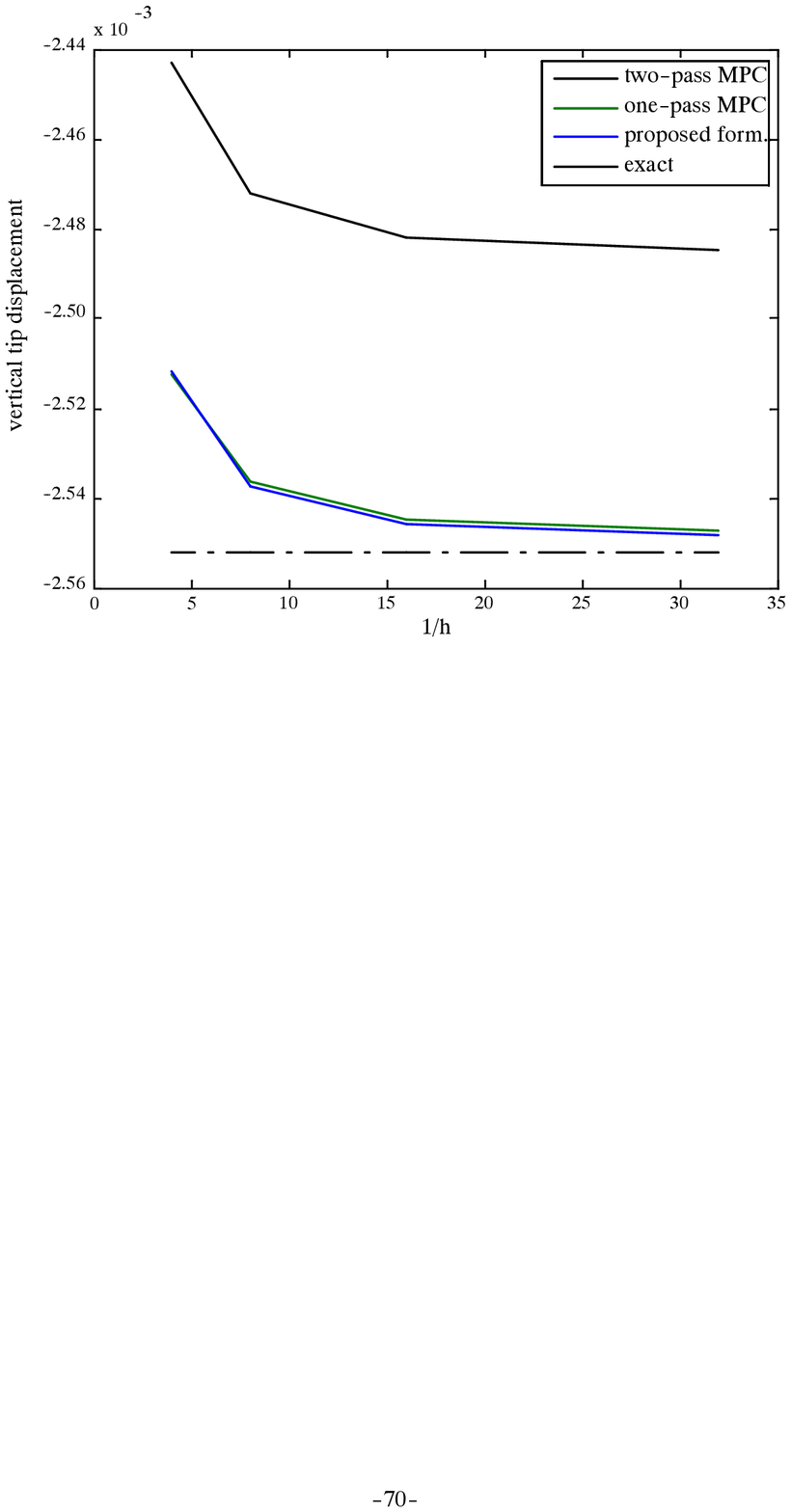}
	\caption{Convergence of the vertical tip displacement with mesh refinement \label{fig:DG_ex_lock_u}}
	\end{figure}

Figure \ref{fig:DG_ex_lock_def} compares the deformed configuration obtained using the proposed method with Q4 elements with the result of the two-pass MPC approach. This figure shows that, despite the large number of degrees of freedom, the two-pass approach leads to a very stiff interface that cannot accommodate the deformation. Figure \ref{fig:DG_ex_lock_u} shows the convergence of the vertical tip displacement with increased mesh refinement for the two-pass MPC, one-pass MPC, and the proposed formulation. It is clear from this figure that, while the single-pass MPC and the proposed formulation converge quickly, the two-pass MPC does not, and shows severe locking.

\begin{figure}
	\centering
	\includegraphics[clip]{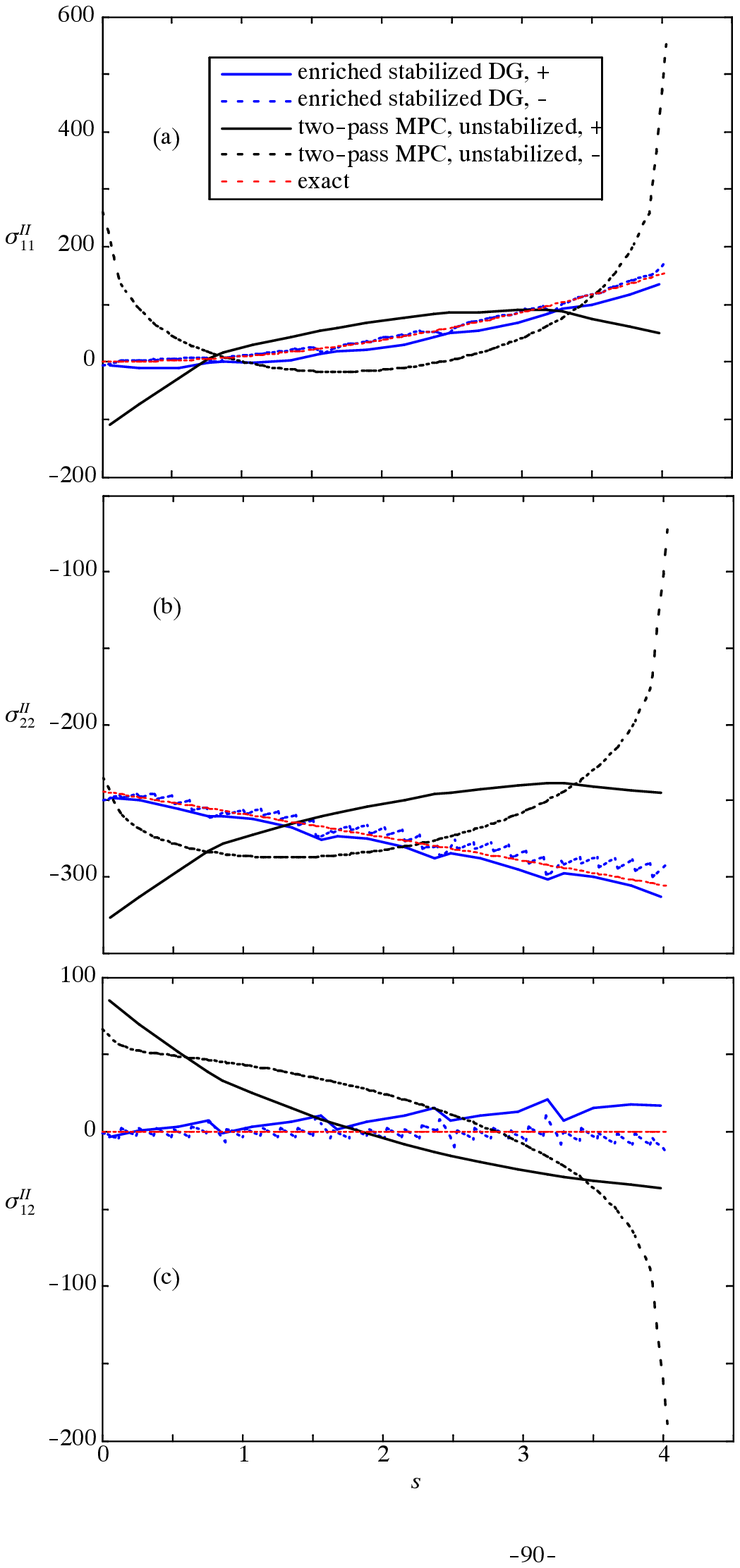}
	\caption{Interface II stresses using Q4 elements (a) $\sigma_{11}$, (b) $\sigma_{22}$, and (c) $\sigma_{12}$ \label{fig:DG_ex_traction_Q4II}}
	\end{figure}

The locking effects can be observed more clearly by comparing the stresses along interface II, as shown in figure \ref{fig:DG_ex_traction_Q4II}. In this figure, we show the variation of the three stress components with $s$, where $s$ is a length variable that parameterizes the interface. We plot the stresses on surfaces $II^-$, and $II^+$ with dotted and solid lines, respectively (the designation of the superscript is based on the direction of the surface normal). It is clear from this figure that, while the stresses obtained using the enriched stabilized formulation are close to the exact solution, the results of the two-pass MPC show significant errors of several orders of magnitude. The error is particularly large at the top of the structure (where deformation is largest) due to the artificial constraint caused by the MPCs.

\begin{figure}
	\centering
	\includegraphics[clip]{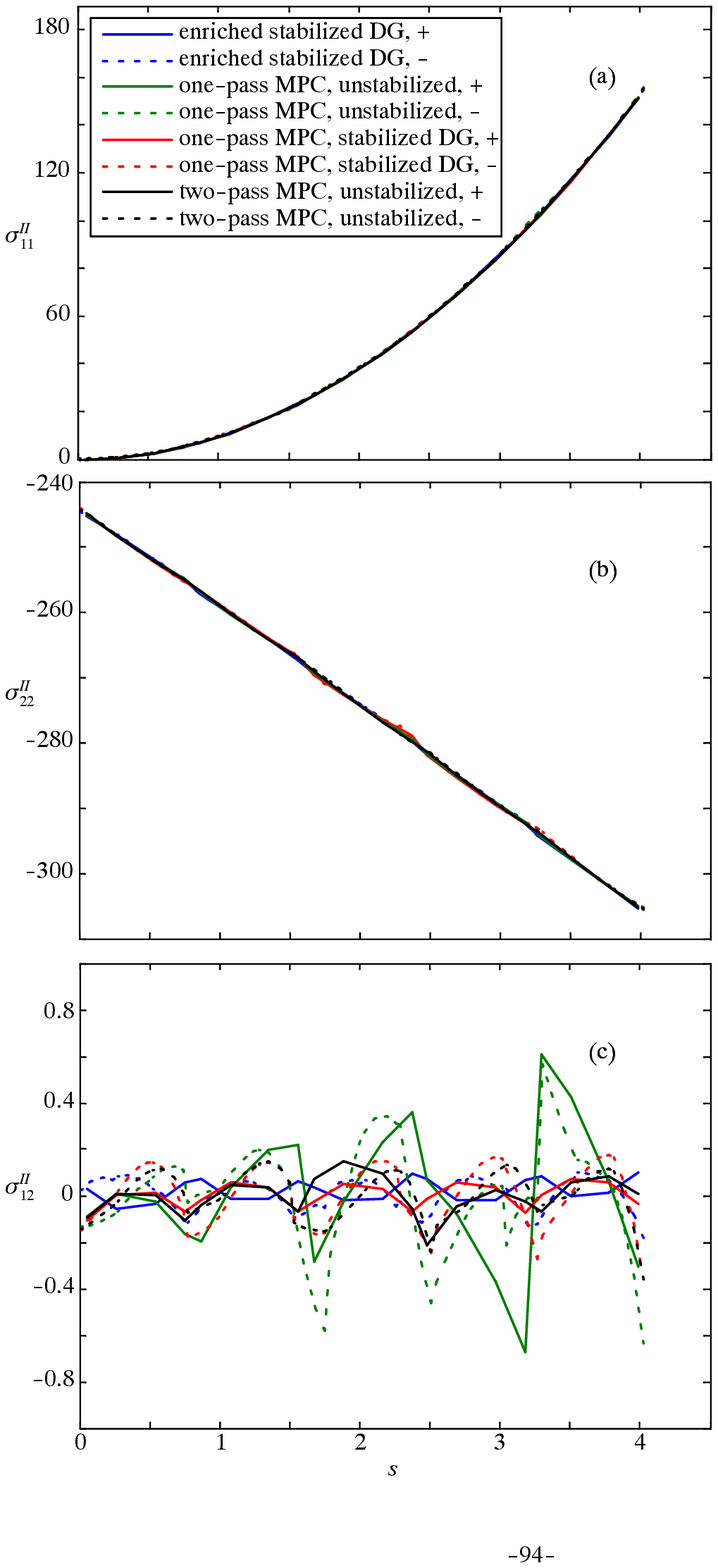}
	\caption{Interface II stresses using Q8 elements (a) $\sigma_{11}$, (b) $\sigma_{22}$, and (c) $\sigma_{12}$ \label{fig:DG_ex_traction_Q8II}}
	\end{figure}

\begin{figure}
	\centering
	\includegraphics[clip]{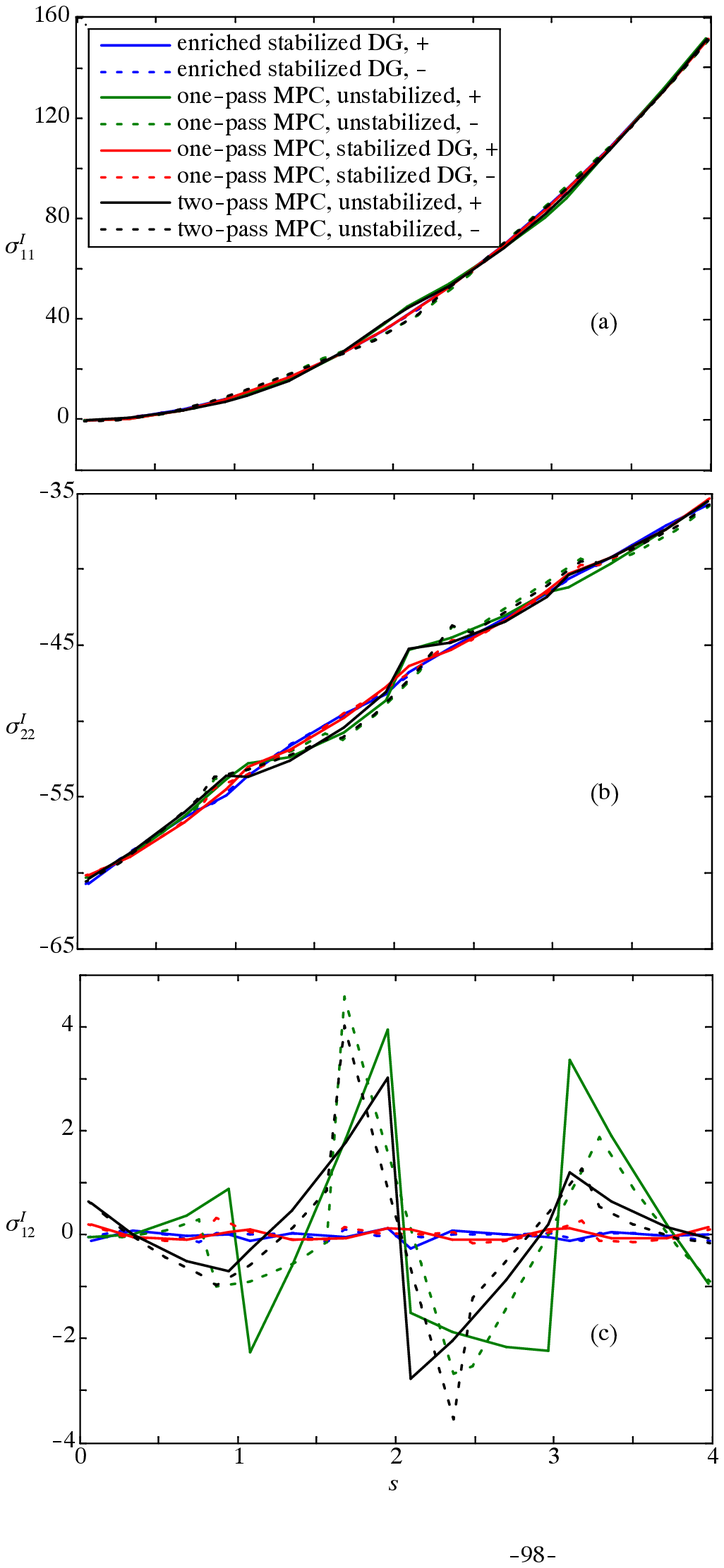}
	\caption{Interface I stresses using Q8 elements (a) $\sigma_{11}$, (b) $\sigma_{22}$, and (c) $\sigma_{12}$ \label{fig:DG_ex_traction_Q8I}}
	\end{figure}

Figures \ref{fig:DG_ex_traction_Q8II} and \ref{fig:DG_ex_traction_Q8I} compare the stress distribution obtained using 1) the proposed enriched stabilized formulation 2) the single-pass MPC combined with DG stabilization of the interface 3) the single-pass MPC without stabilization, and 4) the two-pass MPC, for Q8 elements. Even though the variations among these methods are minimal for normal stresses, the shear stresses obtained using the unstabilized one-pass 3) and two-pass MPCS 4) methods display severe oscillations. These oscillations are greatly reduced with the DG stabilization procedure 2), with further improvement due to enrichment 1). The proposed DG-based stabilization is very successful in projecting the stresses from one side of the interface to the other.

\begin{figure}
	\centering
	\includegraphics[clip]{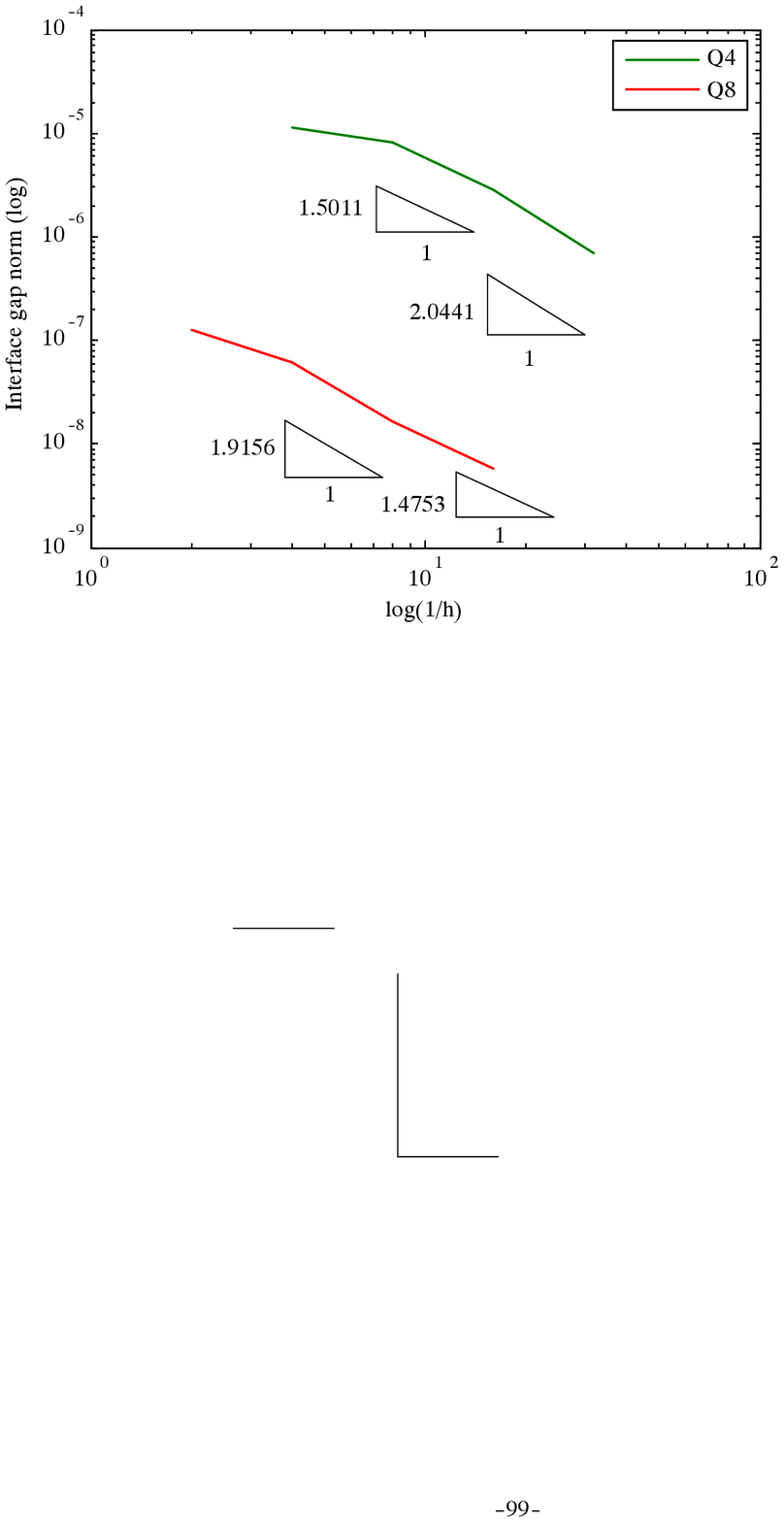}
	\caption{Convergence of displacement gap along the interface \label{fig:DG_ex_gap}}
	\end{figure}

Figure \ref{fig:DG_ex_gap} shows the convergence of the interface gap in $L_2$ norm with mesh refinement using the proposed method. We observe good convergence rates in Q4 and Q8 elements.
%
\section{Conclusions} \label{sec_DG:conclusion}

We have proposed a new primal approach for the coupling of non-conforming finite element meshes. This approach is based on a local enrichment of the non-conforming interface that enables a simple enforcement of the continuity of the displacement field using a set of discrete node-to-node constraints. These constraints can be taken care of automatically by the assembly procedure without the need for additional
variables or Lagrange multipliers. The enrichment enables a two-pass approach that eliminates the need for a master/slave treatment and enables an unbiased enforcement of the continuity condition at all nodes along the interface. This aspect of the formulation is not possible in currently available coupling formulations where the continuity condition is only enforced in a weak sense. The advantage of enabling a strong continuity enforcement at all interface nodes is in the applicability of the method to the treatment of non-smooth contact, as will be discussed in the following chapter. The primal nature of the formulation releases the restrictions of the LBB condition for dual coupling methods. 
We treat the interface using a form of the discontinuous Galerkin
method that guarantees the complete transfer of forces along non-conforming inter-element
boundaries. The proposed interface formulation is consistent and includes the continuous Galerkin as a subset. Results show that the proposed formulation is stable, exhibits good convergence properties, and yields greatly improved estimates of the interface stress fields.











%% file: enrichment.tex

\chapter{A stabilized interface formulation for contact problems with sliding} \label{enrichment}
%
%
%
%
%

\section{Introduction} 

The formulation of contact problems is similar in many ways to the coupling of non-conforming meshes. The modeling task for both involves ensuring geometric compatibility and complete transfer of forces on the interface. One feature that distinguishes the contact problem is that the bodies are allowed to slide tangentially along the interface. Therefore, the continuity condition applies to the normal component of the displacement field only. Another feature is the unilateral nature of the coupling, which transforms the mathematical framework to an inequality constrained optimization, as described in Chapter 2. 

Many of the coupling formulations described in section \ref{sec_DG:intro} have been applied to the contact problem, with a heavy emphasis on dual methods. The first dual contact formulation can be traced to Papadopoulos and Taylor \cite{papadopoulosCMAME92} who introduced the concept of applying a field of Lagrange multipliers representing the contact pressure to enforce a weak interpenetrability constraint. A later development of this work \cite{jonesIJNME01} proposed enforcing equilibrium weakly between two isoparametric interpolations of the contact pressure on either side of the interface, while enforcing the interpenetration constraint strongly at the integration points. Both the mortar \cite{pusoCMAME04,yangIJNME05,pusoCMAME08,brunbenIUTAM07} and, more recently, FETI \cite{avery:ffd} methods have been used to model contact conditions. Of the primal approaches, the penalty-based Nitsche method can also be found in the literature on contact \cite{wriggers2006}.

As discussed in Chapter 4, the disadvantage of dual methods lies in the Ladyzhenskaya-Babu\v{s}ka-Brezzi (LBB) restriction and the  bias due to the choice of the master surface. Unbiased dual two-pass approaches fail the LBB condition and are therefore prone to surface locking \cite{solbergCMAME05}. This is particularly true of the discrete node-to-surface gap function formulation for contact, similar to the MPC approach for coupling non-conforming meshes. The locking issue is more prevalent in contact problems due to the strict enforcement of the no-penetration conditions at locations where contact is detected. Unlike the coupling of nonconforming meshes, the contact interface  cannot be generally determined a-priori, and does not remain the same through deformation. Therefore, enforcing a set of discrete constraints at locations where contact is detected makes sense from a physical point of view. Unfortunately, within the constraints of the numerical problem, such an approach is unstable. 

A mortar-based two-pass algorithm has recently been developed by Solberg et al. \cite{solbergCMAME07}. Unlike the single-pass mortar method, this approach strongly enforces the contact constraints at a number of select points while continuity of pressure is satisfied in a weak sense along the contact surface. As a result, a penalty-based stabilization term needs to be imposed to minimize the pressure jump across the contact interface. It is suggested that \textit{on each surface, nodes [be] a priori identified as active or inactive, via a binary patterning scheme}. This ad-hoc approach is not guaranteed to work for arbitrary meshes.

In this chapter, we extend the stabilized interface formulation discussed in Chapter 4 to the contact problem. The main focus of the chapter is to adapt the enrichment procedure to the nature of the contact problem where the node is allowed to slide on surface of the contact element. This feature is desirable to maintain the node-to-node form of the contact constraints and is most conveniently developed in a quasi-static framework. As for the coupling of non-conforming meshes, we drop the inertia terms to simplify the development of the method. We show that the added degrees of freedom enable an unbiased two-pass approach for the enforcement of the interpenetrability constraint at all interface nodes without inducing surface locking. 

The outline of the chapter is as follows. In Section \ref{sec_enr:EnrProcedure}, we describe the enrichment procedure for the case of a node sliding along the contact surface, and its effect on the formulation of the enriched element and contact constraints. We then discuss the stabilized interface formulation for contact in Section \ref{sec_enr:StabFormulation}. A few numerical examples are presented and discussed in Section \ref{sec_enr:Results}. Section \ref{sec_enr:Concl} provides conclusions.
%
%
%
%
%
\section{Moving enrichment procedure} \label{sec_enr:EnrProcedure}
	\begin{figure}
	\centering
	\includegraphics[clip]{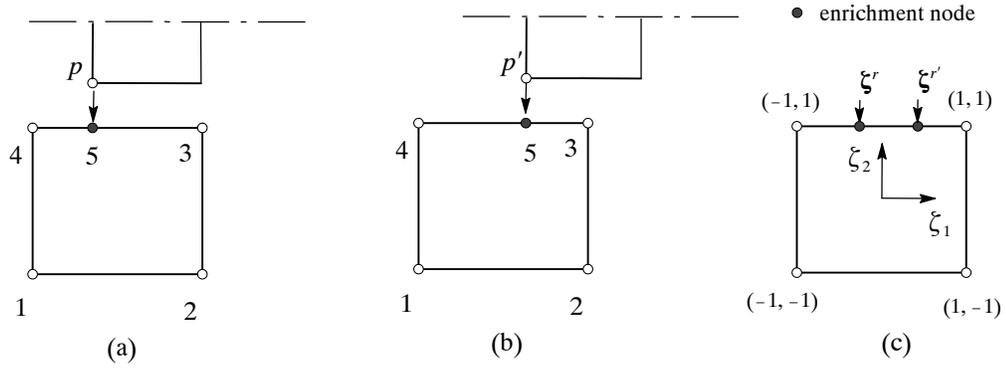}
	\caption{(a) Local enrichment of the interface element (b) enrichment update for sliding, and (c) moving node reference in the parent domain \label{fig:enr_Enr}}
	\end{figure}
As is the case for the coupling problem, the aim of the enrichment is to transform the node-to-surface contact constraint to a node-to-node one without modifying the existing mesh. To achieve this objective, we adopt a similar approach: If contact is detected between a point $p$ and an element surface, a node can be inserted at a location $\mathbf\zeta^{r}$ such that the gap function is reduced to a node-to-node constraint, as shown in Figure \ref{fig:enr_Enr}(a). Completeness of the finite element interpolation in the contact element can be preserved by updating the set of Lagrangian shape functions to account for the additional node. For example, the updated shape functions corresponding to the case illustrated in Figure \ref{fig:enr_Enr}(a) are
		\begin{align} 
		\tilde{N}^{5}(\mathbf{\zeta},\mathbf{\zeta}^{r}) &= \textstyle\frac{1}{2}(\zeta_{1}+1)
		\displaystyle\frac{(\zeta_{2}+1)(\zeta_{2}-1)}{(\zeta_{2}^{r}+1)(\zeta_{2}^{r}-1)} 
		\\ 
		\tilde{N}^{\alpha} &= N_{Q4}^{\alpha}-N_{Q4}^{\alpha}(\mathbf{\zeta}^{r})\tilde{N}^{5}, 
		\end{align}
where $N_{Q4}^{\alpha}$ are the shape functions of a Q4 element. For an element $m$ defined by a set of nodes with spatial coordinates $\mathbf{x}^{\alpha}$ and corresponding shape functions $N^{\alpha}$, where $\alpha = 1,...,n$.

The enriched spatial $\mathbf{x}$, material $\mathbf{X}$, and displacement $\mathbf{u}$ fields in the element are given by (summation from 1 to $n$ is implied on $\alpha$)
		\begin{align} \label{eq_enr:xEnrIso}
		\mathbf{x}(\mathbf{\zeta},\mathbf{\zeta}^{r}) = \tilde{N}^{\alpha}(\mathbf{\zeta},\mathbf{\zeta}^{r}) \mathbf{x}^{\alpha} 
		+ \tilde{N}^{r}(\mathbf{\zeta},\mathbf{\zeta}^{r})\mathbf{x}^{r} , \ \ \  \alpha = 1,\cdots,n \\
		\label{eq_enr:XEnrIso}
		\mathbf{X}(\mathbf{\zeta},\mathbf{\zeta}^{r}) = \tilde{N}^{\alpha}(\mathbf{\zeta},\mathbf{\zeta}^{r}) \mathbf{X}^{\alpha} 
		+ \tilde{N}^{r}(\mathbf{\zeta},\mathbf{\zeta}^{r})\mathbf{X}^{r} , \ \ \  \alpha = 1,\cdots,n \\
		\label{eq_enr:uEnrIso}
		\mathbf{u}(\mathbf{\zeta},\mathbf{\zeta}^{r}) = \tilde{N}^{\alpha}(\mathbf{\zeta},\mathbf{\zeta}^{r}) \mathbf{d}^{\alpha} 
		+ \tilde{N}^{r}(\mathbf{\zeta},\mathbf{\zeta}^{r})\mathbf{d}^{r} , \ \ \  \alpha = 1,\cdots,n,
		\end{align}
where $\tilde{N}^{r}$ is the shape function associated with the additional node and the $\tilde{N}^{\alpha}$ are the modified (enriched) element shape functions defined as follows
		\begin{equation} \label{eq_enr:NEnr}
		\tilde{N}^{\alpha}(\mathbf{\zeta},\mathbf{\zeta}^{r}) = 
		N^{\alpha}(\mathbf{\zeta})-N^{\alpha}(\mathbf\zeta^{r})\tilde{N}^{r}(\mathbf{\zeta},\mathbf{\zeta}^{r}),\ \ \ \alpha = 1,\cdots,n.
		\end{equation}
	\begin{figure}
	\centering
	\includegraphics[clip]{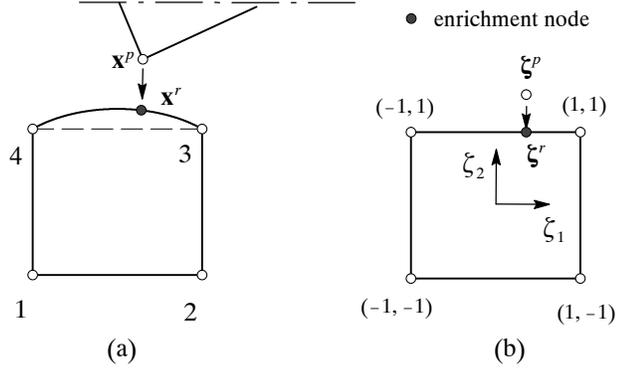}
	\caption{(a) Spatial coordinates of the contact $\mathbf{x}^p$ and enrichment $\mathbf{x}^r$ nodes (c) mapping to the parent domain \label{fig:enr_xpxr}}
	\end{figure}
Note that, in the above equations, we have used the notation $\mathbf{\zeta}^r$ and $\mathbf{x}^r$ to distinguish the enrichment node $r$ from the contact node $p$. This distinction, while irrelevant to the coupling problem, is essential to the contact problem to allow $p$ to move away from the element if the contact constraint gets deactivated, as shown in Figure \ref{fig:enr_xpxr} (a). Although in the limit of contact resolution $\mathbf{x}^p = \mathbf{x}^r$, they are essentially independent quantities. This distinction becomes clearer in the calculation of the enrichment location $\mathbf{\zeta}^r$, as shown in Figure \ref{fig:enr_xpxr} (b). As for the coupling problem, we map the contact node $p$ to the parent domain of the contact element by solving the system of nonlinear equations 
\begin{equation}  \label{eq_enr:enr_loc}
\mathbf{x}^p = N^\alpha\left(\mathbf{\zeta}^p\right)\mathbf{x}^\alpha
\end{equation}
for $\mathbf{\zeta}^p$. If contact has occurred through the surface $\zeta_{j}=c$, the enrichment location $\mathbf{\zeta}^r$ can then be obtained from $\mathbf{\zeta}^p$ by setting $\zeta_{j}^p=c$. As a result, even though the node $p$ can separate from the contact element, the enrichment node $r$ remains glued to the surface. An initial estimate of the spatial coordinates of the added node $r$ can be computed as
\begin{equation}  \label{eq_enr:xenr}
\mathbf{x}^r = N^\alpha\left(\mathbf{\zeta}^r\right)\mathbf{x}^\alpha.
\end{equation}

As in the coupling problem (Section \ref{sec_DG:enr_procedure}), we assume that the mapping between the material domain and the element parent (reference) domain is not affected by the enrichment on the surface, to wit
		\begin{align} 
		\notag
		\mathbf{X}(\mathbf{\zeta}) &= \tilde{N}^{\alpha}(\mathbf{\zeta},\mathbf{\zeta}^{r}) \mathbf{X}^{\alpha} 
		+ \tilde{N}^{r}(\mathbf{\zeta},\mathbf{\zeta}^{r})\mathbf{X}^{r}  \\ \label{eq_enr:XEnrIso2}
		&=  N^{\alpha} \mathbf{X}^{\alpha}  \ \ \  \alpha = 1,\cdots,n 
		\end{align}
Therefore, the material point corresponding to the added node $\mathbf{X}^{r}$ satisfies the constraint
		\begin{equation} \label{eq_enr:Xp}
		\mathbf{X}^{r} = N^{\alpha}(\mathbf\zeta^{r})\mathbf{X}^{\alpha},
		\end{equation}
Unlike the coupling problem, however, the enrichment location, dictated by the spatial variable $\mathbf{x}^p$, does not remain fixed. As the point $p$ slides along the element surface, the enrichment reference $\mathbf{\zeta}^{r}$ has to be updated to accommodate the motion such that the contact remains node-to-node. Consequently, the enrichment reference needs to be computed in the spatial domain. Moreover, the material point associated with the enrichment $\mathbf{X}^{r}$ does not remain fixed through deformation and the finite element discretization in the contact element now includes a node of moving reference. 

The presence of a moving reference implies that the element formulation is no longer purely Lagrangian. Such formulation is typically addressed via the Arbitrary Lagrangian Eulerian framework where the convective effects due to the change in the mapping between the reference and material domains are accounted for. However, given the assumption we made in equation \eqref{eq_enr:Xp}, the enrichment does not affect the mapping between the boundary of the material domain and that of the enriched element. Therefore, no mass convection occurs through the mesh boundary. In the following section we show that if the enrichment location and corresponding material, spatial and displacement fields are updated properly, the presence of a moving reference does not affect the Lagrangian nature of the element formulation. 
%
%
\subsection{Enriched element formulation with a moving reference}
In this section, we show that the enriched element formulation retains the essential features of the underlying Lagrangian mesh regardless of the presence of a node with moving reference. To prove this proposition, we show that the convective terms that result from the moving enrichment reference come out to be exactly zero. 

First, let us rewrite the enriched spatial field in the element as follows
		\begin{align} \label{eq_enr:xEnr2}
		\mathbf{x} &= \tilde{N}^{\alpha}(\mathbf{\zeta},\mathbf{\zeta}^{r}) \, \mathbf{x}^\alpha + \tilde{N}^r (\mathbf{\zeta},\mathbf{\zeta}^{r}) \,\mathbf{x}^r \\
		&= \left[
		N^{\alpha}(\mathbf{\zeta})-N^{\alpha}(\mathbf\zeta^{r})\tilde{N}^{r}(\mathbf{\zeta},\mathbf{\zeta}^{r})\right] \mathbf{x}^\alpha + \tilde{N}^r (\mathbf{\zeta},\mathbf{\zeta}^{r}) \,\mathbf{x}^r \\
		&= N^\alpha (\mathbf{\zeta}) \, \mathbf{x}^\alpha + \tilde {N}^r (\mathbf{\zeta},\mathbf{\zeta}^{r}) \left [ \mathbf{x}^r -N^\alpha (\mathbf{\zeta}^r) \,\mathbf{x}^\alpha \right ] ,\ \ \ \alpha = 1,\cdots,n.
		\end{align}
To compute the rate of change of the spatial field with deformation, we define a pseudo-time variable $t$ and compute the total derivative of $\mathbf{x}$ with respect to $t$ as follows. 
%
%
		\begin{align} 
		\frac{d\mathbf{x}}{dt} &= \frac{d}{dt} \left\{ N^\alpha (\mathbf{\zeta}) \, \mathbf{x}^\alpha \,(t) + \tilde {N}^r (\mathbf{\zeta},\mathbf{\zeta}^{r} \,(t)) \left [ \mathbf{x}^r\,(t) -N^\alpha (\mathbf{\zeta}^r \,(t)) \,\mathbf{x}^\alpha \,(t) \right ]  \right\} \\ \notag
		 &= N^{\alpha}(\mathbf{\zeta}) \frac{d\mathbf{x}^{\alpha}}{dt}
		+\frac{d}{dt} \tilde{N}^{r} (\mathbf{\zeta},\mathbf{\zeta}^{r}) 
		\left[\mathbf{x}^{r}-N^{\alpha}(\mathbf{\zeta}^{r})\mathbf{x}^{\alpha}\right] \\ \label{eq_enr:dxdt}
		&+ \tilde{N}^{r}\left[\frac{d\mathbf{x}^{r}}{dt}-N^{\alpha}(\mathbf{\zeta}^{r})\frac{d\mathbf{x}^{\alpha}}{dt}
		-\frac {d}{dt}N^{\alpha}(\mathbf{\zeta}^{r}) \mathbf{x}^{\alpha} \right].
		\end{align}
Equation \eqref{eq_enr:dxdt} can be simplified as follows. First, note that the terms $\tilde{N}^r(\mathbf{\zeta},\mathbf{\zeta}^{r})$ and  $N^\alpha(\mathbf{\zeta}^r)$ are dependent on $t$ trough $\mathbf{\zeta}^{r}$. The rate of change of these quantities with respect to $t$ is therefore computed as
	  \begin{equation} \label{eq_enr:ALE1}
	  \frac{d}{dt}\tilde{N}^{r}(\mathbf{\zeta},\mathbf{\zeta}^{r})= \frac{\partial}{\partial\zeta^{r}}\tilde{N}^{r}
	  (\mathbf{\zeta},\mathbf{\zeta}^{r}) \cdot\frac{d\mathbf{\zeta}^{r}}{dt}   \\\ , \\\
	  \frac{d}{dt}N^{\alpha}(\mathbf{\zeta}^{r})=\frac{\partial{N}^{\alpha}(\mathbf{\zeta}^{r})}{\partial\mathbf{\zeta}^{r}}\cdot 
	  \frac{d\mathbf{\zeta}^{r}}{dt}.
	  \end {equation}
These convective terms account for the change in the reference $\mathbf{\zeta}^{r}$, and therefore the corresponding material point $\mathbf{X}^{r}$, as a result of deformation. Furthermore, ${d\mathbf{x}^{r}}/{dt}$ represents the total pseudo-time derivative of the enrichment spatial variable and is computed as follows
	  \begin{align}  \notag
	  \frac{d\mathbf{x}^{r}}{dt} &= \frac{\partial\mathbf{x}^{r}}{\partial t}| _{\mathbf{\zeta}^{r}}
	  + \frac{\partial\mathbf{x}^{r}}{\partial \mathbf{\zeta}^{r}} \cdot \frac{d\mathbf{\zeta}^{r}}{dt} \\
	  \label{eq_enr:ALE2}
	  &= \dot{\mathbf{x}}^{r}+ \frac{\partial\mathbf{x}^{r}}{\partial \mathbf{\zeta}^{r}} \cdot \frac{d\mathbf{\zeta}^{r}}{dt}.
	  \end{align}
In this equation $\dot{\mathbf{x}}^r$ is the local rate of change of the spatial variable $\mathbf{x}^r$ assuming $\mathbf{\zeta}^r$ is fixed. The second term in equation \eqref{eq_enr:ALE2} is a convective term that accounts for the change in the value of the enrichment spatial coordinate vector $\mathbf{x}^r$ when the location at which that vector is computed changes while all other kinematic variables remain unchanged (see Figure \ref{fig:enr_dzetar}) (a).
	  \begin{align} \notag
	  \frac{\partial\mathbf{x}^{r}}{\partial\mathbf{\zeta}^{r}} \cdot \frac{d\mathbf{\zeta}^{r}}{dt} &= 
	  \frac {d}{d\epsilon} \mathbf{x}(\mathbf{\zeta}^{r}+\epsilon \frac{d\mathbf{\zeta}^{r}}{dt},\mathbf{\zeta}^{r})|_{\epsilon = 0} \\
	  \notag
	  &= \frac{d}{d\epsilon}\left[ N^{\alpha}(\mathbf{\zeta}^{r}+ \epsilon \frac{d\mathbf{\zeta}^{r}}{dt}) \mathbf{x}^{\alpha} 
	  + \tilde{N}^{r}(\mathbf{\zeta}^{r} + \epsilon \frac{d\mathbf{\zeta}^{r}}{dt}, \mathbf{\zeta}^{r}) \left[\mathbf{x}^{r}
	  -N^{\alpha}(\mathbf{\zeta}^{r})\mathbf{x}^{\alpha}\right]\right]_{\epsilon = 0}       \\
	  \notag
	  &= \frac{\partial}{\partial\mathbf{\zeta}^{r}}{N}^{\alpha}(\mathbf{\zeta}^{r})
	  \cdot\frac{d\mathbf{\zeta}^{r}}{dt}\mathbf{x}^{\alpha}
	  + [\frac{\partial}{\partial\mathbf{\zeta}}\tilde{N}^{r}
	  (\mathbf{\zeta},\mathbf{\zeta}^{r})_{\mathbf{\zeta}^{r}}\cdot\frac{d\mathbf{\zeta}^{r}}{dt} ]
	  \left[\mathbf{x}^{r}-N^{\alpha}(\mathbf{\zeta}^{r})\mathbf{x}^{\alpha}\right]   \\
	  \label{eq_enr:ALE3}
	  &= \frac{d}{dt} {N}^{\alpha}(\mathbf{\zeta}^{r}) \mathbf{x}^{\alpha}+ [\frac{\partial}{\partial\mathbf{\zeta}}\tilde{N}^{r}
	  (\mathbf{\zeta},\mathbf{\zeta}^{r})_{\mathbf{\zeta}^{r}}\cdot\frac{d\mathbf{\zeta}^{r}}{dt} ]
	  \left[\mathbf{x}^{r}-N^{\alpha}(\mathbf{\zeta}^{r})\mathbf{x}^{\alpha}\right]
	  \end{align}
where $\frac{\partial}{\partial\mathbf{\zeta}}\tilde{N}^{r}(\mathbf{\zeta},\mathbf{\zeta}^{r})_{\mathbf{\zeta}^{r}}$ should be read as the partial derivative of $\tilde{N}^{r}(\mathbf{\zeta},\mathbf{\zeta}^{r})$ with respect to $\mathbf{\zeta}$ evaluated at $\mathbf{\zeta}^{r}$.
Substituting Equations \eqref{eq_enr:ALE1}, \eqref{eq_enr:ALE2} and \eqref{eq_enr:ALE3} into \eqref{eq_enr:dxdt} yields
		\begin{align} \notag
		\frac{d\mathbf{x}}{dt} &= \left[N^{\alpha}(\mathbf{\zeta})- N^{\alpha}(\mathbf{\zeta}^{r})\tilde{N}^{r} \right]
		\dot{\mathbf{x}}^{\alpha}+ \tilde{N}^{r} \dot{\mathbf{x}}^{r} \\ \label{eq_enr:dxdt2}
		&+\left[ \frac{\partial}{\partial\zeta^{r}}\tilde{N}^{r}
	  (\mathbf{\zeta},\mathbf{\zeta}^{r})+\tilde{N}^{r}\frac{\partial}{\partial\mathbf{\zeta}}\tilde{N}^{r}
	  (\mathbf{\zeta},\mathbf{\zeta}^{r})_{\mathbf{\zeta}^{r}}\right]\cdot 
	  \frac{d\mathbf{\zeta}^{r}}{dt}\cdot\left[N^{\alpha}(\mathbf{\zeta})- N^{\alpha}(\mathbf{\zeta}^{r})\tilde{N}^{r} \right].
		\end{align}
Note that, at the Lagrangian nodes $\mathbf{x}^\alpha$, the local rate of change $\dot{\mathbf{x}^\alpha}$ is identical to the total rate of change $d\mathbf{x}^\alpha/dt$.
The first two terms in equation \eqref{eq_enr:dxdt2} represent the Lagrangian phase of the solution where the effects of mesh motion are not present. These effects are captured by the third term which contains the derivatives of the enrichment function $N^{r}$. Given that the element shape functions after enrichment are interpolatory, we can assume that the enrichment function is of the form
		\begin{equation}
		\tilde{N}^{r}(\mathbf{\zeta},\mathbf{\zeta}^{r}) = f(\mathbf{\zeta}) / f(\mathbf{\zeta}^{r}),
		\end{equation}
where $f(\mathbf{\zeta}) = 0$ at all nodes $\alpha$. Taking the partial derivative with respect to $\mathbf{\zeta}^r$ yields
		\begin{align} \label{eq_enr:part1}
    \frac{\partial}{\partial\zeta^{r}}\tilde{N}^{r} (\mathbf{\zeta},\mathbf{\zeta}^{r})
	  &= - \frac{f(\mathbf{\zeta})} {f^{2}(\mathbf{\zeta}^{r})} \frac{\partial f(\mathbf{\zeta}^{r})}{\partial\zeta^{r}} \\
	  \label{eq_enr:part2}
    \frac{\partial}{\partial\mathbf{\zeta}}\tilde{N}^{r} (\mathbf{\zeta},\mathbf{\zeta}^{r})_{\mathbf{\zeta}^{r}} 
    &= \frac{1}{f(\mathbf{\zeta}^{r})}\left[\frac{\partial f(\mathbf{\zeta})}{\partial\mathbf{\zeta}}\right] _{\mathbf{\zeta}^{r}} 
    = \frac{1}{f(\mathbf{\zeta}^{r})} \frac{\partial f(\mathbf{\zeta}^{r})} {\partial\mathbf{\zeta}^{r}} \\
    \label{eq_enr:part3}
    \frac{\partial}{\partial\zeta^{r}}\tilde{N}^{r}(\mathbf{\zeta},\mathbf{\zeta}^{r})
    +\tilde{N}^{r}\frac{\partial}{\partial\mathbf{\zeta}}\tilde{N}^{r} (\mathbf{\zeta},\mathbf{\zeta}^{r})_{\mathbf{\zeta}^{r}} 
    &= - \frac{f(\mathbf{\zeta})} {f^{2}(\mathbf{\zeta}^{r})} \frac{\partial f(\mathbf{\zeta}^{r})}{\partial\zeta^{r}}
    + \frac {f(\mathbf{\zeta})} { f(\mathbf{\zeta}^{r})} 
    \frac{1}{f(\mathbf{\zeta}^{r})} \frac{\partial f(\mathbf{\zeta}^{r})} {\partial\mathbf{\zeta}^{r}}
    = \mathbf{0}
		\end{align}
The higher-order terms due to the moving reference vanish due to the interpolatory nature of the enrichment. Therefore, the rate of change of the spatial field in the enriched element can be solely described by its Lagrangian (fixed reference) components
		\begin{align} \label{eq_enr:dxdt3}
		\frac{d\mathbf{x}}{dt} &= \left[N^{\alpha}(\mathbf{\zeta})- N^{\alpha}(\mathbf{\zeta}^{r})\tilde{N}^{r} \right]
		\dot{\mathbf{x}}^{\alpha}+ \tilde{N}^{r} \dot{\mathbf{x}}^{r}
		\end{align}
assuming the location of the enrichment, and the corresponding spatial (and therefore material and displacement) variables, are updated simultaneously. The update procedure is shown in Figure \ref{fig:enr_dzetar} (b) and goes as follows:
	\begin{figure}
	\centering
	\includegraphics[clip]{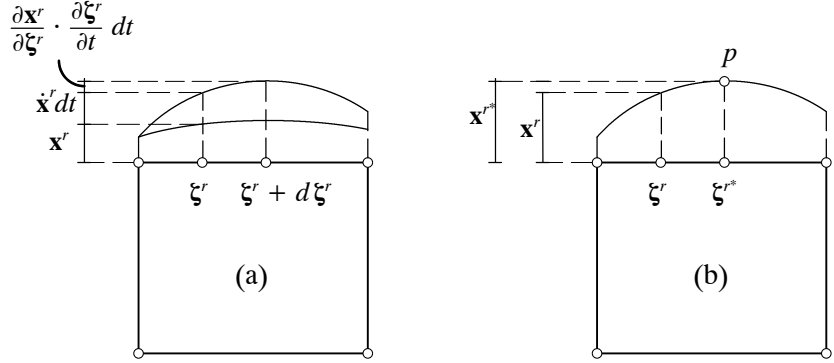}
	\caption{(a) Rate of change of the enrichment spatial variable $\mathbf{x}^r$ (b) Enrichment update procedure \label{fig:enr_dzetar}}
	\end{figure}
\begin{enumerate}\item {Given the previous enrichment reference $\mathbf{\zeta}^r$, spatial element coordinates $\mathbf{x}^r$,$\mathbf{x}^\alpha$ for $\alpha = 1,...,n$, and spatial location of the contact node $\mathbf{x}^p$, find $\mathbf{\zeta}^p$ such that
\begin {equation}
\mathbf{x}^p = \tilde{N}^\alpha (\mathbf{\zeta}^p,\mathbf{\zeta}^r)\, \mathbf{x}^\alpha + \tilde{N}^r (\mathbf{\zeta}^p,\mathbf{\zeta}^r) \mathbf{x}^r,
\end{equation}}
\item Assuming the contact surface $\zeta_j=c$, find the new enrichment location $\mathbf{\zeta}^{r*} = \mathbf{\zeta}^p$ with $\zeta_j^{r*} = c$
\item {Compute the new spatial and material coordinates 
\begin{align}
\mathbf{x}^{r*} &= \tilde{N}^\alpha (\mathbf{\zeta}^{r*},\mathbf{\zeta}^r)\, \mathbf{x}^\alpha + \tilde{N}^r (\mathbf{\zeta}^{r*},\mathbf{\zeta}^r) \mathbf{x}^r,\\ \notag
\mathbf{X}^{r*} &= \tilde{N}^\alpha (\mathbf{\zeta}^{r*},\mathbf{\zeta}^r)\, \mathbf{X}^\alpha + \tilde{N}^r (\mathbf{\zeta}^{r*},\mathbf{\zeta}^r) \mathbf{X}^r,\\
&= N^\alpha (\mathbf{\zeta}^{r*}) \, \mathbf{x}^\alpha
\end{align}}
\item Set $\mathbf{x}^r = \mathbf{x}^{r*}$ and $\mathbf{X}^r = \mathbf{X}^{r*}$
\item Update the displacement $\mathbf{u}^r = \mathbf{x}^{r} - \mathbf{X}^{r}$
\end{enumerate}
\subsection{Formulation of the contact constraints} \label{sec_enr:OrVol}
	\begin{figure}
	\centering
	\includegraphics[clip]{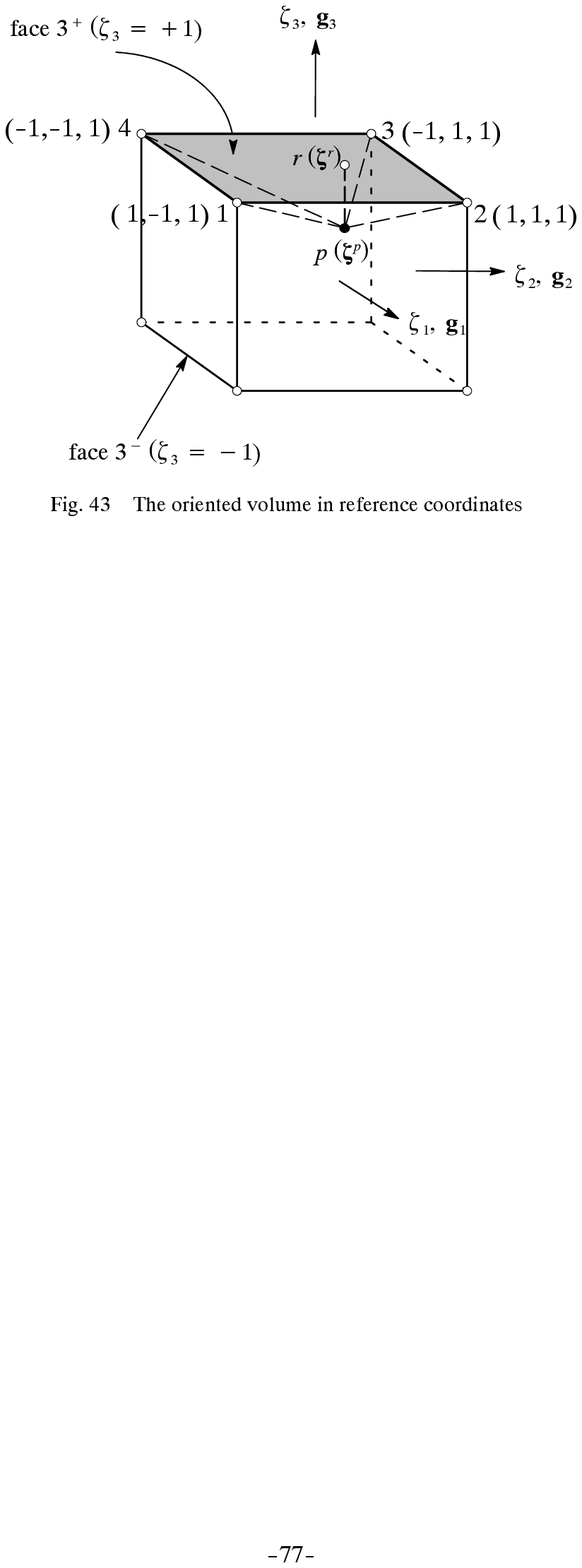}
	\caption{The oriented volume in the reference coordinates of the enriched element \label{fig:enr_OrVol}}
	\end{figure}
After the enrichment procedure described above has been applied, we implement the non-smooth contact constraints proposed in Chapter 2 to enforce geometric compatibility and prevent overlap along the interface. Since the enrichment provides additional degrees of freedom at contact locations, we can safely use a two-pass approach that enforces the interpenetrability constraints at all active interface nodes without inducing locking. The formulation of the constraints is similar to that described in Section \ref{sec_con:OrVolCons}, where we compute the oriented volume between a node $p$ and the contact surface $\zeta_{i}=c$  of an element as follows (see Figure \ref{fig:enr_OrVol})
		\begin{equation} 
		v_{i}^{p}=\zeta_{i}^{p}-c
		\end{equation}
The coordinates $\mathbf{\zeta}^{p}$ can be computed by mapping the node $p$ to the enriched element parent domain, which involves solving the system of nonlinear equations
		\begin{equation} 
		\mathbf{x}^{p}=\tilde{N}^{\alpha}(\mathbf{\zeta}^{p},\mathbf{\zeta}^r)\mathbf{x}^{\alpha} + \tilde{N}^{r}(\mathbf{\zeta}^{p},\mathbf{\zeta}^r) \mathbf{x}^r
		\end{equation}
for $\mathbf{\zeta}^{p}$. Once these coordinates have been found, the decision as to whether the node is inside or outside the element can be judged easily, since for the node to be inside the element, the following condition must be satisfied
		\begin{equation} 
		-1 \leq \zeta_{i}^{p} \leq +1 \textrm{      for     } i=1,2,3
		\end{equation}
Therefore the contact constraints take the form
		\begin{equation} 
		g_{i}^{c}=\zeta_{i}^{p}-c \geq 0 \textrm{      for     } i=1,2,3
		\end{equation}
Linearizing these constraints with respect to the global displacement variable $\mathbf{d}$, leads to (see Appendix A)
		\begin{equation} 		\nabla_{\mathbf{d}}g_{i}^{c}(\mathbf{d})=\left[\mathbf{D}_{p}^{T}-\tilde{N}^{\alpha}(\mathbf{\zeta}^p)\,\mathbf{D}_{\alpha}^{T}-\tilde{N}^{r}(\mathbf{\zeta}^p)\,\mathbf{D}_{r}^{T}\right]\tilde{\mathbf{J}}_{p}^{-T}
		\mathbf{e}_{j}
		\end{equation}
where where $\mathbf{D}_{p}$ is a $3\times3N$ Boolean matrix with value $\mathbf{I}$ at node $p$ and $\mathbf{0}$ otherwise, $\tilde{\mathbf{J}}_{p}$ is the isoparametric transformation Jacobian of the enriched element and $\mathbf{e}_{j}$ is the unit vector associated with the dimension in space along which interpenetration has occurred. As the node $p$ approached the element surface as $\mathbf{\zeta}^p \rightarrow \mathbf{\zeta}^r$ and $\mathbf{D}_{\alpha} \rightarrow 0$ while $\mathbf{D}_{r} \rightarrow 1$. Therefore the contact constraints become node-to-node and result in a set of discrete contact forces applied directly to the enrichment nodes.
%
%
\section{Stabilized interface formulation} \label{sec_enr:StabFormulation}
The discrete contact constraints described above enforce geometric compatibility at the interface nodes. This, however, does not necessarily imply full conformity of the displacement field along the interface, as the displacement remains discontinuous between the nodes. As discussed in Chapter 4, the discontinuity of the displacement (or more precisely the variational displacement) field leads to an incomplete transfer of the traction field across the interface. 
	\begin{figure}
	\centering
	\includegraphics[clip]{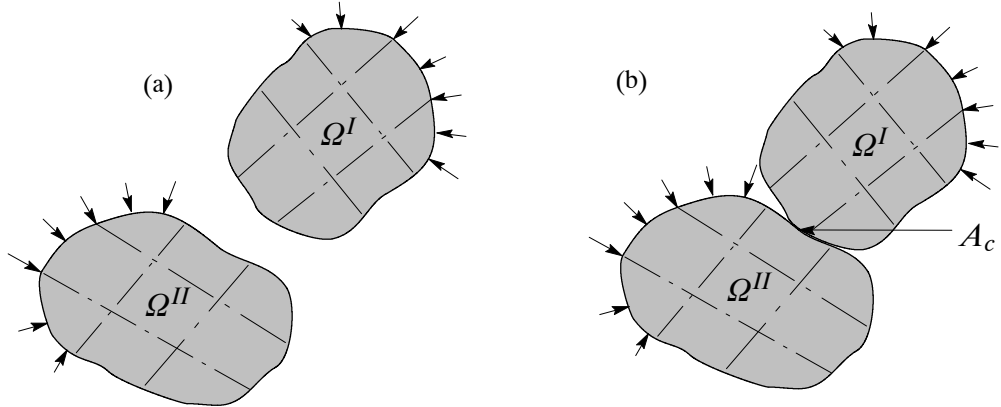}
	\caption{Two solid domains in a (a) no-contact and (b) contact configuration \label{fig:enr_BIBII}}
	\end{figure}
To remedy this issue, we apply the stabilized DG formulation discussed in Chapter 4 to the contact interface. To customize the formulation to the contact problem we note the following:
\begin{enumerate}
\item The continuity of the displacement field applies only in the direction normal to the contact interface, as the contacting bodies should be allowed to slide relative to each other. Therefore only the normal component of the traction field is of concern. 
\item Given the unilateral nature of the contact problem, enforcing the nodal contact constraints is most conveniently done via Lagrange multipliers and an active set approach, as discussed in Chapter 2. 
\end{enumerate}

Consider the domain $\Omega = \Omega^I \cup \Omega^{II}$ defined by the two bodies $\Omega^I$ and $\Omega^{II}$, as shown in Figure \ref{fig:enr_BIBII}. Let $\Gamma_t$ and $\Gamma_u$ be the Neumann and Dirichlet parts of the boundary, respectively. We define the discretization of the domain into a set of finite elements as $\Omega = \cup_e \Omega_e$. Given the body force vector $\mathbf{b}$ and the applied surface tractions $\mathbf{h}$, and if the two bodies are not in contact (Figure \ref{fig:enr_BIBII} (a)), the (discretized) energy functional in the domain is given by
\begin{equation}
\Pi = \sum_e \int_{\Omega_e} \left[ \psi (\mathbf{u}^h) - \mathbf{b}\cdot\mathbf{u}^h\right] d\Omega - \int_{{\Gamma_e}_t} \mathbf{h}\cdot\mathbf{u}^h \, d\Gamma
\end{equation}
where $\mathbf{u}^h = N^\alpha \, \mathbf{d}^\alpha$ is the discretized displacement field in each element $\Omega_e$. In the presence of a set of active contact constraints $g_i^c(\mathbf{d}) \in N_c$ that connect the two bodies to each other (Figure \ref{fig:enr_BIBII} (b)), the modified energy functional becomes
\begin{equation}
\hat{\Pi} = \sum_e \int_{\Omega_e} \left[ \psi (\mathbf{u}^h) - \mathbf{b}\cdot\mathbf{u}^h\right] d\Omega - \sum_e \int_{{\Gamma_e}_t} \mathbf{h}\cdot\mathbf{u}^h \, d\Gamma - \sum_{i\in N_c} \lambda_i^c \, g_i^c(\mathbf{d})
\end{equation}
with the conditions
\begin{equation}
g_i^c(\mathbf{d}) = 0 \ \ \ \textrm{for} \ \ \ i \in N_c.
\end{equation}
In this equation, $\mathbf{d}=\left[\mathbf{d}^1,\mathbf{d}^2,\cdots,\mathbf{d}^N\right]\in\mathbb{R}^{3n_{total}+1}$ is the global displacement vector in the body, as defined previously, with $n_{total}$ being the total number of nodes in the finite element mesh. Extremizing the above-functional yields the necessary conditions for equilibrium
\begin{align} \notag
\frac{d\hat{\Pi}}{dt} = 0 &\Rightarrow \sum_e\int_{V_e}\mathbf{P}\cdot\nabla_{\mathbf{X}}\dot{\mathbf{u}}^{h}\,dV-\sum_e\int_{V_e}\mathbf{b}_0\cdot\dot{\mathbf{u}}^{h}\,d \Omega -
\sum_e\int_{{A_e}_t}\mathbf{h}_0\cdot\dot{\mathbf{u}}^{h}\,d\Gamma  \\  \label{eq_enr:unstab}
&-\sum_{i\in N_c} \lambda_i^c \, \nabla_{\mathbf{d}}\,g_i^c(\mathbf{d})\cdot\dot{\mathbf{d}} = 0 \\
\textrm{and} \ \ \ &g_i^c(\mathbf{d}) = 0 \ \ \ \textrm{for} \ \ \ i \in N_c.
\end{align}
In these equations, $\mathbf{b}_0 = J\mathbf{b} $ is the body force vector acting on the body in its undeformed configuration, $\mathbf{t}_0=\mathbf{t}\,dA/da$ is the traction vector acting on the surface of the body in its undeformed configuration and $\mathbf{P}$ is the first Piola-Kirchhoff stress tensor. 

As discussed in Section \ref{sec_DG:mpc}, unless the enforcement of the discrete contact constraints leads to a point-wise continuity of the (normal) displacement field across the interface, equation \eqref{eq_enr:unstab} does not guarantee a complete transfer of the contact pressure across the interface. The stabilized form of this equation can be obtained by adding the normal components of the interface terms in equation \eqref{eq_DG:Gequil_stabilized_PF}, as follows
\begin{align} \notag
&\sum_e\int_{V_e}\mathbf{P}\cdot\nabla_{\mathbf{X}}\dot{\mathbf{u}}^{h}\,dV-\sum_e\int_{V_e}\mathbf{b}_0\cdot\dot{\mathbf{u}}^{h}\,d \Omega -
\sum_e\int_{{\Gamma_e}_t}\mathbf{h}_0\cdot\dot{\mathbf{u}}^{h}\,d\Gamma  \\  \notag
&-\frac{1}{2}\sum_c\int_{A_{c}^+}(\mathbf{n}^+ \cdot \mathbf{PN}^+)(\mathbf{n}^+\cdot\left[\dot{\mathbf{u}}^{h+}-\dot{\mathbf{u}}^{h-}\right])\,dA \\  &-\frac{1}{2}\sum_c\int_{A_{c}^-}(\mathbf{n}^-\cdot\mathbf{PN}^-)(\mathbf{n}^- \cdot\left[\dot{\mathbf{u}}^{h-}-\dot{\mathbf{u}}^{h+}\right])\,dA -\sum_{i\in N_c} \lambda_i^c \, \nabla_{\mathbf{d}}\,g_i^c(\mathbf{d})\cdot\dot{\mathbf{d}} =0 
\end{align}
where $\mathbf{n}$ and $\mathbf{N}$ are the normals to the deformed and undeformed configurations, respectively, and $A_{c}$ is the set of contact element interfaces. It is interesting to point out here that, even though the interface terms contain the surface normals, these terms are defined over a set of element surfaces over which the normals are uniquely defined. The contact at the nodes, which could potentially happen at non-smooth locations, is handled by the non-smooth contact formulation proposed in Chapter 2.
\begin{remark}
The presence of the Lagrange multipliers is mainly to enable a unilateral coupling approach where the sign of the multiplier is used to determine which contact constraints become inactive and need to be removed from the active set. At contact resolution, the spatial coordinates of the contact node $\mathbf{x}^p$ and that of the enrichment $\mathbf{x}^r$ are identical. Therefore, the enrichment variable $\mathbf{x}^r$ can be eliminated and the element formulation can be written in terms of $\mathbf{x}^p$ without the need for the Lagrange multiplier. However, when the constraint gets deactivated, these two variables become independent.
\end{remark}

%
%
%

\section{Numerical results} \label{sec_enr:Results}
%
%
\subsection{Sliding patch test} \label{sec_enr:ExSlidingPatch} 
\begin{figure}
\includegraphics[clip]{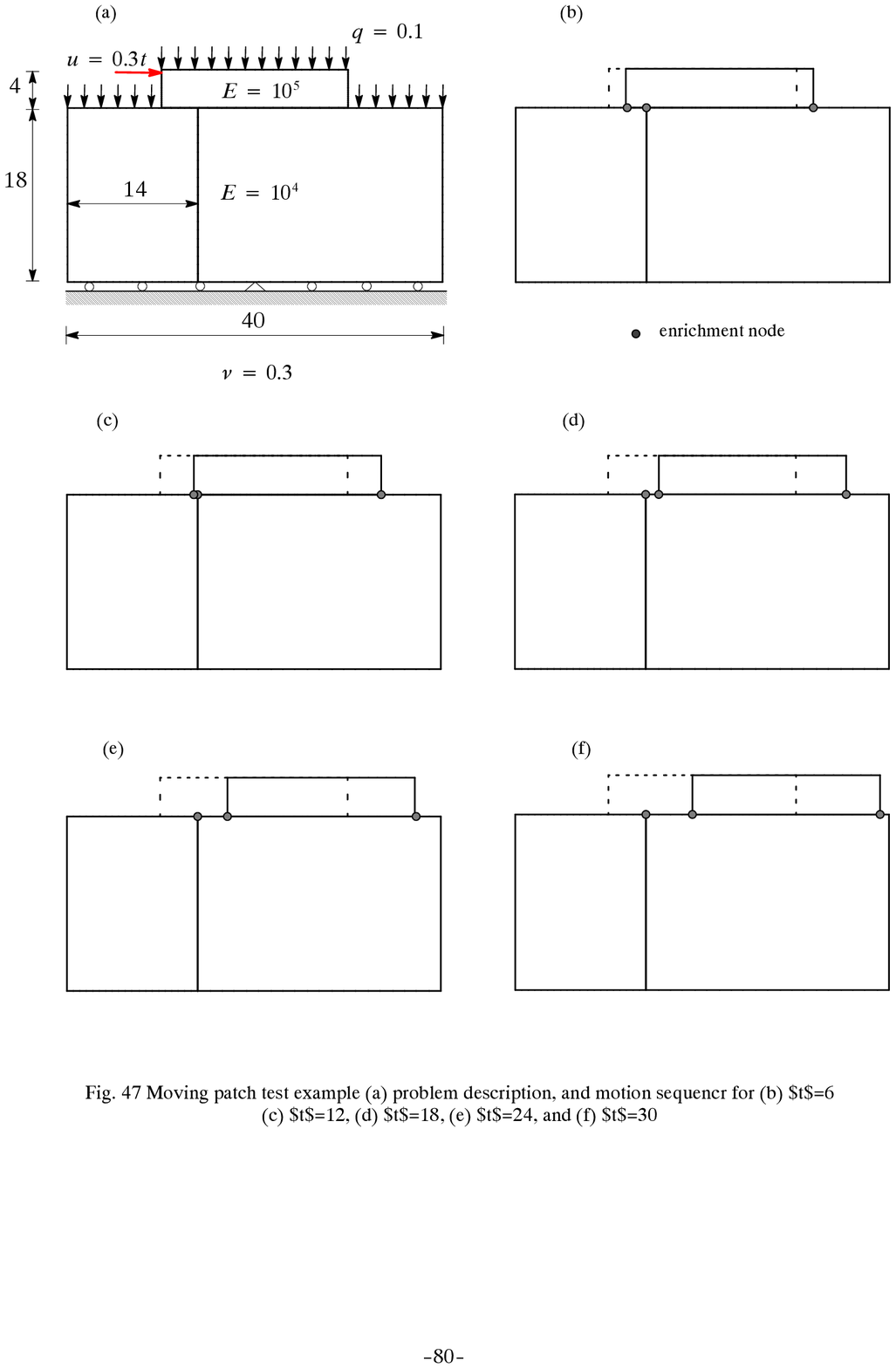}
\caption{Moving patch test example (a) problem description, and motion sequence for (b) $t$=6 (c) $t$=12, (d) $t$=18, (e) $t$=24, and (f) $t$=30 \label{fig:enr_sliding_patch}}
\end{figure}

%
Figure \ref{fig:enr_sliding_patch} (a) shows a modified patch test where a horizontal displacement $u_1 = 0.3 t$ is applied to the punch that causes it to slide along the foundation surface. The parameter $t$ is a pseudo-time variable we use to increment the applied motion. The purpose of this example is to test the performance of the moving enrichment formulation. Furthermore, we specify the material moduli for the punch to be $E=10^5$ and $\nu=0.3$, whereas $E=10^4$ and $\nu=0.3$ in the foundation, making the punch substantially stiffer that the foundation.

Figures \ref{fig:enr_sliding_patch} (b) through (f) show the motion sequence with increasing $t$. It is obvious that, even in the presence of different material properties, the DG formulation ensures a complete transfer of forces across the interface. Furthermore, the moving enrichment is successful in tracking the location of the contact nodes, even across the element subdivision in the foundation, and the contact constraints remain node-to-node throughout the motion. It is useful to point out here that if the enrichment location coincides exactly with an existing node, the Jacobian of the enriched element goes to zero. To eliminate this issue, we keep a tolerance of $|\Delta \mathbf{\zeta}| = 10^{-2}$ between the enrichment location and any existing node. 

%
\subsection{Beam bending} \label{sec_enr:ExBeamBending}
	\begin{figure}
	\centering
	\includegraphics[clip]{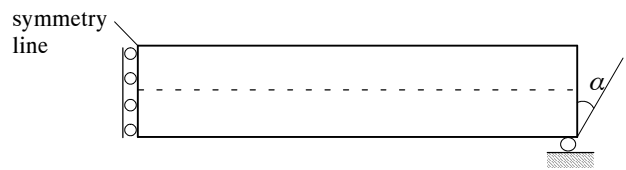}
	\caption{Double beam bending problem: problem definition  \label{fig:enr_bending}}
	\end{figure}
We apply the stabilized interface formulation to the beam bending problem discussed in Chapter 3 \cite{pusoCMAME04}.
The problem consists of the two 10 x 1 beams shown in Figure \ref{fig:enr_bending}. The material properties are $E =1$, $\nu =0$. The top and bottom surface of the structure are subjected to a pressure of $p=0.1$. A rotation $\alpha$ is then applied at the right end of the beams. The rotation angle is incremented quasi-statically from 0 to 90 degrees.
	\begin{figure}
	\centering
	\includegraphics[clip]{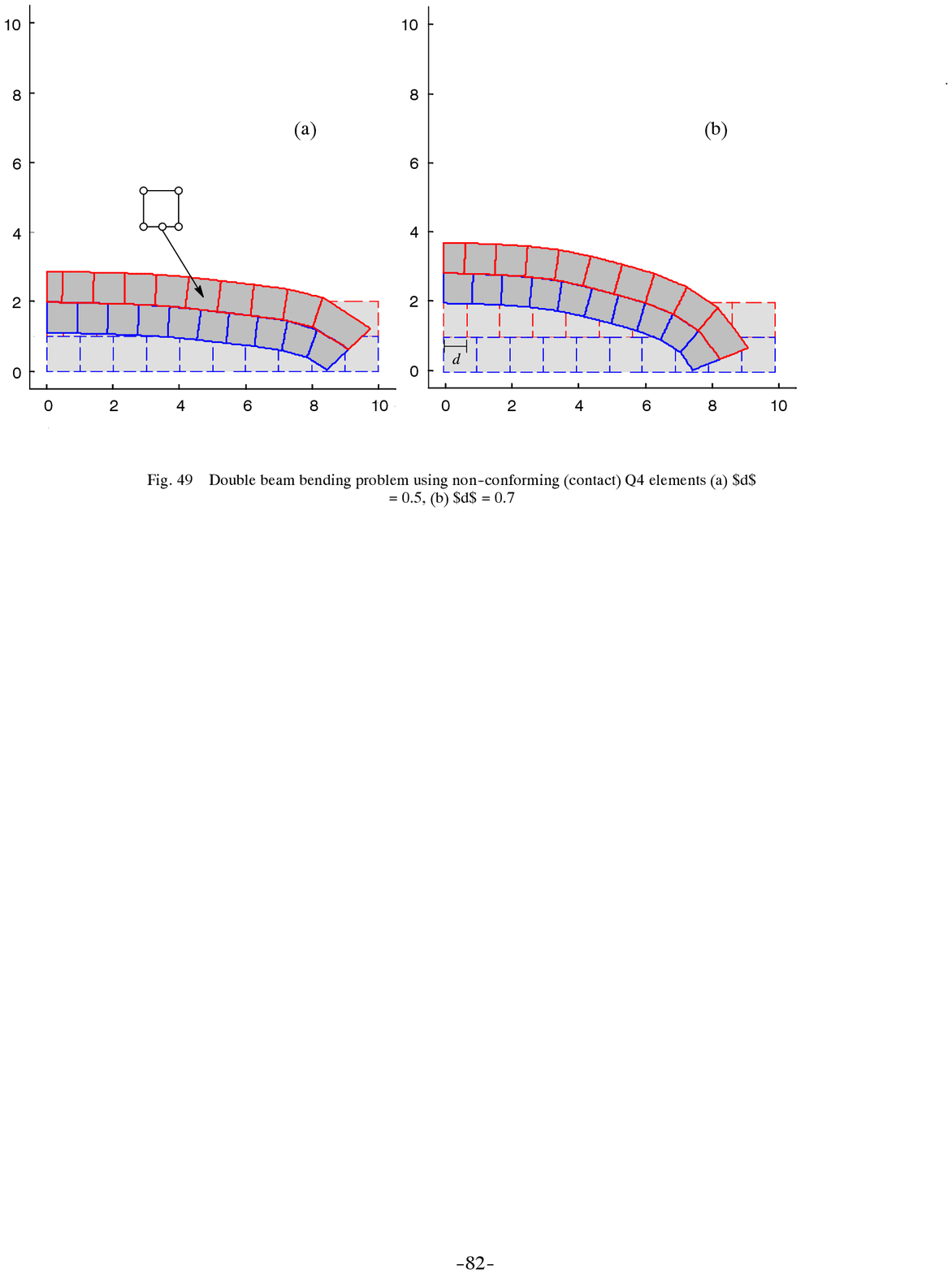}
	\caption{Double beam bending problem using non-conforming (contact) Q4 elements with no sliding (a) $d$ = 0.5, (b) $d$ = 0.7  \label{fig:enr_lock}}
	\end{figure}

\begin{figure}
	\centering
	\includegraphics[clip]{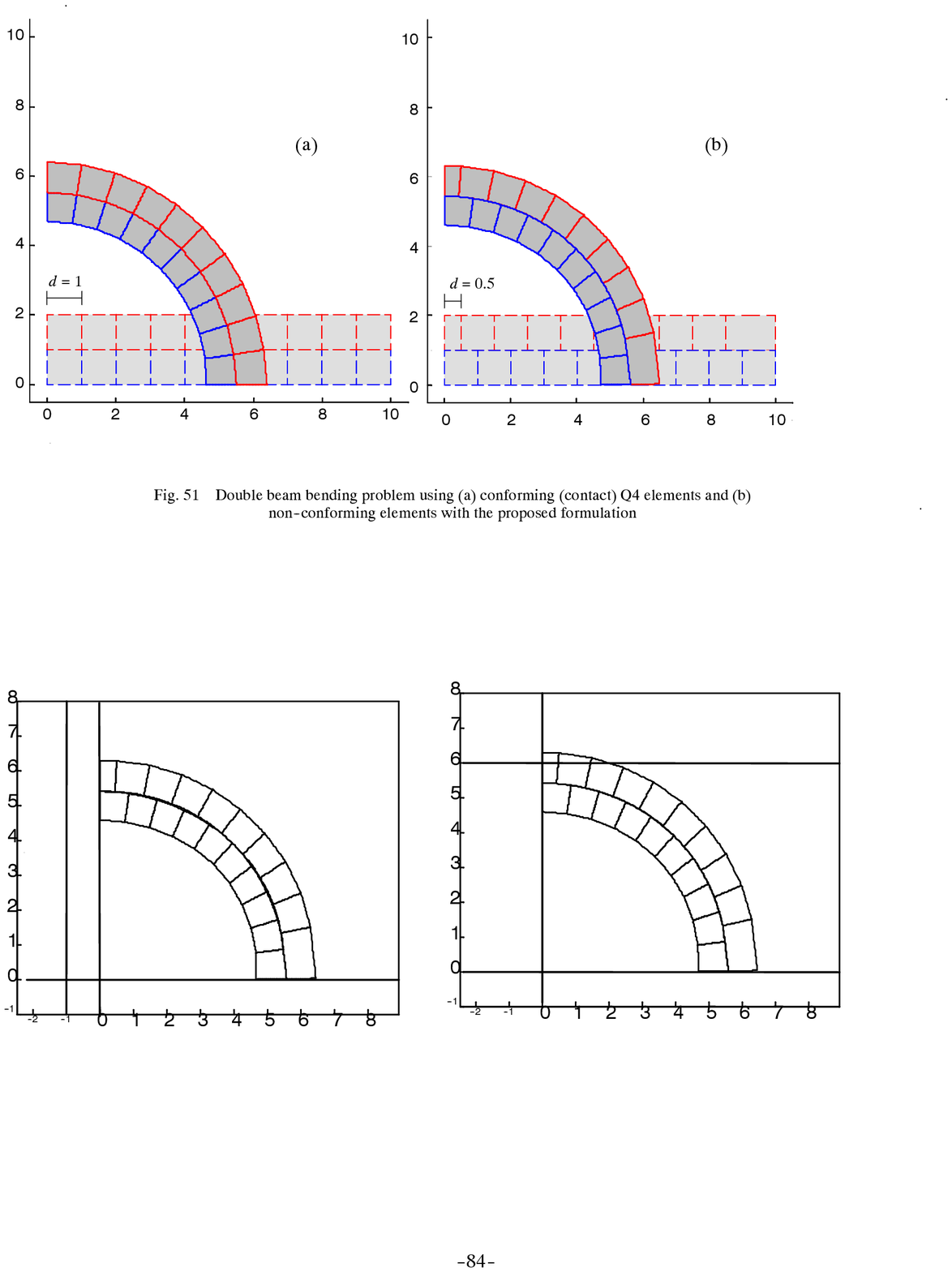}
	\caption{Double beam bending problem using (a) conforming (contact) Q4 elements and (b) non-conforming elements with the proposed formulation  \label{fig:enr_bending_sol}}
	\end{figure}

Figure \ref{fig:enr_lock} shows the configuration we obtained previously using the two-pass node-to-surface formulation, which exhibits severe locking. Figure \ref{fig:enr_bending_sol} (b) shows the result obtained using the proposed enrichment and stabilized interface formulation for the case $d=0.5$\footnote{We added tying nodes to top surface of the second and ninth elements of the lower beam to further restrict sliding. The two-pass node-to-surface formulation still locks in this configuration.}.  It is clear from this figure that the locking is eliminated and the result matches that of the conforming mesh, as shown in Figure \ref{fig:enr_bending_sol} (a).

\begin{figure}
	\centering
	\includegraphics[clip]{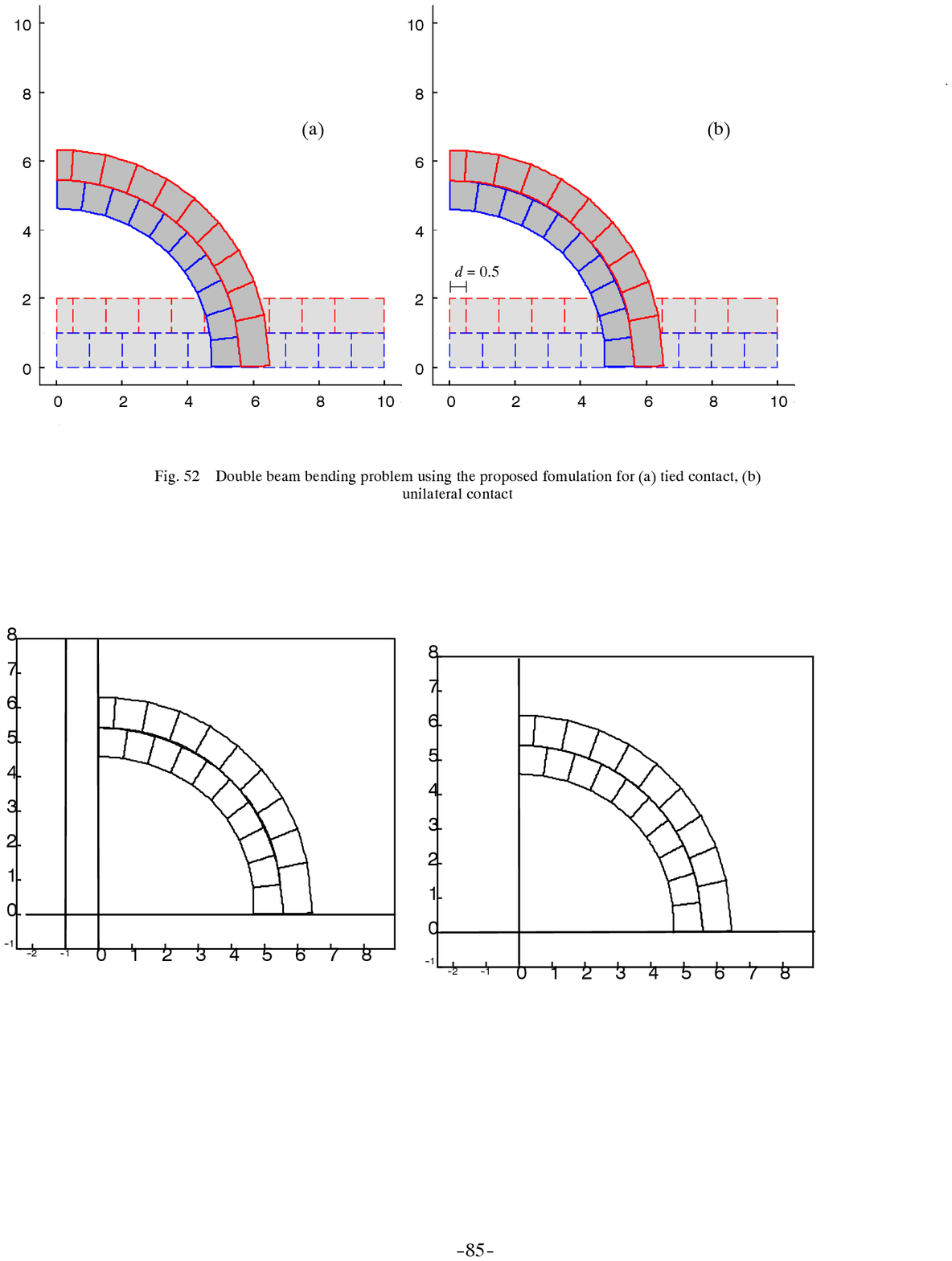}
	\caption{Double beam bending problem using the proposed formulation for (a) tied contact, (b) unilateral contact  \label{fig:enr_bending_tied}}
	\end{figure}

It is interesting to point out here that, even though all contact events get activated initially, many constraints are released before the final configuration is reached (the stabilization terms are applied nonetheless). This feature is not present in the solution provided by Puso and Laursen \cite{pusoCMAME04}, where the contact was assumed to be tied and all contact events were kept active throughout the motion. To investigate whether the lack of locking is due to the release of the contact constraints, and not to the interface formulation, we repeat the solution assuming tied contact events and keep all constraints active regardless of the sign of the Lagrange multiplier. The result is shown in Figure \ref{fig:enr_bending_tied} (a) and clearly matches that obtained earlier with unilateral constraints, shown in Figure \ref{fig:enr_bending_tied} (b). Therefore we conclude that the elimination of the interface locking is due to the enrichment and stabilization of the interface.

\section{Conclusions} \label{sec_enr:Concl}

Surface locking is a phenomenon that hampers the two-pass node-to-surface contact formulation. To address this issue, we propose a stabilization procedure that makes use of a local enrichment of the contact surface to strongly enforce the non-penetration constraint at contact locations without inducing surface locking. This kind of formulation is much needed since the methods available in the literature are either biased single-pass methods  that enforce the interpenetration condition only in a weak sense, or two-pass methods that employ ad-hoc procedures for the selection of active contact constraints. 









%% file: appendix.tex
%
%
%
%
\setcounter{chapter}{0}
\def\thechapter{\Alph{chapter}}
\setcounter{section}{0} \setcounter{figure}{0}
\def\thesection{\thechapter.\arabic{section}}
%
%
\chapter{Linearization of the contact constraint} \label{AppA:ConsLin}
%
%
\section{Calculation of the Jacobian of the contact constraint function} \label{AppSec:Jac}
Suppose that contact was detected between node $p$ and element $k$ in the direction $j$. The rate of change of the contact constraint with respect to $\mathbf{d}$ is
\begin{equation}
\nabla_{\mathbf{d}}g^c = \frac{\partial g^c}{\partial d_i}\mathbf{e}_i
\end{equation}
Since $g$ is a function of $\mathbf{\zeta}^p$, which in turn is an implicit function of the current nodal coordinates, the partial derivatives $\partial g_c / \partial d_i$ are obtained via the chain rule. To wit,
\begin{align}   \notag
\nabla_{\mathbf{d}}g^c&=\frac{\partial g^c}{\partial\zeta_1}\nabla_{\mathbf{d}}\zeta_1 + \frac{\partial g^c}{\partial\zeta_2}\nabla_{\mathbf{d}}\zeta_2 +\frac{\partial g^c}{\partial\zeta_3}\nabla_{\mathbf{d}}\zeta_3     \\
\label{eq:graddgc}
 &= \left[ 
\begin{array}{ccc}
  \nabla_{\mathbf{d}}\zeta_1^p & \nabla_{\mathbf{d}}\zeta_2^p & \nabla_{\mathbf{d}}\zeta_3^p	
\end{array}
 \right] \left[  
\begin{array}{c}
	\delta_{1j} \\ \delta_{2j} \\ \delta_{3j}
\end{array}
 \right]
\end{align}
where $\delta_{ij} = 1$ for $i = j$ (no summation implied) and $\delta_{ij} = 0$ otherwise. Eq. \eqref{eq:graddgc} can be alternatively expressed as
\begin{equation}  \label{eq:delgcT}
\nabla_{\mathbf{d}}g^c=\nabla_{\mathbf{d}}{\mathbf{\zeta}^p}^T \mathbf{e}_j\Rightarrow \nabla_{\mathbf{d}}{g^c}^T=\mathbf{e}_j^T\nabla_{\mathbf{d}}\mathbf{\zeta}^p
\end{equation}
where $\mathbf{e}_j$ is the unit vector associated with the dimension in space along which interpenetration has occurred. 

The term $\nabla_{\mathbf{d}}\mathbf{\zeta}^p$ in \eqref{eq:delgcT} can be computed as follows. Recall that, from the isoparametric formulation of the element at point $p$,
\begin{equation}
\mathbf{x}^p = N^\alpha\left(\mathbf{\zeta}^p\right)\mathbf{x}^\alpha
\end{equation}
The gradient of Eq. \eqref{eq:E} with respect to $\mathbf{d}$ yields
\begin{equation}
\nabla_\mathbf{d}\mathbf{x}^p =\mathbf{x}^\alpha \otimes\nabla_{\mathbf{d}}N^\alpha\left(\mathbf{\zeta}^p\right) + N^\alpha\left(\mathbf{\zeta}^p\right)\nabla_{\mathbf{d}}\mathbf{x}^\alpha
\end{equation}
\begin{equation} \label{eq:Pp}
\Rightarrow\mathbf{x}^\alpha \otimes\nabla_{\mathbf{d}}N^\alpha\left(\mathbf{\zeta}^p\right) 
= \mathbf{D}_p- N^\alpha\left(\mathbf{\zeta}^p\right)\mathbf{D}_\alpha
\end{equation}
where $\mathbf{D}_q=\nabla_{\mathbf{d}}\mathbf{x}^q$ is a Boolean matrix with value $\mathbf{I}$ at the node $q$ and $\mathbf{0}$ otherwise. The left-hand side of Eq. \eqref{eq:Pp} can be computed using the chain rule
\begin{equation} \label{eq:Jpdeluzp}
\mathbf{x}^\alpha\otimes\nabla_{\mathbf{d}}N^\alpha\left(\mathbf{\zeta}^p\right)=\mathbf{x}^\alpha\otimes\frac{\partial N^\alpha}{\partial\mathbf{\zeta}^p}\frac{\partial\mathbf{\zeta}^p}{\partial\mathbf{d}} = \mathbf{J}_p\nabla_{\mathbf{d}}\mathbf{\zeta}^p
\end{equation}
where $\mathbf{J}_p=\partial\mathbf{x}^p/\partial\mathbf{\zeta}^p=\nabla_{\mathbf{\zeta}}\mathbf{x}^p$ is the isoparametric transformation Jacobian of the element, evaluated at node $p$. Combining Eqs. \eqref{eq:Pp} and \eqref{eq:Jpdeluzp}
\begin{equation} \label{eq:graddzp}
\mathbf{D}_p-N^\alpha\mathbf{D}_\alpha=\mathbf{J}_p\nabla_{\mathbf{d}}\mathbf{\zeta}^p
\Rightarrow\nabla_{\mathbf{d}}\mathbf{\zeta}^p = \mathbf{J}_p^{-1}\left[\mathbf{D}_p-N^\alpha\mathbf{D}_\alpha\right]
\end{equation}
Consequently, Eq. \eqref{eq:delgcT} reduces to
\begin{equation}
\nabla_{\mathbf{d}}{g^c}^T = \mathbf{e}_j^T\mathbf{J}_p^{-1}\left[\mathbf{D}_p-N^\alpha\mathbf{D}_\alpha\right]
\end{equation}
Therefore
\begin{equation}
\nabla_{\mathbf{d}}g^c = \left[\mathbf{D}_p^T-N^\alpha\mathbf{D}_\alpha^T\right]\mathbf{J}_p^{-T}\mathbf{e}_j
\end{equation}
More generally, the gradient of a constraint at a given node can be expressed as follows:
\begin{equation}
\nabla_{\mathbf{d}}g^c = \left[\mathbf{D}_p^T-N^\alpha\mathbf{D}_\alpha^T\right]\mathbf{f}_p^c
\end{equation}
The term $\left[\mathbf{D}_p^T-N^\alpha\mathbf{D}_\alpha^T\right]$ actually plays the role of applying a discrete contact force vector $\mathbf{f}_p^c=\mathbf{J}_p^{-T}\mathbf{e}_j$ at node $p$ and distributing an equal reaction among the nodes of the contact element. When contact is resolved, the reaction applies to the nodes of the contact surface only, since the value of $N^\alpha$ is zero at other locations. Therefore, conservation of linear momentum is maintained.
%
%
\section{Calculation of the Hessian of the contact constraint function} \label{AppSec:Hes}
The Hessian of the contact constraints can be calculated by taking the directional derivative of the Jacobian in the direction of a displacement increment $\Delta\mathbf{d}$
\begin{align}  \notag
\mathbf{H}^c\Delta\mathbf{d} &\equiv D\left(\nabla_{\mathbf{d}}g_c\right)\cdot\Delta\mathbf{d} \\
\label{eq:HcDeltad1}
&= -\mathbf{D}_\alpha^T\left[\mathbf{J}_p^{-T}\mathbf{e}_j \otimes D N_\alpha\cdot\Delta\mathbf{d}\right]+\left[\mathbf{D}_p^T-N^\alpha\mathbf{D}_\alpha^T\right]\left[D\mathbf{J}_p^{-T}\cdot\Delta\mathbf{d}\right]\mathbf{e}_j
\end{align}
To economize the notation, let us define $\mathbf{P}_p\equiv\left[\mathbf{D}_p-N^\alpha\mathbf{D}_\alpha\right]$. Applying the chain rule to the first term in Eq. \eqref{eq:HcDeltad1}, we get
\begin{equation} \label{eq:HcDeltad2}
\Rightarrow\mathbf{H}^c\Delta\mathbf{d}=-\mathbf{D}_\alpha^T\left[\mathbf{J}_p^{-T}\mathbf{e}_j\otimes\nabla_{\mathbf{\zeta}}N^\alpha\nabla_{\mathbf{d}}\mathbf{\zeta}^p\Delta\mathbf{d}\right] - \mathbf{P}_p\mathbf{J}_p^{-T}\left[D\mathbf{J}_p^{-T}\cdot\Delta\mathbf{d}\right]\mathbf{J}_p^{-T}\mathbf{e}_j
\end{equation}
Recall that $\mathbf{x}^p = N^{\alpha}(\mathbf{\zeta}^p)\mathbf{x}^\alpha$ and $\mathbf{J}_p=\partial\mathbf{x}^p/\partial\mathbf{\zeta}^p=\mathbf{x}^\alpha\otimes\nabla_{\mathbf{\zeta}}N^\alpha$. Thus,
\begin{align}
D\mathbf{J}_p^T\cdot\Delta\mathbf{d} &= \nabla_{\mathbf{\zeta}}N^\alpha\otimes\Delta\mathbf{d}^\alpha+\nabla_{\mathbf{d}}\left(\nabla_{\mathbf{\zeta}}N^\alpha\right)\Delta\mathbf{d}\otimes\mathbf{x}^\alpha \\
&=\nabla_{\mathbf{\zeta}}N^\alpha\otimes\Delta\mathbf{d}^\alpha+\nabla_{\mathbf{\zeta}}\left(\nabla_{\mathbf{\zeta}}N^\alpha\right)\nabla_{\mathbf{d}}\mathbf{\zeta}^p\Delta\mathbf{d}\otimes\mathbf{x}^\alpha
\end{align}
From Eq. \eqref{eq:graddzp}, $\nabla_{\mathbf{d}}\mathbf{\zeta}^p=\mathbf{J}_p^{-1}\mathbf{P}_p$. To wit,
\begin{equation} \label{eq:DJpTDeltad}
D\mathbf{J}_p^T\cdot\Delta\mathbf{d} = \nabla_{\mathbf{\zeta}}N^\alpha\otimes\Delta\mathbf{d}^\alpha+\nabla_{\mathbf{\zeta}}\left(\nabla_{\mathbf{\zeta}}N^\alpha\right)\mathbf{J}_p^{-1}\mathbf{P}_p\Delta\mathbf{d}\otimes\mathbf{x}^\alpha
\end{equation}
Substituting Eq. \eqref{eq:DJpTDeltad} into Eq. \eqref{eq:HcDeltad2} leads to
\begin{align}  \notag
\mathbf{H}^c\Delta\mathbf{d} = &-\mathbf{D}_\alpha^{T}\left[\mathbf{J}_p^{-T}\mathbf{e}_j \otimes \nabla_{\mathbf{\zeta}}N^\alpha\mathbf{J}_p^{-1}\mathbf{P}_p\right]\Delta\mathbf{d}                          \\
\label{eq:ODEv}
&- \mathbf{P}_p^T\mathbf{J}_p^{-T}\left[\nabla_{\mathbf{\zeta}}N^\alpha\otimes \Delta\mathbf{d}^\alpha 
+ \nabla_{\mathbf{\zeta}}\left(\nabla_{\mathbf{\zeta}}N^\alpha\right)\mathbf{J}_p^{-1}\mathbf{P}_p\Delta\mathbf{d}\otimes\mathbf{x}^\alpha \right]\mathbf{J}_p^{-T}\mathbf{e}_j                        \\
\notag
= &-\left\{\mathbf{D}_\alpha^{T}\left[\mathbf{J}_p^{-T}\mathbf{e}_j\otimes\nabla_{\mathbf{\zeta}}N^\alpha\right]\mathbf{J}_p^{-1}\mathbf{P}_p  
+\mathbf{P}_p^T \mathbf{J}_p^{-T}\left[\nabla_{\mathbf{\zeta}}N^\alpha\otimes\mathbf{J}_p^{-T}\mathbf{e}_j\right]\mathbf{D}_\alpha\right\}\Delta\mathbf{d}    \\
&-\left\{\mathbf{P}_p^T \mathbf{J}_p^{-T}\left[\gamma^\alpha\nabla_{\mathbf{\zeta}}\left(\nabla_{\mathbf{\zeta}}N^\alpha\right)\right]\mathbf{J}_p^{-1} \mathbf{P}_p\right\}\Delta\mathbf{d}                         \\
= &- \left\{ \left[ \mathbf{T}^c + {\mathbf{T}^c}^T \right]+ \mathbf{P}_p^T \mathbf{J}_p^{-T}\left[\gamma^\alpha\nabla_{\mathbf{\zeta}}\left(\nabla_{\mathbf{\zeta}}N^\alpha\right)\right]\mathbf{J}_p^{-1} \mathbf{P}_p\right\}\Delta\mathbf{d}
\end{align}
in which $\mathbf{T}^c\equiv\mathbf{D}_\alpha^T\left[\mathbf{J}_p^{-T}\mathbf{e}_j\otimes\nabla_{\mathbf{\zeta}}N^\alpha\right]\mathbf{J}_p^{-1}\mathbf{P}_p$ and $\gamma^\alpha\equiv\left(\mathbf{x}^\alpha\cdot\mathbf{J}_p^{-T}\mathbf{e}_j\right)$. The Hessian is clearly symmetric.